\documentclass[12pt]{article}
\usepackage{amsmath,amsthm,amssymb,amsfonts}

\newtheorem{theorem}{Theorem}[section]
\newtheorem{proposition}[theorem]{Proposition}
\newtheorem{corollary}[theorem]{Corollary}
\newtheorem{lemma}[theorem]{Lemma}

\newtheorem{remark}[theorem]{Remark}
\newtheorem{example}[theorem]{Example}

\begin{document}

\title{On Genus Two Riemann Surfaces Formed from Sewn Tori}
\author{Geoffrey Mason\thanks{%
Support provided by the National Science Foundation DMS-0245225, and the
Committee on Research at the University of California, Santa Cruz } \\
Department of Mathematics, \\
University of California Santa Cruz, \\
CA 95064, U.S.A. \and Michael P. Tuite\thanks{%
Supported by The Millenium Fund, National University of Ireland, Galway } \\
Department of Mathematical Physics, \\
National University of Ireland, \\
Galway, Ireland.}
\maketitle

\begin{abstract}
We describe the period\ matrix and other data on a higher genus Riemann
surface in terms of data coming from lower genus surfaces via an explicit
sewing procedure. We consider in detail the construction of a genus two
Riemann surface by either sewing two punctured tori together or by sewing a
twice-punctured torus to itself. In each case the genus two period matrix is
explicitly described\ as a holomorphic map from a suitable domain
(parameterized by genus one moduli and sewing parameters) to the Siegel
upper half plane\ $\mathbb{H}_{2}$. Equivariance of these maps under certain
subgroups of $Sp(4,\mathbb{Z)}$ is shown. The invertibility of both maps in
a particular domain of $\mathbb{H}_{2}$ is also shown.
\end{abstract}

\section{Introduction}

\label{Sect_Intro}

This paper is the second of a series intended to develop a mathematically
rigorous theory of chiral partition and $n$-point functions on Riemann
surfaces at all genera, based on the theory of vertex operator algebras. The
purpose of the paper is to provide a rigorous and explicit description of
the period \ matrix and other data on a higher genus Riemann surface in
terms of data coming from lower genus surfaces via an explicit sewing
procedure. In particular, we consider and compare in some detail the
construction of a genus two Riemann surface in two separate ways: either by
sewing two tori together or by sewing a torus to itself. Although our
primary motivation is to lay the foundations for the explicit construction
of the partition and $n$-point functions for a vertex operator algebra on
higher genus Riemann surfaces \cite{MT2}, we envisage that the results
herein may also be of wider interest.

As is well-known (cf. \cite{H} for a systematic development), the axioms for
a vertex (operator) algebra $V$ amount to an algebraicization of aspects of
the theory of gluing spheres, i.e. compact Riemann surfaces of genus zero. A
quite complete theory of (bosonic) $n$-point functions at genus one was
developed by Zhu in his well-known paper \cite{Z}. For particularly
well-behaved vertex operator algebras (we have in mind \emph{rational}
vertex operator algebras in the sense of \cite{DLM}), Zhu essentially showed
that the $n$-point functions, defined as particular graded traces over $V$,
are elliptic and have certain $SL(2,\mathbb{Z})$ modular-invariance
properties\footnote{%
The general conjecture that the partition function is a modular function if $%
V$ is rational remains open.} with respect to the torus modular parameter $%
\tau $. A complete description of bosonic Heisenberg and lattice VOA $n$%
-point functions is given in \cite{MT1}. Apart from its intrinsic interest,
modular-invariance is an important feature of conformal field theory and its
application to string theory e.g. \cite{P}, \cite{GSW}. Mathematically, it
has been valuable in recent developments concerning the structure theory of
vertex operator algebras \cite{DM1}, \cite{DM2}, and it may well play an
important role in geometric applications such as elliptic cohomology and
elliptic genus. These are all good reasons to anticipate that a theory of $n$%
-point functions at higher genus will be valuable; another is the
consistency of string theory at higher loops (genera). Our ultimate goal,
then, is to emulate Zhu's theory by first defining and then understanding
the automorphic properties of $n$-point functions on a higher genus Riemann
surface where the modular variable $\tau $ is replaced by the period matrix $%
\Omega^{(g)}$, a point in the Siegel upper half-space $\mathbb{H}_{g}$ at
genus $g$.

When one tries to implement the vision outlined above at genus $g\geq 2$,
difficulties immediately arise which have no analog at lower genera \cite{T}. 
In Zhu's theory, there is a clear relation between the variable $\tau $
and the relevant vertex operators: it is not the definition of $n$-point
functions, but rather the elucidation of their properties, that is
difficult. At higher genus, the very definition of $n$-point function is
less straightforward and raises interesting issues. Many (but not all) of
these are already present at genus $2$, and it is this case that we mainly
deal with in the present paper.

A very general approach to conformal field theory on higher genus Riemann
surfaces has been discussed in the physics literature \cite{MS}, 
\cite{So},\cite{D'HP},\cite{VV},\cite{P}. In particular, there has been 
much progress in recent years in understanding genus two superstring 
theory \cite{D'HPI},\cite{D'HPVI},\cite{D'HGP}.
The basic idea we follow is that by cutting a Riemann surface
along various cycles it can be reduced to (thrice) punctured spheres, and
conversely one can construct Riemann surfaces by sewing punctured spheres.
It is interesting to note that Zhu's $g=1$ theory is \emph{not} formulated
by sewing punctured spheres per se, but rather by implementing a conformal
map of the complex plane onto a cylinder. Tracing over $V$ has the effect of
sewing the ends of the cylinder to obtain a torus. This idea does not
generalize to the case when $g=2$, and we must hew more closely to the
sewing approach of conformal field theory. Roughly speaking, what we do is
sew tori together in order to obtain a compact Riemann surface $S$ of genus
2 and which is endowed with certain genus 1 data encoded by $V$ via the Zhu
theory \cite{T}. There are two essentially different ways to obtain $S$
(which, for simplicity, we take to have no punctures in this work) from
genus 1 data: either by sewing a pair of once-punctured complex tori, or by
sewing a twice punctured torus to itself (attaching a handle). These two
sewing schemes will give rise to seemingly different theories and different
definitions of $n$-point functions. These issues are discussed in detail in
the sequel \cite{MT2}.

We now give a more technical introduction to the contents of the present
paper. Notwithstanding our earlier discussion of the role of vertex
operators, they do not appear explicitly in the present work! We are
concerned here exclusively with setting up foundations so that the ideas we
have been discussing are rigorous and computationally effective. Section \ref%
{Sect_ellipticfunctions} records the many modular and elliptic-type
functions that we will need. In the paper \cite{Y}, which is very important
for us, Yamada developed a general approach to computing the period matrix
of a Riemann surface $S$ obtained by sewing Riemann surfaces $S_{1},S_{2}$
(which may coincide) of smaller genus. In Sections \ref{Sect_Epsilon_g} and %
\ref{Sect_Epsilon_Torii} we develop the theory in the case that $S_{1}$ and $%
S_{2}$ are \emph{distinct}. We refer to this as the \emph{$\epsilon $%
-formalism}, $\epsilon $ itself being a complex number which is a part of
the data according to which the sewing is performed. (See Figure 1 below).
We begin in Section \ref{Sect_Epsilon_g} with some of the details of
Yamada's general theory, and make some explicit computations. In particular,
we introduce infinite matrices $A_{a}$ for $a=1,2$, whose entries are
certain weighted moments of the normalized differential of the second kind
on $S_{a}$. These matrices determine another infinite matrix $X$ whose
entries are weighted moments of the normalized differential of the second
kind on $S$ (Proposition \ref{PropXij}), and this in turn determines the
period matrix $\Omega ^{(g)}$ of $S$ (Theorem \ref{Theoremperiodgeps}). In
particular, the infinite matrix 
\begin{equation}
I-A_{1}A_{2}  \label{A1A2matrix}
\end{equation}%
plays an important role ($I$ is the infinite unit matrix). The entries of
this matrix depend on data coming from $S_{a}$, and in particular they are
power series in $\epsilon $. We show (Theorem \ref{Theorem_Det}) that (\ref%
{A1A2matrix}) has a well-defined determinant $\det (I-A_{1}A_{2})$ which is 
\emph{holomorphic} for small enough $\epsilon $. Section \ref{Sect_Epsilon_g}
ends with some additional results concerning the holomorphy of $\det
(I-A_{1}A_{2})$ and $\Omega ^{(g)}$ in various domains.

In Section \ref{Sect_Epsilon_Torii} we study in more detail the case in
which the $S_{a}$ have genus 1, so that they have a modulus $\tau _{a}\in 
\mathbb{H}_{1}$. The triple $(\tau _{1},\tau _{2},\epsilon )$ determines a
genus 2 surface as long as the three parameters in question satisfy a
certain elementary inequality. This defines a manifold $\mathcal{D}%
^{\epsilon }(\tau _{1},\tau _{2},\epsilon )\subseteq \mathbb{H}_{1}\times 
\mathbb{H}_{1}\times \mathbb{C}$ consisting of all such \emph{admissible
triples}. Associating to this data the genus two period matrix $\Omega
^{(2)}=\Omega ^{(2)}(\tau _{1},\tau _{2},\epsilon )$ of $S$ defines a map 
\begin{eqnarray}
F^{\epsilon }:\mathcal{D}^{\epsilon }(\tau _{1},\tau _{2},\epsilon )
&\rightarrow &\mathcal{H}_{2}  \label{Fdef} \\
(\tau _{1},\tau _{2},\epsilon ) &\mapsto &\Omega ^{(2)}(\tau _{1},\tau
_{2},\epsilon )  \notag
\end{eqnarray}%
which is important for everything that follows. When we introduce partition
functions in the $\epsilon $-formalism at $g=2$ in the sequel to the present
paper \cite{MT2}, they will be functions on $\mathcal{D}^{\epsilon }$, not $%
\mathbb{H}_{2}$. The map $F^{\epsilon }$ interpolates between the two
domains. We obtain (Theorem \ref{Theorem_epsperiod}) an explicit expression
for the genus 2 period matrix 
\begin{equation}
\Omega ^{(2)}(\tau _{1},\tau _{2},\epsilon )=\left( 
\begin{array}{cc}
\Omega _{11}^{(2)} & \Omega _{12}^{(2)} \\ 
\Omega _{12}^{(2)} & \Omega _{22}^{(2)}%
\end{array}%
\right)  \label{Omageepsilonmatrix}
\end{equation}%
determined by an admissible triple. Each $\Omega _{ij}^{(2)}$ turns out to
be essentially a power series in $\epsilon $ with coefficients which are%
\footnote{%
Notation for modular and elliptic-type functions is covered in Section 2.}
quasimodular forms, i.e. certain polynomials in the Eisenstein series $%
E_{2}(\tau _{i})$, $E_{4}(\tau _{i})$, $E_{6}(\tau _{i})$ for $i=1,2$.
Moreover $F^{\epsilon }$ is an analytic map, and we show (Theorem \ref%
{TheoremGequiv}) that it is equivariant with respect to the action of a
group $G\cong (SL(2,\mathbb{Z})\times SL(2,\mathbb{Z}))\wr \mathbb{Z}_{2}$
(the wreathed product of $SL(2,\mathbb{Z})$ and $\mathbb{Z}_{2}$). $G$
embeds into $Sp(4,\mathbb{Z})$ in a standard way, and this defines the
action of $G$ on $\mathbb{H}_{2}$. The action on $\mathcal{D}^{\epsilon }$
is explained in Subsection \ref{Subsec_Feps}. These calculations are
facilitated by an alternate description (Proposition \ref{Propepsperiodgraph}%
) of the entries of (\ref{Omageepsilonmatrix}) in terms of combinatorial
gadgets that we call \emph{chequered necklaces}. They are certain kinds of
graphs with nodes labelled by positive integers and edges labelled by
quasimodular forms, and they play a critically important role in the sequel
to the present paper. We show (Proposition \ref{Prop_Fepsinverse}) that
about any degeneration point $p$ where $\epsilon =0$ (i.e. the two tori $%
S_{1},S_{2}$ touch at a point), there is a $G$-invariant neighborhood of $p$
throughout which $F^{\epsilon }$ is invertible.

Sections \ref{Sect_Rho_g} and \ref{Sect_Rho_Torus} are devoted to
development of the corresponding formalism in the case that $S$ is obtained
by self-sewing (i.e. attaching a handle to) a surface $S_{1}$ of one lower
genus. We refer to this as the $\rho $-formalism. Although we are able to
achieve results that parallel the development of the $\epsilon $-formalism
outlined in the previous paragraph, it is fair to say that of the two, the $%
\rho $-formalism is the more complicated. In Section \ref{Sect_Rho_g} we
first discuss the results of Yamada (loc. cit.) in a general $\rho $%
-formalism, and calculate weighted moments as before. This leads us to
introduce the analog of (\ref{A1A2matrix}), namely 
\begin{equation}
I-R,  \label{Amatrix}
\end{equation}%
where $R$ is an infinite matrix whose entries are $2\times 2$ block matrices
determined by weighted moments of the normalized differential of the second
kind on $S_{1}$. As before, the entries of $R$ are holomorphic in $\rho $, and
we show (Theorem \ref{Theorem_Det_rho}) that $\det (I-R)$ is defined and
holomorphic in a certain $\rho $-domain. The matrix $R$ then determines the
period matrix on $S$ (Theorem \ref{Theorem_periodg_rho}). We also discuss
several sewing scenarios for self-sewing a sphere, including one
(Proposition \ref{Prop_torus2}) where the Catalan series, familiar from
combinatorics \cite{St}, plays an unexpected role.

In Section \ref{Sect_Rho_Torus} we investigate in detail the self-sewing of
a twice-punctured torus with modulus $\tau \in \mathbb{H}_{1}$ to form a
genus two Riemann surface. As before sewing determines a map 
\begin{eqnarray}
F^{\rho }:\mathcal{D}^{\rho }(\tau ,w,\rho ) &\rightarrow &\mathbb{H}_{2}
\label{rhomap} \\
(\tau ,w,\rho ) &\mapsto &\Omega ^{(2)}(\tau ,w,\rho ),  \notag
\end{eqnarray}%
where now $\mathcal{D}^{\rho }(\tau ,w,\rho )\subseteq \mathbb{H}_{1}\times 
\mathbb{C}\times \mathbb{C}$ determines the admissible sewing parameters ($w$
describes the relative position of the punctures). Again we obtain (Theorem %
\ref{Theorem_rhoperiod}) explicit formulas for the entries of the matrix 
\begin{equation}
\Omega ^{(2)}(\tau ,w,\rho )=\left( 
\begin{array}{cc}
\Omega _{11}^{(2)} & \Omega _{12}^{(2)} \\ 
\Omega _{12}^{(2)} & \Omega _{22}^{(2)}%
\end{array}%
\right)  \label{Omegarhomatrix}
\end{equation}%
and show that $F^{\rho }$ is holomorphic. Roughly speaking, the entries of $%
R $ and the $\Omega _{ij}^{(2)}$ in this case are power series in $\rho $
with coefficients which are quasimodular and elliptic functions in the
variables $\tau ,w$. We provide a combinatorial description of (\ref%
{Omegarhomatrix}) in terms of a notion of chequered necklace suitably
modified compared to the $\epsilon $ case. A complicating factor is that $%
\Omega _{22}^{(2)}$ involves a logarithm of (the inverse square of) the
prime form on $S$. Because of this, it is necessary to pass to a covering
space $\hat{\mathcal{D}}^{\rho }$ of $\mathcal{D}^{\rho }$ before
equivariance properties can be considered. This is carried-out in Section %
\ref{Subsec_Frho}, where we construct a diagram

\begin{equation*}
\begin{array}{ccc}
\mathcal{D}^{\rho } & \overset{F^{\rho }}{\longrightarrow } & \ \mathbb{H}%
_{2} \\ 
\nwarrow &  & \nearrow \hat{F}^{\rho } \\ 
& \hat{\mathcal{D}}^{\rho } & 
\end{array}%
\end{equation*}

We show (Theorem \ref{Theorem_rho_equiv}) that the map $\hat{F}^{\rho }$ is
equivariant with respect to a group $L$ described as a semi-direct product
of $SL(2,\mathbb{Z})$ and the (nonabelian) \emph{Heisenberg} group $H\cong 
\mathbb{Z}^{1+2}$ (a $2$-step nilpotent group with center $\mathbb{Z}$).
Again $L$ acts on $\mathbb{H}_{2}$ via an embedding into $Sp(4,\mathbb{Z})$
and on $\hat{\mathcal{D}}^{\rho }$ in a manner prescribed in Theorem \ref%
{Theorem_l_L}. In Subsection \ref{Subsec_Frhoinvert} we obtain the expected
local invertibility of $F^{\rho }$ about a point of degeneration, which is a
bit more subtle than degeneration in the $\epsilon $-formalism. One of the
reasons for establishing the local invertibility results is that once
obtained, we have a way of comparing the two sewing domains $\mathcal{D}%
^{\epsilon }$ and $\mathcal{D}^{\rho }$, at least in some regions, by
looking at 
\begin{equation}
(F^{\rho {\ }}{)}^{-1}\circ F^{\epsilon }.  \label{funccomp}
\end{equation}%
In the final Section \ref{Sect_Epsilon_Rho}, (\ref{funccomp}) is briefly
discussed, where we observe that it is equivariant with respect to a common
subgroup of $G$ and $L$ isomorphic to $SL(2,\mathbb{Z})$. The Appendix
concludes with the explicit formulas for $\Omega ^{(2)}$ to $O(\epsilon
^{9}) $ in the $\epsilon $-formalism and to $O(\rho ^{5})$ in the $\rho $%
-formalism.

\section{Some Elliptic Functions}

\label{Sect_ellipticfunctions}

We briefly discuss a number of modular and elliptic-type functions that we
will need. The notation we introduce will be in force throughout the paper.
The Weierstrass elliptic function with periods\footnote{%
The period basis is more usually denoted by $\omega _{1},\omega _{2}$.} $%
\sigma ,\varsigma \in \mathbb{C}^{\ast }$ is defined by 
\begin{equation}
\wp (z,\sigma ,\varsigma )=\frac{1}{z^{2}}+\sum_{m,n\in \mathbb{Z},(m,n)\neq
(0,0)}[\frac{1}{(z-m\sigma -n\varsigma )^{2}}-\frac{1}{(m\sigma +n\varsigma
)^{2}}].  \label{Weierstrass}
\end{equation}%
Choosing $\varsigma =2\pi i$ and $\sigma =2\pi i\tau $ ($\tau $ will \emph{%
always} lie in the complex upper half-plane $\mathbb{H}$), we define 
\begin{eqnarray}
P_{2}(\tau ,z) &=&\wp (z,2\pi i\tau ,2\pi i)+E_{2}(\tau )  \notag \\
&=&\frac{1}{z^{2}}+\sum_{k=2}^{\infty }(k-1)E_{k}(\tau )z^{k-2}.  \label{P2}
\end{eqnarray}%
Here, $E_{k}(\tau )$ is equal to $0$ for $k$ odd, and for $k$ even is the
Eisenstein series \cite{Se} 
\begin{equation*}
E_{k}(\tau )=-\frac{B_{k}}{k!}+\frac{2}{(k-1)!}\sum_{n\geq 1}\sigma
_{k-1}(n)q^{n}.
\end{equation*}%
Here and below, we take $q=\exp (2\pi i\tau )$; $\sigma
_{k-1}(n)=\sum_{d\mid n}d^{k-1}$, and $B_{k}$ is the $k$th Bernoulli number
defined by 
\begin{eqnarray*}
\frac{t}{e^{t}-1}-1+\frac{t}{2} &=&\sum_{k\geq 2}B_{k}\frac{t^{k}}{k!} \\
&=&{\frac{1}{12}}{t}^{2}-{\frac{1}{720}}{t}^{4}+{\frac{1}{30240}}{t}^{6}+O({t%
}^{8}).
\end{eqnarray*}%
$P_{2}$ can be alternatively expressed as 
\begin{equation}
P_{2}(\tau ,z)=\frac{q_{z}}{(q_{z}-1)^{2}}+\sum_{n\geq 1}\frac{nq^{n}}{%
1-q^{n}}(q_{z}^{n}+q_{z}^{-n}),  \label{P2exp}
\end{equation}%
where $q_{z}=\exp (z)$. If $k\geq 4$ then $E_{k}(\tau )$ is a holomorphic
modular form of weight $k$ on $SL(2,\mathbb{Z})$. That is, it satisfies 
\begin{equation*}
E_{k}(\gamma \tau )=(c\tau +d)^{k}E_{k}(\tau )
\end{equation*}%
for all $\gamma =\left( 
\begin{array}{cc}
a & b \\ 
c & d%
\end{array}%
\right) \in SL(2,\mathbb{Z})$, where we use the standard notation

\begin{equation}
\gamma \tau =\frac{a\tau +b}{c\tau +d}.  \label{mobiusaction}
\end{equation}%
On the other hand, $E_{2}(\tau )$ has an exceptional transformation law 
\begin{equation}
E_{2}(\gamma \tau )=(c\tau +d)^{2}E_{2}(\tau )-\frac{c(c\tau +d)}{2\pi i}.
\label{gammaE2}
\end{equation}%
The first three Eisenstein series $E_{2}(\tau ),E_{4}(\tau ),E_{6}(\tau )$
are algebraically independent and generate a weighted polynomial algebra $Q=%
\mathbb{C}[E_{2}(\tau ),E_{4}(\tau ),E_{6}(\tau )]$ which, following \cite%
{KZ}, we call the algebra of \emph{quasimodular forms}.

We define $P_{1}(\tau ,z)$ by 
\begin{equation}
P_{1}(\tau ,z)=\frac{1}{z}-\sum_{k\geq 2}E_{k}(\tau )z^{k-1},  \label{P1}
\end{equation}%
where $P_{2}=-\frac{d}{dz}P_{1}$ and $P_{1}+zE_{2}$ is the classical
Weierstrass zeta function. $P_{1}$ is quasi-periodic with 
\begin{eqnarray}
P_{1}(\tau ,z+2\pi i) &=&P_{1}(\tau ,z),  \notag \\
P_{1}(\tau ,z+2\pi i\tau ) &=&P_{1}(\tau ,z)-1.  \label{P1quasiperiod}
\end{eqnarray}%
We also define $P_{0}(\tau ,z)$, up to a choice of the logarithmic branch,
by 
\begin{equation}
P_{0}(\tau ,z)=-\log (z)+\sum_{k\geq 2}\frac{1}{k}E_{k}(\tau )z^{k},
\label{P0}
\end{equation}%
where $P_{1}=-\frac{d}{dz}P_{0}$. We define the elliptic prime form $K(\tau
,z)\,$\ by \cite{Mu} 
\begin{equation}
K(\tau ,z)=\exp (-P_{0}(\tau ,z)),  \label{Primeform}
\end{equation}%
so that $P_{2}=\frac{d^{2}}{dz^{2}}\log K$. ($\exp (z^{2}E_{2}/2)K(\tau ,z)$
is the classical Weierstrass sigma function).\ $K(\tau ,z)$ is
quasi-periodic with 
\begin{eqnarray}
K(\tau ,z+2\pi i) &=&-K(\tau ,z),  \notag \\
K(\tau ,z+2\pi i\tau ) &=&-q_{z}^{-1}q^{-1/2}K(\tau ,z)  \label{Kquasiperiod}
\end{eqnarray}%
$K(\tau ,z)$ is an odd function of $z$ and can be expressed as 
\begin{equation}
K(\tau ,z)=-\frac{i\theta_{1}(\tau ,z)}{\eta (\tau )^{3}}=z+O(z^{3}),
\label{Kthetaeta}
\end{equation}%
where $\theta_{1}(\tau ,z)=\sum_{n\in \mathbb{Z}}\exp (\pi i\tau
(n+1/2)^{2}+(n+1/2)(z+i\pi ))$ and 
\begin{equation}
\eta (\tau )=q^{\frac{1}{24}}\prod_{n\geq 1}(1-q^{n})  \label{etafun}
\end{equation}%
is the Dedekind eta function.

Define elliptic functions $P_{k}(\tau ,z)\,$ for $k\geq 3\,$ from the
analytic expansion 
\begin{equation}
P_{1}(\tau ,z-w)=\sum_{k\geq 1}P_{k}(\tau ,z)w^{k-1}  \label{P1Pnexpansion}
\end{equation}
where 
\begin{equation}
P_{k}(\tau ,z)=\frac{(-1)^{k-1}}{(k-1)!}\frac{d^{k-1}}{dz^{k-1}}P_{1}(\tau
,z)=\frac{1}{z^{k}}+E_{k}+O(z).  \label{Pkdef}
\end{equation}

\noindent Finally, it is convenient to define for $k,l=1,2,\ldots $ 
\begin{eqnarray}
C(k,l) &=&C(k,l,\tau )=(-1)^{k+1}\frac{(k+l-1)!}{(k-1)!(l-1)!}E_{k+l}(\tau ),
\label{Ckldef} \\
D(k,l,z) &=&D(k,l,\tau ,z)=(-1)^{k+1}\frac{(k+l-1)!}{(k-1)!(l-1)!}%
P_{k+l}(\tau ,z).  \label{Dkldef}
\end{eqnarray}%
$\,$Note that $C(k,l)=C(l,k)$ and $D(k,l,z)=(-1)^{k+l}D(l,k,z)$. These
naturally arise in the analytic expansions (in appropriate domains) 
\begin{equation}
P_{2}(\tau ,z-w)=\frac{1}{(z-w)^{2}}+\sum_{k,l\geq 1}C(k,l)z^{l-1}w^{k-1},
\label{P2expansion}
\end{equation}%
and for $k\geq 1$ 
\begin{eqnarray}
P_{k+1}(\tau ,z) &=&\frac{1}{z^{k+1}}+\frac{1}{k}\sum_{l\geq 1}C(k,l)z^{l-1},
\label{Pkexpansion} \\
P_{k+1}(\tau ,z-w) &=&\frac{1}{k}\sum_{l\geq 1}D(k,l,w)z^{l-1}.
\label{Pkzwexpansion}
\end{eqnarray}

\section{The $\protect\epsilon $ Formalism for Sewing\ Together Two Riemann
Surfaces}

\label{Sect_Epsilon_g}In this section we review a general construction due
to Yamada \cite{Y} for "sewing" together two Riemann surfaces of genus $%
g_{1} $ and $g_{2}$ to form a surface of genus $g_{1}+g_{2}$. The principle
aim is to describe various structures such as the genus $g_{1}+g_{2}$ period
matrix in terms of data coming from the genus $g_{1}$ and genus $g_{2}$
surfaces. The basic method described below follows that of Yamada. However,
a significant number of changes have been made in order to express the final
formulas more neatly. We also discuss the holomorphic properties of the
period matrix and of a certain infinite dimensional determinant. In the next
section, this general formalism will be applied to the construction of a
genus two Riemann surface.

\subsection{The Bilinear Form $\protect\omega ^{(g)}$ and the Period Matrix $%
\Omega ^{(g)}$}

Consider a compact Riemann surface $\mathcal{S}$ of genus $g$ with canonical
homology basis $a_{1,}\ldots a_{g},b_{1},\ldots b_{g}$. In general there
exists $g$ holomorphic 1-forms $\nu_{i}^{(g)}$, $i=1,\ldots g$ which we may
normalize by \cite{FK1}, \cite{Sp} 
\begin{equation}
\oint_{a_{i}}\nu_{j}^{(g)}=2\pi i\delta_{ij}{.}  \label{norm}
\end{equation}%
These forms can be neatly encapsulated in a unique singular bilinear two
form $\omega ^{(g)}$, known as the \emph{normalized differential of the
second kind}. It is defined by the following properties \cite{Sp}, \cite{Mu}%
, \cite{Y}: 
\begin{equation}
\omega ^{(g)}(x,y)=(\frac{1}{(x-y)^{2}}+\text{regular terms})dxdy
\label{omegag}
\end{equation}%
for any local coordinates $x,y$, with normalization 
\begin{equation}
\int_{a_{i}}\omega ^{(g)}(x,\cdot )=0{\ }  \label{nugnorm}
\end{equation}%
{for }$i=1,\ldots ,g$. Using the Riemann bilinear relations, one finds that 
\begin{equation}
\nu_{i}^{(g)}(x)=\oint_{b_{i}}\omega ^{(g)}(x,\cdot ),  \label{nui}
\end{equation}%
with $\nu_{i}^{(g)}$normalized as in (\ref{norm}). The genus $g$ period
matrix $\Omega ^{(g)}$ is then defined by 
\begin{equation}
\Omega_{ij}^{(g)}=\frac{1}{2\pi i}\oint_{b_{i}}\nu_{j}^{(g)}\quad
\label{period}
\end{equation}%
for $i,j=1,\ldots ,g$. It is useful to also introduce the \emph{normalized
differential of the third kind} \cite{Mu}, \cite{Y}%
\begin{equation}
\omega_{p_{2}-p_{1}}^{(g)}(x)=\int_{p_{1}}^{p_{2}}\omega ^{(g)}(x,\cdot ),
\label{omp2p1}
\end{equation}%
for\ which $\oint_{a_{i}}\omega_{p_{2}-p_{1}}^{(g)}=0$. For a local
coordinate $x$ near $p_{a}$ for $a=1,2$ we have 
\begin{equation*}
\omega_{p_{2}-p_{1}}^{(g)}(x)=(\frac{(-1)^{a}}{x-p_{a}}+\text{regular terms}%
)dx.
\end{equation*}%
Both $\omega_{\ }^{(g)}(x,y)$ and $\omega_{p_{2}-p_{1}}^{(g)}(x)$\ can be
expressed in terms of the \emph{prime form} $%
K^{(g)}(x,y)(dx)^{-1/2}(dy)^{-1/2}$, a holomorphic form of weight $(-\frac{1%
}{2},-\frac{1}{2})$ with \cite{Mu} 
\begin{eqnarray}
\omega ^{(g)}(x,y) &=&\partial_{x}\partial_{y}\log K^{(g)}(x,y)dxdy,
\label{omgprime} \\
\omega_{p_{2}-p_{1}}^{(g)}(x) &=&\partial_{x}\log \frac{K^{(g)}(x,p_{2})}{%
K^{(g)}(x,p_{1})}dx.  \label{omp2p1prime}
\end{eqnarray}%
We also note that $K^{(g)}(x,y)=-K^{(g)}(y,x)$ and that $%
K^{(g)}(x,y)=x-y+O((x-y)^{3})$.

\begin{example}
\label{Exampleom1xy} For the genus one Riemann torus with periods $2\pi i$
and $2\pi i\tau $ along the $a$ and $b$ cycles, the holomorphic $1$-form
satisfying (\ref{norm}) in the usual parameterization is $\nu
_{1}^{(1)}=dz\, $. The normalized differential of the second kind is
determined by $P_{2}(\tau ,z)$ via 
\begin{equation}
\omega ^{(1)}(x,y)=P_{2}(\tau ,x-y)dxdy,.  \label{omega1}
\end{equation}%
In this case, (\ref{nugnorm}) and (\ref{nui}) follow from (\ref%
{P1quasiperiod}), and $\Omega_{11}^{(1)}=\tau $. The normalized differential
of the third kind is $\omega _{p_{2}-p_{1}}^{(1)}(x)=(P_{1}(\tau
,x-p_{2})-P_{1}(\tau ,x-p_{1}))dx$ and the prime form is $%
K^{(1)}(x,y)=K(\tau ,x-y)$ of (\ref{Primeform}).
\end{example}

It is well-known that $\Omega ^{(g)}$ is a complex symmetric matrix with
positive-definite imaginary part, i.e. $\Omega ^{(g)}\in \mathbb{H}_{g}$,
the genus $g$ Siegel complex upper half-space. The intersection form $\Xi $
is a natural non-degenerate symplectic bilinear form on the first homology
group $H_{1}(\mathcal{S},\mathbb{Z})\cong \mathbb{Z}^{2g}$, satisfying 
\begin{equation*}
\Xi (a_{i},a_{j})=\Xi (b_{i},b_{j})=0,\quad \Xi (a_{i},b_{j})=\delta
_{ij},\quad i,j=1,\ldots ,g.
\end{equation*}%
The genus $g$ symplectic group\footnote{Here and elsewhere, the transpose 
of a matrix or vector is denoted by T} is\ 
\begin{eqnarray*}
Sp(2g,\mathbb{Z}) &=&\{\gamma =\left( 
\begin{array}{ll}
A & B \\ 
C & D%
\end{array}%
\right) \in SL(2g,\mathbb{Z})| \\
AB^{T} &=&BA^{T},CD^{T}=D^{T}C,AD^{T}-BC^{T}=I_{g}\}.
\end{eqnarray*}%
It acts on $\mathbb{H}_{g}$ via 
\begin{equation}
\gamma .\Omega ^{(g)}{=(A\Omega ^{(g)}+B)(C\Omega ^{(g)}+D)^{-1},}
\label{eq: modtrans}
\end{equation}%
and naturally on $H_{1}(\mathcal{S},\mathbb{Z})$, where it preserves $\Xi $.

\subsection{The $\protect\epsilon $ Formalism for Sewing Two Riemann Surfaces%
}

\label{Subsec_Omega_eps_g}We now discuss a general method described by
Yamada \cite{Y} for calculating the bilinear form (\ref{omegag}) and hence
the period matrix on the surface formed by sewing together two Riemann
surfaces. Consider two Riemann surfaces $\mathcal{S}_{a}$ of genus $g_{a}$
for $a=1,2$. Choose a local coordinate $z_{a}$ on $\mathcal{S}_{a}$ in the
neighborhood of a point $p_{a}$, and consider the closed disk $\left\vert
z_{a}\right\vert \leq r_{a}$ for $r_{a}>0$ sufficiently small. (Note that
the choice $r_{a}=1$ is made in ref. \cite{Y}). Introduce a complex sewing
parameter $\epsilon $ where $|\epsilon |\leq r_{1}r_{2}$, and excise the
disk 
\begin{equation*}
\{z_{a},\left\vert z_{a}\right\vert \leq |\epsilon |/r_{\bar{a}}\}\subset 
\mathcal{S}_{a}
\end{equation*}%
to form a punctured surface 
\begin{equation*}
\hat{\mathcal{S}}_{a}=\mathcal{S}_{a}\backslash \{z_{a},\left\vert
z_{a}\right\vert \leq |\epsilon |/r_{\bar{a}}\}.
\end{equation*}%
Here and below, we use the convention 
\begin{equation}
\overline{1}=2,\quad \overline{2}=1  \label{bardef}
\end{equation}%
Define the annulus

\begin{equation*}
\mathcal{A}_{a}=\{z_{a},|\epsilon |/r_{\bar{a}}\leq \left\vert
z_{a}\right\vert \leq r_{a}\}\subset \hat{\mathcal{S}}_{a},
\end{equation*}%
and identify $\mathcal{A}_{1}$ and $\mathcal{A}_{2}$ as a single region $%
\mathcal{A}=\mathcal{A}_{1}\simeq \mathcal{A}_{2}$ via the sewing relation 
\begin{equation}
z_{1}z_{2}=\epsilon .  \label{pinch}
\end{equation}

\begin{center}
\begin{picture}(300,100)

\put(50,52){\qbezier(-30,18)(-10,10)(10,18)}
\put(50,52){\qbezier(10,18)(50,35)(90,18)}
\put(50,48){\qbezier(-30,-18)(-10,-10)(10,-18)}
\put(50,48){\qbezier(10,-18)(50,-35)(90,-18)}

\put(45,50){\qbezier(25,0)(45,17)(60,0)}
\put(45,50){\qbezier(20,2)(45,-17)(65,2)}

\put(175,52){\qbezier(90,18)(110,10)(130,18)}
\put(175,52){\qbezier(10,18)(50,35)(90,18)}
\put(175,48){\qbezier(90,-18)(110,-10)(130,-18)}
\put(175,48){\qbezier(10,-18)(50,-35)(90,-18)}

\put(200,50){\qbezier(25,0)(45,17)(60,0)}
\put(200,50){\qbezier(20,2)(45,-17)(65,2)}

\put(140,50){\circle{16}}
\put(140,50){\circle{40}}

\put(140,50){\vector(-1,-2){0}}
\put(50,50){\qbezier(90,0)(100,15)(90,30)}%
\put(140,90){\makebox(0,0){$z_1=0$}}

\put(140,50){\line(-1,1){14.1}}
\put(127,55){\makebox(0,0){$r_1$}}

\put(140,50){\line(1,0){8}}
\put(145,50){\vector(1,4){0}}
\put(55,20){\qbezier(90,4)(85,17)(90,30)}%
\put(150,15){\makebox(0,0){$|\epsilon|/r_2$}}


\put(50,50){\makebox(0,0){$\mathcal{S}_1$}}


\put(185,50){\circle{16}}
\put(185,50){\circle{40}}

\put(185,50){\vector(1,-2){0}}
\put(95,50){\qbezier(90,0)(80,15)(90,30)}%
\put(185,90){\makebox(0,0){$z_2=0$}}

\put(185,50){\line(-1,-1){14.1}}
\put(171,45){\makebox(0,0){$r_2$}}

\put(185,50){\line(1,0){8}}
\put(190,50){\vector(1,4){0}}
\put(90,20){\qbezier(100,4)(95,17)(100,30)}%
\put(190,15){\makebox(0,0){$|\epsilon|/r_1$}}

\put(280,50){\makebox(0,0){$\mathcal{S}_2$}}

\end{picture}

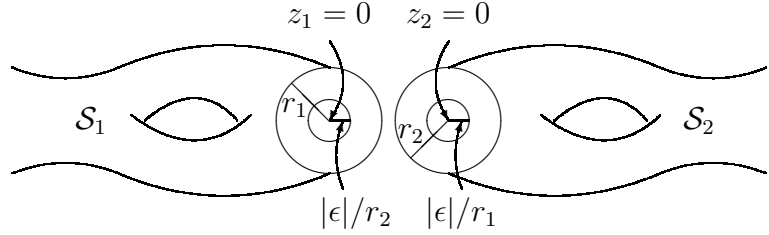
{\small Fig. 1 Sewing Two Riemann Surfaces}
\end{center}

In this way we obtain a compact Riemann surface $\hat{\{\mathcal{S}}%
_{1}\backslash \mathcal{A}_{1}\}\cup \{\hat{\mathcal{S}}_{2}\backslash 
\mathcal{A}_{2}\}\cup \mathcal{A}$ of genus $g_{1}+g_{2}$. The sewing
relation (\ref{pinch}) can be considered to be a parameterization of a
cylinder connecting the two punctured Riemann surfaces. Noting the
notational differences with ref. \cite{Y}, the genus $g_{1}+g_{2}$
normalized differential of the second kind $\omega ^{(g_{1}+g_{2})}\,$of (%
\ref{omegag}) enjoys the following properties:

\begin{theorem}[Ref. \protect\cite{Y}, Theorem 1, Theorem 4]
\label{Theoremomg1g2holo}

\begin{description}
\item[ ] 

\item[(a)] $\omega ^{(g_{1}+g_{2})}$ is holomorphic in $\epsilon $ for $%
|\epsilon |<r_{1}r_{2}$;

\item[(b)] $\lim_{\epsilon \rightarrow 0}\omega ^{(g_{1}+g_{2})}(x,y)=\omega
^{(g_{a})}(x,y)\delta_{ab}$ for $x\in \hat{\mathcal{S}}_{a},y\in \hat{%
\mathcal{S}}_{b}$, $a,b=1,2$. $\qed$
\end{description}
\end{theorem}

Regarded as a power series in $\epsilon$, the coefficients of\ $%
\omega^{(g_{1}+g_{2})}$ can be calculated from $\omega ^{(g_{1})}$ and $%
\omega ^{(g_{2})}$ as follows. Let $\mathcal{C}_{a}(z_{a}) \subset \mathcal{A%
}_a$ denote a simple, closed, anti-clockwise oriented contour parameterized
by $z_{a}$, surrounding the puncture at $z_{a}=0$. Note that $\mathcal{C}%
_{1}(z_{1})$ may be deformed to $-\mathcal{C}_{2}(z_{2})$ via (\ref{pinch}).
Then one finds \cite{Y}:

\begin{lemma}[op.cite., Lemma 4]
\label{Lemma_inteqn1} 
\begin{equation}
\omega ^{(g_{1}+g_{2})}(x,y)=\omega ^{(g_{a})}(x,y)\delta _{ab}+\frac{1}{%
2\pi i}\oint_{\mathcal{C}_{a}(z_{a})}(\omega ^{(g_{1}+g_{2})}(y,z_{a})\text{ 
}\int^{z_{a}}\omega ^{(g_{a})}(x,\cdot ))  \label{omegagint}
\end{equation}%
for $x\in \hat{\mathcal{S}}_{a}$, $y\in \hat{\mathcal{S}}_{b}$ and $a,b=1,2$%
. $\qed$
\end{lemma}

Define weighted moments for $\omega ^{(g_{1}+g_{2})}$ for $k,l=1,2,\ldots $
by 
\begin{eqnarray}
&&X_{ab}(k,l)=X_{ab}(k,l,\epsilon )  \notag \\
&=&\frac{\epsilon ^{(k+l)/2}}{\sqrt{kl}}\frac{1}{(2\pi i)^{2}}\oint_{%
\mathcal{C}_{a}(u)}\oint_{\mathcal{C}_{b}(v)}u^{-k}v^{-l}\omega
^{(g_{1}+g_{2})}(u,v).  \label{Xijdef}
\end{eqnarray}%
The $\epsilon ^{(k+l)/2}/\sqrt{kl}$ $\,$factor is introduced for later
convenience. Note that 
\begin{equation}
X_{ab}(k,l)=X_{ba}(l,k)  \label{Xabsym}
\end{equation}
and that $\epsilon ^{-(k+l)/2}X_{ab}(k,l,\epsilon )$ is holomorphic in $%
\epsilon $ for $|\epsilon |<r_{1}r_{2}$ from Theorem \ref{Theoremomg1g2holo}. 
We define $X_{ab}=(X_{ab}(k,l))$ to be the infinite matrix indexed by $k,l$.

Next define a set of holomorphic $1$-forms on $\hat{\mathcal{S}}_{a}$ by 
\begin{equation}
a_{a}(k,x)=a_{a}(k,x,\epsilon )=\frac{\epsilon ^{k/2}}{2\pi i\sqrt{k}}\oint_{%
\mathcal{C}_{a}(z_{a})}z_{a}^{-k}\omega ^{(g_{a})}(x,z_{a}),  \label{akdef}
\end{equation}%
and define $a_{a}(x)=(a_{a}(k,x))$ to be the infinite row vector indexed by $%
k$. Note from (\ref{omegag}) that for $x,y\in \hat{\mathcal{S}}_{a}$ with $%
x\neq 0$ we have 
\begin{eqnarray}
\omega ^{(g_{a})}(x,y) &=&\sum_{k\geq 1}[\frac{1}{2\pi i}\oint_{\mathcal{C}%
_{a}(z_{a})}z_{a}^{-k}\omega ^{(g_{a})}(x,z_{a})]y^{k-1}dy,  \notag \\
&=&\sum_{k\geq 1}\sqrt{k}\epsilon ^{-k/2}a_{a}(k,x)y^{k-1}dy.
\label{omegagiyexp}
\end{eqnarray}

Using Lemma \ref{Lemma_inteqn1} we have:

\begin{lemma}
\label{LemmaaiXijaj} $\omega ^{(g_{1}+g_{2})}(x,y)$ is given by 
\begin{equation}
\omega ^{(g_{1}+g_{2})}(x,y)=\left\{ 
\begin{array}{cc}
\omega ^{(g_{a})}(x,y)+a_{a}(x)X_{\bar{a}\bar{a}}a_{a}^{T}(y) & x,y\in \hat{%
\mathcal{S}}_{a}, \\ 
a_{a}(x)(-I+X_{\bar{a}a})a_{\bar{a}}^{T}(y) & x\in \hat{\mathcal{S}}%
_{a},\quad y\in \hat{\mathcal{S}}_{\bar{a}}.%
\end{array}%
\right.  \label{omgsisj}
\end{equation}
\end{lemma}

\noindent \textbf{Proof.} From (\ref{omegagiyexp}) it follows that 
\begin{equation}
\int_{0}^{z_{a}}\omega ^{(g_{a})}(x,\cdot )=\sum_{k\geq 1}\frac{\epsilon
^{-k/2}}{\sqrt{k}}a_{a}(k,x)z_{a}^{k}.  \label{omegagint1}
\end{equation}%
Let $x,y\in \hat{\mathcal{S}}_{1}$. Using (\ref{pinch}), (\ref{omegagint})
and (\ref{omegagint1}) we find that $\omega ^{(g_{1}+g_{2})}(x,y)-\omega
^{(g_{1})}(x,y)$ is given by 
\begin{eqnarray}
&&\sum_{k\geq 1}\frac{\epsilon ^{-k/2}}{\sqrt{k}}a_{1}(k,x)(-\frac{\epsilon
^{k}}{2\pi i}\oint_{\mathcal{C}_{2}(z_{2})}z_{2}^{-k}\omega
^{(g_{1}+g_{2})}(y,z_{2}))\text{ }  \label{om2minusom1} \\
&=&\sum_{k,l\geq 1}\frac{\epsilon ^{(k+l)/2}}{\sqrt{kl}}a_{1}(k,x)a_{1}(l,y)%
\frac{1}{(2\pi i)^{2}}\oint_{\mathcal{C}_{2}(w_{2})}\oint_{\mathcal{C}%
_{2}(z_{2})}z_{2}^{-k}w_{2}^{-l}\omega ^{(g_{1}+g_{2})}(z_{2},w_{2}),  \notag
\end{eqnarray}%
giving (\ref{omgsisj}) for $x,y\in \hat{\mathcal{S}}_{1}$.

For $x\in \hat{\mathcal{S}}_{1},$ $y\in \hat{\mathcal{S}}_{2}$ we find that $%
\omega ^{(g_{1}+g_{2})}(x,y)$ is given by 
\begin{eqnarray}
&&\sum_{k\geq 1}\frac{\epsilon ^{-k/2}}{\sqrt{k}}a_{1}(k,x)(-\frac{\epsilon
^{k}}{2\pi i}\oint_{\mathcal{C}_{2}(z_{2})}z_{2}^{-k}\omega
^{(g_{1}+g_{2})}(y,z_{2}))  \label{omegags2s1} \\
&=&-\sum_{k\geq 1}a_{1}(k,x)a_{2}(k,y)  \notag \\
&&+\sum_{k,l\geq 1}\frac{\epsilon ^{(k+l)/2}}{\sqrt{kl}}a_{1}(k,x)a_{2}(l,y)%
\frac{1}{(2\pi i)^{2}}\oint_{\mathcal{C}_{1}(z_{1})}\oint_{\mathcal{C}%
_{2}(z_{2})}z_{1}^{-l}z_{2}^{-k}\omega ^{(g_{1}+g_{2})}(z_{1},z_{2}),  \notag
\end{eqnarray}%
giving (\ref{omgsisj}). A similar analysis follows for $x,y\in \hat{\mathcal{%
S}}_{1}$ and $x\in \hat{\mathcal{S}}_{2},$ $y\in \hat{\mathcal{S}}_{1}$. $%
\qed$

We next compute the explicit form of the moment matrix $X_{ab}$ in terms of
the moments of $\omega ^{(g_{a})}$, which we denote by 
\begin{eqnarray}
A_{a}(k,l) &=&A_{a}(k,l,\epsilon )=\frac{\epsilon ^{(k+l)/2}}{(2\pi i)^{2}%
\sqrt{kl}}\oint_{\mathcal{C}_{a}(x)}\oint_{\mathcal{C}_{a}(y)}x^{-k}y^{-l}%
\omega ^{(g_{a})}(x,y)  \notag \\
&=&\frac{\epsilon ^{k/2}}{2\pi i\sqrt{k}}\oint_{\mathcal{C}%
_{a}(x)}x^{-k}a_{a}(l,x).  \label{Adef}
\end{eqnarray}%
Note from (\ref{omegag}) that for $x,y\in \hat{\mathcal{S}}_{a}$ we have 
\begin{eqnarray}
&&\omega ^{(g_{a})}(x,y)-\frac{dxdy}{(x-y)^{2}}  \notag \\
&=&\sum_{k,l\geq 1}\{\frac{1}{(2\pi i)^{2}}\oint_{\mathcal{C}_{a}(u)}\oint_{%
\mathcal{C}_{a}(v)}u^{-k}v^{-l}\omega ^{(g_{a})}(u,v)\}x^{k-1}y^{l-1}dxdy, 
\notag \\
&=&\sum_{k,l\geq 1}\sqrt{kl}\epsilon ^{-(k+l)/2}A_{a}(k,l,\epsilon
)x^{k-1}y^{l-1}dxdy.  \label{omegagi}
\end{eqnarray}

\begin{proposition}
\label{PropXij}\ The matrices $X_{ab}\,$are given in terms of $A_{a}$ by 
\begin{eqnarray}
X_{aa} &=&A_{a}(I-A_{\bar{a}}A_{a})^{-1},  \label{Xaa} \\
X_{a\bar{a}} &=&I-(I-A_{a}A_{\bar{a}})^{-1}.  \label{Xaabar}
\end{eqnarray}%
Here, 
\begin{equation}
(I-A_{a}A_{\bar{a}})^{-1}=\sum_{n\geq 0}(A_{a}A_{\bar{a}})^{n}
\label{inverseexp}
\end{equation}%
and is convergent as a power series in $\epsilon $ for $|\epsilon
|<r_{1}r_{2}$.
\end{proposition}

\noindent \textbf{Proof.} Compute $X_{11}$ from (\ref{om2minusom1}) to find 
\begin{eqnarray*}
X_{11}(k,l) &=&\frac{\epsilon ^{(k+l)/2}}{\sqrt{kl}}\frac{1}{(2\pi i)^{2}}%
\oint_{\mathcal{C}_{1}(x)}\oint_{\mathcal{C}_{1}(y)}x^{-k}y^{-l}\omega
^{(g_{1})}(x,y) \\
&&-\sum_{m\geq 1}[\frac{\epsilon ^{k/2}}{2\pi i}\frac{1}{\sqrt{k}}\oint_{%
\mathcal{C}_{1}(x)}x^{-k}a_{1}(m,x). \\
&&\frac{\epsilon ^{(m+l)/2}}{\sqrt{ml}}\frac{1}{(2\pi i)^{2}}\oint_{\mathcal{%
C}_{1}(y)}\oint_{\mathcal{C}_{2}(z_{2})}y^{-l}z_{2}^{-m}\omega
^{(g_{1}+g_{2})}(y,z_{2})],
\end{eqnarray*}%
and similarly for $X_{22}$. Thus using (\ref{Adef}) and recalling (\ref%
{Xabsym}) we have 
\begin{equation}
X_{aa}=A_{a}(I-X_{\bar{a}a}).  \label{Xaarel}
\end{equation}%
We may find $X_{12}$ from (\ref{omegags2s1}) as follows: 
\begin{eqnarray*}
X_{12}(k,l) &=&-\sum_{m\geq 1}\frac{\epsilon ^{k/2}}{2\pi i}\frac{1}{\sqrt{k}%
}\oint_{\mathcal{C}_{1}(x)}a_{1}(m,x)x^{-k} \\
&&(\frac{\epsilon ^{(m+l)/2}}{\sqrt{ml}}\frac{1}{(2\pi i)^{2}}\oint_{%
\mathcal{C}_{2}(y)}\oint_{\mathcal{C}_{2}(z_{2})}y^{-l}z_{2}^{-m}\omega
^{(2)}(y,z_{2}))
\end{eqnarray*}%
and similarly for $X_{21}$ i.e. 
\begin{equation}
X_{a\bar{a}}=-A_{a}X_{\bar{a}\bar{a}}.  \label{Xaabarrel}
\end{equation}%
\qquad\ 

Define infinite block matrices 
\begin{equation}
X=\left[ 
\begin{array}{cc}
X_{11} & X_{12} \\ 
X_{21} & X_{22}%
\end{array}%
\right] ,\qquad A=\left[ 
\begin{array}{cc}
A_{1} & 0 \\ 
0 & A_{2}%
\end{array}%
\right] ,\qquad Q=\left[ 
\begin{array}{cc}
0 & -A_{1} \\ 
-A_{2} & 0%
\end{array}%
\right]  \label{XAQ_def}
\end{equation}%
so that (\ref{Xaarel}, \ref{Xaabarrel}) can be combined as%
\begin{equation}
X=A+QX,  \label{XQrel}
\end{equation}%
so that 
\begin{equation}
X=(I-Q)^{-1}A.  \label{XQsol}
\end{equation}%
Here $(I-Q)^{-1}=\sum_{n\geq 0}Q^{n}$ which we now show converges for $%
|\epsilon |<r_{1}r_{2}$.

Consider $X=A+AX+Q^{2}X$. Then since $Q^{2}=\mathrm{diag}%
(A_{1}A_{2},A_{2}A_{1})$ we obtain iterative relations 
\begin{eqnarray}
X_{aa} &=&A_{a}(I+A_{\bar{a}}X_{aa}),  \label{Xaaself} \\
I-X_{a\bar{a}} &=&I+A_{a}A_{\bar{a}}(I-X_{a\bar{a}}).  \label{Xaabarself}
\end{eqnarray}%
Now $\epsilon ^{-(k+l)/2}X_{ab}(k,l)$ is holomorphic in $\epsilon $ for $%
|\epsilon |<r_{1}r_{2}$ by Theorem \ref{Theoremomg1g2holo}. Therefore, $%
X_{ab}(k,l)$ has a series expansion in $\epsilon ^{1/2}$ convergent for $%
|\epsilon |<r_{1}r_{2}$. Then (\ref{Xaaself}) implies $X_{aa}(k,l)=(%
\sum_{n=0}^{N}(A_{a}A_{\bar{a}})^{n}A_{a})(k,l)+O(\epsilon ^{(k+l)/2+2N+1})$
where the coefficient of each power of $\epsilon $ consists of a finite sum
of finite products of $A_{1}$ and $A_{2}$. Hence 
\begin{equation*}
X_{aa}=\sum_{n=0}^{\infty }A_{a}(A_{\bar{a}}A_{a})^{n}=A_{a}(I-A_{\bar{a}%
}A_{a})^{-1},
\end{equation*}%
converges for $|\epsilon |<r_{1}r_{2}$. A similar argument holds for $X_{a%
\bar{a}}$ where one finds 
\begin{equation*}
X_{a\bar{a}}=\sum_{n=1}^{\infty }(A_{a}A_{\bar{a}})^{n}=I-(I-A_{a}A_{\bar{a}%
})^{-1}
\end{equation*}%
converges for $|\epsilon |<r_{1}r_{2}$. Finally 
\begin{equation*}
(I-Q)^{-1}=\sum_{n\geq 0}Q^{n}=\sum_{n\geq 0}Q^{2n}(I+Q)
\end{equation*}%
is therefore also convergent for $|\epsilon |<r_{1}r_{2}$.$\qed$

The invertibility of the infinite matrix $I-A_{1}A_{2}$ for $|\epsilon
|<r_{1}r_{2}$ is crucial in the $\epsilon $ sewing formalism. We now define
an infinite determinant of $I-A_{1}A_{2}$ which we show is a 
holomorphic function of $\epsilon $ for $|\epsilon |<r_{1}r_{2}$. This
determinant plays a dominant role in the sequel to this work \cite{MT2}.
Firstly, since $A_{1}(k,m)A_{2}(m,l)=O(\epsilon ^{m+(k+l)/2})$ we may define
a $(2N-3)\times (2N-3)$ matrix 
\begin{equation}
T_{N}(k,l)=\sum_{1\leq m\leq N-(k+l)/2}A_{1}(k,m)A_{2}(m,l),  \label{TN}
\end{equation}%
for $1\leq k,l\leq 2N-3$. $T_{N}$ is a truncated approximation for $%
A_{1}A_{2}$ to $O(\epsilon ^{N})$

\begin{equation*}
A_{1}A_{2}=\left( 
\begin{array}{lll}
T_{N} & 0 & \cdots \\ 
0 & 0 & \cdots \\ 
\vdots & \vdots & \ddots%
\end{array}%
\right) +O(\epsilon ^{N+1}).
\end{equation*}%
We may then define as formal power series in $\epsilon $ to $O(\epsilon^{N})$ 
the expressions\footnote{For the sake of notational simplicity we denote 
{\em both} the usual finite dimensional and the defined infinite dimensional 
determinants by $\det$.}
\begin{eqnarray}
\det (I-A_{1}A_{2}) &=&\det (I_{N}-T_{N})+O(\epsilon ^{N+1}),  \label{Det} \\
Tr\log (I-A_{1}A_{2}) &=&Tr\log (I_{N}-T_{N})+O(\epsilon ^{N+1}),
\label{Trace_Log}
\end{eqnarray}%
where $Tr\log (I_{N}-T_{N})=-\sum_{n=1}^{N/2}\frac{1}{n}Tr(T_{N}^{n})+O(%
\epsilon ^{N+1})$. Comparing the finite matrix contributions of (\ref{Det})
and (\ref{Trace_Log}) we have

\begin{lemma}
\label{Lemma_logdet} As formal power series in $\epsilon $ 
\begin{equation}
\log \det (I-A_{1}A_{2})=Tr\log (I-A_{1}A_{2}).\qed  \label{logdet}
\end{equation}
\end{lemma}

We now show that $Tr\log (I-A_{1}A_{2})$ is holomorphic in $\epsilon $ for $%
|\epsilon |<r_{1}r_{2}$ so that:

\begin{theorem}
\label{Theorem_Det} $\det (I-A_{1}A_{2})$ is non-vanishing and holomorphic
in $\epsilon $ for $|\epsilon |<r_{1}r_{2}$.
\end{theorem}

\noindent \textbf{Proof.} Let $\omega
^{(g_{1}+g_{2})}=f(z_{1},z_{2},\epsilon )dz_{1}dz_{2}$ for $|z_{a}|\leq
r_{a} $. Then $f(z_{1},z_{2},\epsilon )$ is holomorphic in $\epsilon $ for $%
|\epsilon |\leq r$ for $r<r_{1}r_{2}$ from Theorem \ref{Theoremomg1g2holo}.
Apply Cauchy's inequality to the coefficient functions for $%
f(z_{1},z_{2},\epsilon )=\sum_{n\geq 0}f_{n}(z_{1},z_{2})\epsilon ^{n}$ to
find 
\begin{equation}
|f_{n}(z_{1},z_{2})|\leq \frac{M}{r^{n}},  \label{Cauchy}
\end{equation}%
for $M=\sup_{|z_{a}|\leq r_{a},|\epsilon |\leq r}|f(z_{1},z_{2},\epsilon )|$%
. Consider 
\begin{equation}
\mathcal{I}=\frac{1}{(2\pi i)^{2}}\oint_{\mathcal{C}_{r_{1}}(z_{_{1}})}%
\oint_{\mathcal{C}_{r_{2}}(z_{_{2}})}\omega
^{(g_{1}+g_{2})}(z_{1},z_{2})\log (1-\frac{\epsilon }{z_{1}z_{2}}),
\label{Siberian}
\end{equation}%
for $\mathcal{C}_{r_{a}}(z_{a})$ the contour with $|z_{a}|=r_{a}$. Then
using (\ref{Cauchy}) we find 
\begin{eqnarray*}
|\mathcal{I}| &\leq &\sum_{n\geq 0}\frac{|\epsilon |^{n}}{(2\pi )^{2}}\oint_{%
\mathcal{C}_{r_{1}}(z_{_{1}})}\oint_{\mathcal{C}%
_{r_{2}}(z_{_{2}})}|f_{n}(z_{1},z_{2})\log (1-\frac{\epsilon }{z_{1}z_{2}}%
)dz_{1}dz_{2}| \\
&\leq &\sum_{n\geq 0}M.\frac{|\epsilon |^{n}}{r^{n}}.|\log (1-\frac{%
|\epsilon |}{r_{1}r_{2}})|.r_{1}r_{2},
\end{eqnarray*}%
i.e. $\mathcal{I}$ is absolutely convergent and thus holomorphic in $%
\epsilon $ for $|\epsilon |\leq r<r_{1}r_{2}$. Since $%
|z_{1}z_{2}|=r_{1}r_{2} $ we may alternatively expand in $\epsilon
/z_{1}z_{2}$ to obtain 
\begin{eqnarray*}
\mathcal{I} &=&-\sum_{k\geq 1}\frac{\epsilon ^{k}}{k}\frac{1}{(2\pi i)^{2}}%
\oint_{\mathcal{C}_{r_{1}}(z_{_{1}})}\oint_{\mathcal{C}_{r_{2}}(z_{_{2}})}%
\omega ^{(g_{1}+g_{2})}(z_{1},z_{2})z_{1}^{-k}z_{2}^{-k} \\
&=&-TrX_{12},
\end{eqnarray*}%
where $TrX_{12}=\sum_{k\geq 1}X_{12}(k,k)$ for $X_{12}$ of (\ref{Xaabar}).
But (\ref{Xaabarself}) implies 
\begin{equation*}
TrX_{12}=-\sum_{n\geq 1}Tr((A_{1}A_{2})^{n}),
\end{equation*}%
is absolutely convergent for $|\epsilon |<r_{1}r_{2}$. Hence we find 
\begin{equation*}
Tr\log (I-A_{1}A_{2})=-\sum_{n\geq 1}\frac{1}{n}Tr((A_{1}A_{2})^{n}),
\end{equation*}%
is also absolutely convergent for $|\epsilon |<r_{1}r_{2}$. Thus from Lemma %
\ref{Lemma_logdet}, $\det (I-A_{1}A_{2})$ is non-vanishing and holomorphic
for $|\epsilon |<r_{1}r_{2}$. $\qed$

These determinant properties can also be expressed in terms of the block
matrix $Q$ using\footnote{See the previous footnote.}

\begin{lemma}
\label{Lemma_detQdetA1A2} $\det (I\pm Q)=$ $\det (I-A_{1}A_{2})$.
\end{lemma}

\textbf{Proof.} Let $Q_{N}$ be the truncated approximation for $Q$ to $%
O(\epsilon ^{N})$. Then one finds $\det (I+Q_{N})=\det (I-Q_{N})$. But $\det
(I+Q_{N})\det (I-Q_{N})=\det (I-Q_{N}^{2})=\det (I-T_{N})^{2}$ for $T_{N}$
of (\ref{TN}) and the result follows. $\qed$

The sewn genus $g_{1}+g_{2}$ Riemann surface naturally inherits a basis of
cycles labelled $\{a_{s_{1}},b_{s_{1}}|s_{1}=1,\ldots ,g_{1}\}$ and $%
\{a_{s_{2}},b_{s_{2}}|s_{2}=g_{1}+1,\ldots ,g_{1}+g_{2}\}$ from the genus $%
g_{1}\,$and genus $g_{2}$ surfaces respectively. Integrating $\omega
^{(g_{1}+g_{2})}$ along the $b$ cycles then gives us the holomorphic 1-forms
and period matrix. For $a=1,2$ we define 
\begin{equation}
\alpha_{s_{a}}^{(g_{a})}(k)=\oint_{b_{s_{a}}}a_{a}(k),  \label{alphadef}
\end{equation}%
and the infinite vector $\alpha_{s_{a}}^{(g_{a})}=(\alpha
_{s_{a}}^{(g_{a})}(k))$. Then we find using (\ref{nui}) and (\ref{period})
together with Lemma \ref{LemmaaiXijaj} and Proposition \ref{PropXij} that 
\cite{Y}:

\begin{theorem}[op. cite. Theorem 4.]
\label{Theoremperiodgeps} $\Omega ^{(g_{1}+g_{2})}$ is holomorphic in $%
\epsilon $ for $|\epsilon |<r_{1}r_{2}$, and is given by 
\begin{eqnarray}
2\pi i\Omega_{s_{1}t_{1}}^{(g_{1}+g_{2})} &=&2\pi i\Omega
_{s_{1}t_{1}}^{(g_{1})}+\alpha
_{s_{1}}^{(g_{1})}(A_{2}(I-A_{1}A_{2})^{-1})\alpha_{t_{1}}^{(g_{1})T},
\label{Om11X} \\
2\pi i\Omega_{s_{2}t_{2}}^{(g_{1}+g_{2})} &=&2\pi i\Omega
_{s_{2}t_{2}}^{(g_{2})}+\alpha
_{s_{2}}^{(g_{2})}(A_{1}(I-A_{1}A_{2})^{-1})\alpha_{t_{2}}^{(g_{2})T},
\label{Om22X} \\
2\pi i\Omega_{s_{1}s_{2}}^{(g_{1}+g_{2})} &=&-\alpha
_{s_{1}}^{(g_{1})}(I-A_{2}A_{1})^{-1}\alpha_{s_{2}}^{(g_{2})T}  \notag \\
&=&-\alpha_{s_{2}}^{(g_{2})}(I-A_{1}A_{2})^{-1}\alpha_{s_{1}}^{(g_{1})T}.
\label{Om12X}
\end{eqnarray}%
with $s_{1},t_{1}\in \{1,\ldots ,g_{1}\}$ and $s_{2},t_{2}\in
\{g_{1}+1,\ldots ,g_{1}+g_{2}\}$. $\qed$
\end{theorem}

\begin{example}
\textbf{\label{Example_g_to_sphere}}Let $\mathcal{S}_{1}$ be a genus $g$
surface and $\mathcal{S}_{2}$ the Riemann sphere $\mathbb{C\cup \{\infty \}}$
with bilinear form 
\begin{equation}
\omega ^{(0)}(x,y)=\frac{dxdy}{(x-y)^{2}},\quad x,y\in \mathcal{S}_{2}.
\label{omega0}
\end{equation}%
Choose $p_{2}=0$ with $z_{2}\in \mathbb{C}$ as the local coordinate on $%
\mathcal{S}_{2}$. Then $a_{2}(k,x)=\sqrt{k}\epsilon ^{k/2}x^{-k-1}dx$ from (%
\ref{akdef}) and hence $A_{2}=0$ from (\ref{Adef}). Thus $X_{22}=A_{1}$ and $%
X_{11}=X_{12}=X_{21}=0$. Then one can check that the RHS of (\ref{omgsisj})
reproduces $\omega ^{(g)}$ directly for $a=b=1$ whereas, using (\ref{pinch}%
), it follows for $a=b=2$ from (\ref{omegagi}) and for $a\neq b$ from (\ref%
{omegagiyexp}).
\end{example}

Let us now consider the holomorphic properties of $\Omega ^{(g_{1}+g_{2})}$.
From Theorem \ref{Theoremomg1g2holo}, $\omega ^{(g_{1}+g_{2})}(x,y)$ is
holomorphic in $\epsilon $ for $|\epsilon |<r_{1}r_{2}$ and therefore $%
\Omega ^{(g_{1}+g_{2})}$ is also. We now show that if $\omega ^{(g_{a})}$ is
holomorphic with respect to a complex parameter (such as one of the modular
parameters of the Riemann surface $\mathcal{S}_{a}$) then $\omega
^{(g_{1}+g_{2})}$ and therefore $\Omega ^{(g_{1}+g_{2})}$ is also
holomorphic with respect to that parameter. To this end we firstly prove the
following elementary lemma:

\begin{lemma}
\label{Lem_fxyhol} Let $f(x,y)$ be a complex function holomorphic in $x$ for 
$|x|<R$ with expansion $f(x,y)=\sum_{m\geq 0}c_{m}(y)x^{m}$. Suppose that
each $c_{m}(y)$ is holomorphic and that $f(x,y)$ is continuous in $y$ for $%
|y|<S$. Then $f(x,y)$ is also holomorphic in $y$ for $|y|<S$.
\end{lemma}

\noindent \textbf{Proof.} Define the compact region $\mathcal{R}%
=\{(x,y):|x|\leq R_{-},|y|\leq S_{-}\}$ for $R_{-}=R-\delta _{1}$ and $%
S_{-}=S-\delta _{2}$ for $\delta _{1},\delta _{2}>0$. $f$ is continuous in
the compact region $\mathcal{R}$ and hence $|f(x,y)|\leq M\equiv \sup_{%
\mathcal{R}}|f|$. Apply Cauchy's inequality to $f$ as a holomorphic function
of $x$ for $|x|\leq R_{-}$ to find 
\begin{equation*}
|c_{m}(y)|\leq \frac{\sup_{|x|\leq R_{-}}|f(x,y)|}{R_{-}^{m}}\leq \frac{M}{%
R_{-}^{m}}.
\end{equation*}%
But $c_{m}(y)$ is holomorphic for $|y|\leq S_{-}$ with expansion $%
c_{m}(y)=\sum_{n\geq 0}c_{mn}y^{n}$ so that applying Cauchy's inequality
again gives%
\begin{equation*}
|c_{mn}|\leq \frac{\sup_{|y|\leq S_{-}}|c_{m}(y)|}{S_{-}^{n}}\leq \frac{M}{%
R_{-}^{m}S_{-}^{n}}.
\end{equation*}%
Thus $f(x,y)=\sum_{m\geq 0}\sum_{n\geq 0}c_{mn}x^{m}y^{n}$ is absolutely
convergent for $|x|<R_{-}$, $|y|<S_{-}$. Hence $c_{n}(x)=\sum_{m\geq
0}c_{mn}x^{m}$ converges for $|x|<R_{-}$ and $f$ is holomorphic in $y$ with
convergent expansion $f(x,y)=\sum_{n\geq 0}c_{n}(x)y^{n}$ for $|y|<S_{-}$.
We may then choose $\delta _{1},\delta _{2}$ sufficiently small to show that
the result follows for all $|x|<R$, $|y|<S$. $\qed$

\begin{proposition}
\label{Prop_hol_omega} Suppose that $\omega ^{(g_{a})}$ is a holomorphic
function of a complex parameter $\mu $ for $\left\vert \mu \right\vert <S$.
Then for $|\epsilon |<r_{1}r_{2}$, $\omega ^{(g_{1}+g_{2})}$ is also
holomorphic in $\mu $ for $\left\vert \mu \right\vert <S$.
\end{proposition}

\noindent \textbf{Proof.} Suppose that $\omega ^{(g_{1})}$ is holomorphic in 
$\mu $ wlog. Then $a_{1}(k)$ and $A_{1}(k,l)$ are holomorphic (and
continuous) in $\mu $ for $\left\vert \mu \right\vert <S$. We now show that $%
X_{ab}$ is holomorphic in $\mu $ for $\left\vert \mu \right\vert <S$ using
Lemma \ref{Lem_fxyhol}. Using continuity of $A_{1}$ and (\ref{Xaaself}) of
Proposition \ref{PropXij} we find that for $|\mu +\delta |<S$ 
\begin{equation*}
(I-A_{2}A_{1})\lim_{\delta \rightarrow 0}(X_{11}(\mu +\delta )-X_{11}(\mu
))=0.
\end{equation*}%
But $(I-A_{2}A_{1})$ is invertible for $|\epsilon |<r_{1}r_{2}$ from
Proposition \ref{PropXij} and so $X_{11}(\mu )$ is continuous for $%
\left\vert \mu \right\vert <S$. A similar result holds for $X_{12}$, $X_{21}$
and $X_{22}$. From Theorem \ref{Theoremomg1g2holo}, $\epsilon
^{-(k+l)/2}X_{ab}(k,l)$ is holomorphic in $\epsilon $ for $|\epsilon
|<r_{1}r_{2}$. Furthermore, as explained in Proposition \ref{PropXij}, the $%
\epsilon $ expansion coefficients consist of a finite sum of finite products
of $A_{1}$ and $A_{2}$ terms and thus they are holomorphic in $\mu $ for $%
\left\vert \mu \right\vert <S$. We may therefore apply Lemma \ref{Lem_fxyhol}
to $\epsilon ^{-(k+l)/2}X_{ab}(k,l)$ which is continuous in $\mu $ for $%
\left\vert \mu \right\vert <S$ and holomorphic in $\epsilon $ for $|\epsilon
|<R=r_{1}r_{2}$ with $\epsilon $ expansion coefficients holomorphic in $\mu $
for $\left\vert \mu \right\vert <S$. Thus $\epsilon ^{-(k+l)/2}X_{ab}(k,l)$
and therefore $X_{ab}(k,l)$ is holomorphic in $\mu $ for $\left\vert \mu
\right\vert <S$.

Finally consider $\omega ^{(g_{1}+g_{2})}$ as given in (\ref{omgsisj}) of
Lemma \ref{LemmaaiXijaj}. Using arguments similar to those above we see that 
$\omega ^{(g_{1}+g_{2})}$ is continuous in $\mu $ for $\left\vert \mu
\right\vert <S$ and holomorphic in $\epsilon $ for $|\epsilon |<r_{1}r_{2}$
with $\epsilon $ expansion coefficients holomorphic in $\mu $ for $%
\left\vert \mu \right\vert <S$. Thus $\omega ^{(g_{1}+g_{2})}$ is
holomorphic in $\mu $ for $\left\vert \mu \right\vert <S$. $\qed$

\begin{corollary}
\label{Cor_Omg1g2hol} Given the previous conditions then $\Omega
^{(g_{1}+g_{2})}$ is holomorphic in $\mu $ for $\left\vert \mu \right\vert
<S $ where $|\epsilon |<r_{1}r_{2}$. $\qed$
\end{corollary}

In conclusion, we similarly find by applying Proposition \ref{Prop_hol_omega}
to (\ref{Siberian}) that

\begin{proposition}
\label{Prop_Det_holomorphic} Suppose that $\omega ^{(g_{a})}$ is holomorphic
in $\mu $ for $\left\vert \mu \right\vert <S$. Then $\det (I-A_{1}A_{2})$ is
non-vanishing and holomorphic in $\mu $ for $\left\vert \mu \right\vert <S$
and $|\epsilon |<r_{1}r_{2}$. $\qed$
\end{proposition}

\section{Sewing Two Tori to form a Genus Two Riemann Surface}

\label{Sect_Epsilon_Torii}We now specialize to the case of two tori sewn
together to form a genus two Riemann surface. We first consider an
elementary description of a disk on a torus compatible with $SL(2,\mathbb{Z}%
) $ modular-invariance. We then apply the $\epsilon $ formalism in order to
sew two punctured tori with modular parameters $\tau _{1}$ and $\tau _{2}$
together to form a genus two Riemann surface with period matrix $\Omega
^{(2)}(\tau _{1},\tau _{2},\epsilon )\in \mathbb{H}_{2}$, where $\Omega
^{(2)}$ is holomorphic in $(\tau _{1},\tau _{2},\epsilon )\in \mathcal{D}%
^{\epsilon }$ for a suitably defined domain $\mathcal{D}^{\epsilon }$. We
provide an alternative description of $\Omega ^{(2)}$ in terms of the sum of
weights of particular "necklace" graphs. We then describe the equivariance
properties of this holomorphic mapping from $\mathcal{D}^{\epsilon }$ to $%
\mathbb{H}_{2}$ with respect to a certain subgroup $G\subseteq Sp(4,\mathbb{Z%
})$, and prove that it is invertible in a certain $G$-invariant domain.

\subsection{A Closed Disk on a Torus}

\label{Subsect_closeddisk}

A complex torus $\mathcal{S}$ (that is, a compact Riemann surface of genus $%
1 $), can be represented as a quotient $\mathbb{C}/\Lambda $, where $\Lambda 
$ is a lattice in $\mathbb{C}$. Moreover, two such tori $\mathcal{S}%
_{a},a=1,2$ are isomorphic if, and only if, the lattices are \emph{homothetic%
}, that is, there is a $\xi \in \mathbb{C}^{\ast }$ such that $\Lambda
_{2}=\xi \Lambda _{1}$.

A \emph{framing} of $\mathcal{S}=\mathbb{C}/\Lambda $ is a choice of basis $%
(\sigma ,\varsigma )$ such that the modulus $\tau =\sigma /\varsigma $
satisfies $\tau \in \mathbb{H}_{1}$. We say that the basis $(\sigma
,\varsigma )$ is \emph{positively oriented} in this case. A pair of framed
tori $(\mathbb{C}/\Lambda _{a},\sigma _{a},\varsigma _{a}),a=1,2$ are
isomorphic if, and only if, there is a $\xi $ as above such that $(\sigma
_{2},\varsigma _{2})=\xi (\sigma _{1},\varsigma _{1})$. The modulus $\tau $
depends only on the isomorphism class of the framed torus, and there is a 
\emph{bijection} 
\begin{eqnarray}
\{\mbox{isomorphism classes of framed tori}\} &\rightarrow &\mathbb{H}_{1},
\label{framedtoriisom} \\
(\mathbb{C}/\Lambda ,\sigma ,\varsigma ) &\mapsto &\sigma /\varsigma . 
\notag
\end{eqnarray}

$SL(2,\mathbb{Z})$ is the group of automorphisms of $\Lambda $ which
preserves oriented bases. It acts on isomorphism classes of framed tori via 
\begin{equation*}
\left( 
\begin{array}{cc}
a & b \\ 
c & d%
\end{array}%
\right) :(\mathbb{C}/\Lambda ,\sigma ,\varsigma )\mapsto (\mathbb{C}/\Lambda
,a\sigma +b\varsigma ,c\sigma +d\varsigma ),
\end{equation*}%
and on $\mathbb{H}_{1}$ via fractional linear transformations 
\begin{equation*}
\left( 
\begin{array}{cc}
a & b \\ 
c & d%
\end{array}%
\right) :\tau \mapsto \frac{a\tau +b}{c\tau +d}.
\end{equation*}%
With respect to these two actions, the bijection (\ref{framedtoriisom}) is $%
SL(2,\mathbb{Z})$-equivariant.

In the following it is convenient to identify $\mathcal{S}$ with the
standard fundamental region for $\Lambda$ determined by the basis, and with
appropriate identifications of boundary. To describe a well-defined disk in $%
\mathcal{S}$, define the minimal length of $\Lambda$ as 
\begin{eqnarray}  \label{minlength}
D(\Lambda) = \min_{0 \not = \lambda \in \Lambda} |\lambda|.
\end{eqnarray}
It obviously satisfies 
\begin{eqnarray}  \label{Dhomothety}
D(\xi \Lambda) = |\xi|D(\Lambda). \ \ \ (\xi \not = 0)
\end{eqnarray}

We may now describe a closed disk on $\mathcal{S}$. The proof follows from
the triangle inequality.

\begin{lemma}
\label{Lemmadisk}\ For $p\in \mathcal{S}$, the points $z\in \mathcal{S}$
satisfying $|z-p|\leq kD(\Lambda )$ define a closed disk centered at $p$
provided $k<\frac{1}{2}$. $\qed$
\end{lemma}

Let $\mathcal{S}$ be a complex torus of modulus $\tau$. Among all homothetic
lattices $\Lambda$ for which define $\mathcal{S} \cong \mathbb{C}/\Lambda,$
it will be convenient to work with the lattice $\Lambda_{\tau}$ which has
basis $(2\pi i\tau ,2\pi i)$. Note that for $\gamma = \left(%
\begin{array}{cc}
a & b \\ 
c & d%
\end{array}%
\right) \in SL(2, \mathbb{Z})$ we have 
\begin{eqnarray}  \label{Dgammaactontau}
D(\Lambda_{\gamma \tau}) = \frac{1}{|c \tau + d|} D(\Lambda_{\tau}).
\end{eqnarray}

\subsection{The Genus Two Period Matrix in the $\protect\epsilon $ Formalism}

We now apply the $\epsilon $-formalism to a pair of tori $\mathcal{S}_{a}=%
\mathbb{C}/\Lambda _{a}$ with local co-ordinates $z_{a}$, where $\Lambda
_{a} $ has oriented basis $(\sigma _{a},\varsigma _{a})$ and $\tau
_{a}=\sigma _{a}/\varsigma _{a}\in \mathbb{H}_{1}$ for $a=1,2$. We shall
sometimes refer to $\mathcal{S}_{1}$ and $\mathcal{S}_{2}$ as the left and
right torus respectively. After Lemma \ref{Lemmadisk} we may consider the
annuli $\mathcal{A}_{a}$ centred at the origin of $\mathcal{S}_{a}$
described above, with outer radius $r_{a}<\frac{1}{2}D(\Lambda _{a})$.
Following the prescription of Subsection \ref{Subsec_Omega_eps_g}, we sew
the two tori by identifying the annuli $\mathcal{A}_{1}$ and $\mathcal{A}%
_{2}\ $via the relation $z_{1}z_{2}=\epsilon $ as in (\ref{pinch}), where $%
|\epsilon |\leq r_{1}r_{2}<\frac{1}{4}D(\Lambda _{1})D(\Lambda _{2})$.

As discussed in Subsection \ref{Subsect_closeddisk} we take $(\sigma
_{a},\varsigma _{a})=(2\pi i\tau _{a},2\pi i)$ and $q_{a}=\exp (2\pi i\tau
_{a})$ for $a=1,2$. Define the domain\footnote{The superscript $\epsilon $ 
merely denotes that we are working in the $\epsilon $-formalism, and should 
not be interpreted as a variable of any kind} 
\begin{equation}
\mathcal{D}^{\epsilon }=\{(\tau _{1},\tau _{2},\epsilon )\in \mathbb{H}_{1}%
\mathbb{\times H}_{1}\mathbb{\times C}\ |\ |\epsilon |<\frac{1}{4}D(\Lambda
_{\tau _{1}})D(\Lambda _{\tau _{2}})\}.  \label{Deps}
\end{equation}

We now explicitly determine the period matrix\footnote{%
The genus two period matrix is also described in \cite{T} without proof and
in a different notation.}.

\begin{theorem}
\label{Theorem_epsperiod}Sewing determines a holomorphic map 
\begin{eqnarray}
F^{\epsilon }:\mathcal{D}^{\epsilon } &\rightarrow &\mathbb{H}_{2},  \notag
\\
(\tau _{1},\tau _{2},\epsilon ) &\mapsto &\Omega ^{(2)}(\tau _{1},\tau
_{2},\epsilon ).  \label{Fepsmap'}
\end{eqnarray}%
Moreover $\Omega ^{(2)}=\Omega ^{(2)}(\tau _{1},\tau _{2},\epsilon )$ is
given by 
\begin{eqnarray}
2\pi i\Omega _{11}^{(2)} &=&2\pi i\tau _{1}+\epsilon
(A_{2}(I-A_{1}A_{2})^{-1})(1,1),  \label{Om11eps} \\
2\pi i\Omega _{22}^{(2)} &=&2\pi i\tau _{2}+\epsilon
(A_{1}(I-A_{2}A_{1})^{-1})(1,1),  \label{Om22eps} \\
2\pi i\Omega _{12}^{(2)} &=&-\epsilon (I-A_{1}A_{2})^{-1}(1,1).
\label{Om12eps}
\end{eqnarray}%
Notation here is as follows: $A_{a}(\tau _{a},\epsilon )$ is the infinite
matrix with $(k,l)$-entry 
\begin{equation}
A_{a}(k,l,\tau _{a},\epsilon )=\frac{\epsilon ^{(k+l)/2}}{\sqrt{kl}}%
C(k,l,\tau _{a});  \label{Aki1}
\end{equation}%
$(1,1)$ refers to the $(1,1)$-entry of a matrix with $C(k,l,\tau )$ of (\ref%
{Ckldef}).
\end{theorem}

\noindent \textbf{Proof.} \ The bilinear two form $\omega ^{(1)}$ is given
by (\ref{omega1}). Using (\ref{P1Pnexpansion}), the basis of $1$-forms (\ref%
{akdef}) with periods $(2\pi i\tau _{a},2\pi i)$ is then given by 
\begin{eqnarray}
a_{a}(k,x) &=&\frac{\epsilon ^{k/2}dx}{2\pi i\sqrt{k}}\oint_{\mathcal{C}%
_{a}(z)}z^{-k}P_{2}(\tau _{a},x-z)dz,  \notag \\
&=&\sqrt{k}\epsilon ^{k/2}P_{k+1}(\tau _{a},x)dx.  \label{akdef1}
\end{eqnarray}%
Now (\ref{Aki1}) follows from (\ref{Adef}), (\ref{Pkexpansion}) and (\ref%
{Ckldef}). Note that (\ref{P1quasiperiod}) implies that $\alpha _{a}(k)$ of (%
\ref{alphadef}) is 
\begin{equation}
\alpha _{a}(k)=\epsilon ^{1/2}\delta _{k,1}.  \label{akbperiod}
\end{equation}%
We therefore find $\Omega ^{(2)}$ to be given by (\ref{Om11eps})-(\ref%
{Om12eps}) for $|\epsilon |<r_{1}r_{2}<\frac{1}{4}D(\Lambda _{\tau
_{1}})D(\Lambda _{\tau _{2}})$.

By Theorem \ref{Theoremperiodgeps}, $\Omega ^{(2)}$ is holomorphic in $%
\epsilon $ for $|\epsilon |<r_{1}r_{2}<\frac{1}{4}D(\Lambda _{\tau
_{1}})D(\Lambda _{\tau _{2}})$. The left torus bilinear form $\omega
^{(1)}(x,y,\tau _{1})$ is holomorphic in some neighborhood $|\tau _{1}-\tau
_{1}^{0}|<S$ of any point $\tau _{1}^{0}\in \mathbb{H}_{1}$. Hence we may
apply Corollary \ \ref{Cor_Omg1g2hol} with $\mu =\tau _{1}-\tau _{1}^{0}$
for $|\mu |<S$ to conclude that $\Omega ^{(2)}$ is holomorphic in $\tau _{1}$%
. Similarly $\Omega ^{(2)}$ is holomorphic in $\tau _{2}$, and by Hartog's
theorem (e.g. \cite{Gu}) $\Omega ^{(2)}$ is holomorphic on $\mathcal{D}%
^{\epsilon }$. $\qed$

The infinite matrices $A_{a}(\tau_{a},\epsilon )$ will play a crucial r\^{o}%
le in the further analysis of the $\epsilon $ formalism. Dropping the
subscript, they are symmetric and have the form:

\begin{equation*}
A(\tau ,\epsilon )=\left( 
\begin{array}{ccccc}
\epsilon E_{2}(\tau ) & 0 & \sqrt{3}\epsilon ^{2}E_{4}(\tau ) & 0 & \cdots
\\ 
0 & -3\epsilon ^{2}E_{4}(\tau ) & 0 & -5\sqrt{2}\epsilon ^{3}E_{6}(\tau ) & 
\cdots \\ 
\sqrt{3}\epsilon ^{2}E_{4}(\tau ) & 0 & 10\epsilon ^{3}E_{6}(\tau ) & 0 & 
\cdots \\ 
0 & -5\sqrt{2}\epsilon ^{3}E_{6}(\tau ) & 0 & -35\epsilon ^{4}E_{8}(\tau ) & 
\cdots \\ 
\vdots & \vdots & \vdots & \vdots & \ddots%
\end{array}%
\right)
\end{equation*}

\subsection{Chequered Necklace Expansion for $\Omega ^{(2)}$}

It is useful to introduce an interpretation for the expressions for $\Omega
^{(2)}$ found above in terms of the sum of weights of certain graphs. Let us
introduce the set of \emph{chequered necklaces} $\mathcal{N}$. By
definition, these are connected graphs with $n\geq 2$ nodes, $(n-2)$ of
which have valency $2$ and two of which have valency $1$ (these latter are
the \emph{end nodes}), together with an orientation, say from left to right,
on the edges. Moreover vertices are labelled by positive integers and edges
are labelled alternatively by $1$ or $2$ as one moves along the graph e.g.

\begin{equation*}
\ \overset{k_{1}}{\bullet }\overset{1}{\longrightarrow }\overset{k_{2}}{%
\bullet }\overset{2}{\longrightarrow }\overset{k_{3}}{\bullet }\overset{1}{%
\longrightarrow }\overset{k_{4}}{\bullet }\overset{2}{\longrightarrow }%
\overset{k_{5}}{\bullet }\overset{1}{\longrightarrow }\overset{k_{6}}{%
\bullet }
\end{equation*}%
We also define the \emph{degenerate necklace} $N_{0}$ to be a single node
with no edges. Define a weight function 
\begin{equation*}
\omega :\mathcal{N}\longrightarrow \mathbb{C}[E_{2}(\tau_{a}),E_{4}(\tau
_{a}),E_{6}(\tau_{a}),\epsilon \ | \ a = 1,2 ],
\end{equation*}%
as follows: if a chequered necklace $N$ has edges $E$ labelled as \ $\overset%
{k}{\bullet }\overset{a}{\longrightarrow }\overset{l}{\bullet }$ \ then we
define 
\begin{eqnarray}
\omega (E) &=&A_{a}(k,l,\tau_{a},\epsilon ),  \notag \\
\omega (N) &=&\prod \omega (E),  \label{wtedge}
\end{eqnarray}%
where $A_{a}(k,l,\tau_{a},\epsilon)$ is given by (\ref{Aki1}) and the
product is taken over all edges $E$ of $N$. We further define $%
\omega(N_{0})=1$.

Among all chequered necklaces there is a distinguished set for which both
end nodes are labelled by $1$. There are four types of such chequered
necklaces, which may be further distinguished by the labels of the two edges
at the extreme left and right. We use the notation (\ref{bardef}) for $a=1,2$%
, and say that the chequered necklace 
\begin{equation*}
\overset{1}{\bullet }\overset{a}{\longrightarrow }\overset{i}{\bullet }\ldots%
\overset{j}{\bullet }\overset{{b}}{\longrightarrow }\overset{1}{\bullet }%
\hspace{10mm}
\end{equation*}%
is of \emph{type $\overline{a}\overline{b}$ }. We then set 
\begin{eqnarray*}
\mathcal{N}_{ab} &=&\{%
\mbox{isomorphism classes of chequered
necklaces of type}\ ab\}, \\
\omega_{ab} &=&\sum_{N\in \mathcal{N}_{ab}}\omega (N),
\end{eqnarray*}%
where $\omega_{ab}$ is considered as an element in $\mathbb{C}%
[E_{2}(\tau_{a}),E_{4}(\tau_{a}),E_{6}(\tau_{a}),\epsilon \ |\ a=1,2]$. It
is clear that we may use this formalism to represent matrix expressions like
those appearing earlier. Then we have

\begin{lemma}
\label{Lemma_cheq}\noindent\ For $a=1,2$ we have 
\begin{eqnarray*}
\omega _{a\bar{a}} &=&\omega _{\bar{a}a}=(I-A_{a}A_{\bar{a}})^{-1}(1,1) \\
\omega _{aa} &=&(A_{\bar{a}}(I-A_{a}A_{\bar{a}})^{-1})(1,1).\qed
\end{eqnarray*}
\end{lemma}

Thus we may conclude from Theorem \ref{Theorem_epsperiod} that

\begin{proposition}
\label{Propepsperiodgraph} For $a=1,2$ we have 
\begin{eqnarray*}
\Omega _{aa}^{(2)} &=&\tau _{a}+\frac{\epsilon }{2\pi i}\omega _{aa}, \\
\Omega _{a\bar{a}}^{(2)} &=&-\frac{\epsilon }{2\pi i}\omega _{a\bar{a}}.\ \ %
\qed
\end{eqnarray*}
\end{proposition}

\subsection{Equivariance of $F^{\protect\epsilon}$}

\label{Subsec_Feps}In Theorem \ref{Theorem_epsperiod} we established the
existence of the analytic map $F^{\epsilon }$. Here we establish the
equivariance of this map with respect to a certain subgroup $G$ of $Sp(4,%
\mathbb{Z})$. We will employ the graphical representation for $\Omega
=\Omega ^{(2)}(\tau _{1},\tau _{2},\epsilon )$ in terms of chequered
necklaces discussed in the last subsection.

As an abstract group, $G$ is isomorphic to $(SL(2,\mathbb{Z})$ $\times SL(2,%
\mathbb{Z}))\rtimes \mathbb{Z}_{2}$, i.e. the direct product of two copies
of $SL(2,\mathbb{Z})$ which are interchanged upon conjugation by an
involution. There is a natural injection $G\rightarrow Sp(4,\mathbb{Z})$ in
which the two $SL(2,\mathbb{Z})$ subgroups are mapped to 
\begin{equation}
\Gamma_{1}=\left\{ \left[ 
\begin{array}{cccc}
a_{1} & 0 & b_{1} & 0 \\ 
0 & 1 & 0 & 0 \\ 
c_{1} & 0 & d_{1} & 0 \\ 
0 & 0 & 0 & 1%
\end{array}%
\right] \right\} ,\;\Gamma_{2}=\left\{ \left[ 
\begin{array}{cccc}
1 & 0 & 0 & 0 \\ 
0 & a_{2} & 0 & b_{2} \\ 
0 & 0 & 1 & 0 \\ 
0 & c_{2} & 0 & d_{2}%
\end{array}%
\right] \right\} ,  \label{Gamma1Gamma2}
\end{equation}%
and the involution is mapped to 
\begin{equation}
\beta =\left[ 
\begin{array}{cccc}
0 & 1 & 0 & 0 \\ 
1 & 0 & 0 & 0 \\ 
0 & 0 & 0 & 1 \\ 
0 & 0 & 1 & 0%
\end{array}%
\right] .  \label{betagen}
\end{equation}
In this way we obtain a natural action of $G$ on $\mathbb{H}_{2}$. The
action on $\mathcal{D}^{\epsilon}$ is described in the next Lemma.

\begin{lemma}
\label{LemGDeps} $G$ has a left action on $\mathcal{D}^{\epsilon }$ as
follows: 
\begin{eqnarray}
\gamma_{1}.(\tau_{1},\tau_{2},\epsilon ) &=&(\gamma_{1}\tau_{1},\tau _{2},%
\frac{{\epsilon }}{c_{1}\tau_{1}+d_{1}}),  \label{gam1eps} \\
\gamma_{2}.(\tau_{1},\tau_{2},\epsilon ) &=&(\tau_{1},\gamma_{2}\tau _{2},%
\frac{{\epsilon }}{c_{2}\tau_{2}+d_{2}}),  \label{gam2eps} \\
\beta .(\tau_{1},\tau_{2},\epsilon ) &=&(\tau_{2},\tau_{1},{\epsilon }),
\label{betaeps}
\end{eqnarray}%
for $(\gamma_{1},\gamma_{2})\in SL(2,\mathbb{Z})\times SL(2,\mathbb{Z})$.
\end{lemma}

\noindent \textbf{Proof. } It is straightforward (and quite standard) to see
that (\ref{gam1eps}) - (\ref{betaeps}) formally define a (left) action of $G$
on $\mathbb{H}_{1}\times \mathbb{H}_{1}\times \mathbb{C}$. What we must show
is that this action preserves the domain $\mathcal{D}^{\epsilon }$. For $%
\beta $ this is obvious, and for elements $(\gamma _{1},\gamma _{2})$ it
follows from (\ref{Dgammaactontau}). \ \ \ \ \ \ $\qed$

We now establish the following result:

\begin{theorem}
\label{TheoremGequiv} $F^{\epsilon }$ is equivariant with respect to the
action of $G,$ i.e. there is a commutative diagram for $\gamma \in G$, 
\begin{equation*}
\begin{array}{ccc}
\mathcal{D}^{\epsilon } & \overset{F^{\epsilon }}{\rightarrow } & \mathbb{H}%
_{2} \\ 
\gamma \downarrow &  & \downarrow \gamma \\ 
\mathcal{D}^{\epsilon } & \overset{F^{\epsilon }}{\rightarrow } & \mathbb{H}%
_{2}%
\end{array}%
\end{equation*}
\end{theorem}

\noindent \textbf{Proof. } Fix $(\tau _{1},\tau _{2},\epsilon )\in \mathcal{D%
}^{\epsilon }$, with $\Omega =F^{\epsilon }(\tau _{1},\tau _{2},\epsilon
)=\left( 
\begin{array}{cc}
\Omega _{11} & \Omega _{12} \\ 
\Omega _{12} & \Omega _{22}%
\end{array}%
\right) $. Of course, each $\Omega _{ij}$ is a function of $(\tau _{1},\tau
_{2},\epsilon )$. The action of $G$ on $\mathbb{H}_{2}$ is given in (\ref%
{eq: modtrans}), and in particular $\beta :\Omega \mapsto \left( 
\begin{array}{cc}
\Omega _{22} & \Omega _{12} \\ 
\Omega _{12} & \Omega _{11}%
\end{array}%
\right) $. Therefore from (\ref{betaeps}) we have 
\begin{equation*}
F^{\epsilon }(\beta (\tau _{1},\tau _{2},\epsilon ))=F^{\epsilon }(\tau
_{2},\tau _{1},\epsilon )=\left( 
\begin{array}{cc}
\Omega _{22} & \Omega _{12} \\ 
\Omega _{12} & \Omega _{11}%
\end{array}%
\right) =\beta (F^{\epsilon }(\tau _{1},\tau _{2},\epsilon )).
\end{equation*}%
So the Theorem is true in case $\gamma =\beta $. To complete the proof of
the Theorem, it suffices to consider the case when $\gamma =\gamma _{1}$
lies in the `left' modular group acting on $\tau _{1}$. From (\ref{eq:
modtrans}), 
\begin{equation}
\gamma _{1}:\Omega \mapsto \left( 
\begin{array}{cc}
\frac{a_{1}\Omega _{11}+b_{1}}{c_{1}\Omega _{11}+d_{1}} & \frac{\Omega _{12}%
}{c_{1}\Omega _{11}+d_{1}} \\ 
\frac{\Omega _{12}}{c_{1}\Omega _{11}+d_{1}} & \Omega _{22}-\frac{%
c_{1}\Omega _{12}^{2}}{c_{1}\Omega _{11}+d_{1}}%
\end{array}%
\right) ,  \label{gamma1Omega}
\end{equation}%
and we are obliged to show that the matrix in display (\ref{gamma1Omega})
coincides with $F^{\epsilon }(\gamma _{1}(\tau _{1},\tau _{2},\epsilon
))=F^{\epsilon }(\gamma _{1}\tau _{1},\tau _{2},\frac{{\epsilon }}{c_{1}\tau
_{1}+d_{1}})$. In other words, we must establish the following identities:

\begin{eqnarray}
\frac{a_{1}\Omega_{11}+b_{1}}{c_{1}\Omega_{11}+d_{1}} &=&\Omega
_{11}(\gamma_{1}\tau_{1},\tau_{2},\frac{{\epsilon }}{c_{1}\tau_{1}+d_{1}}),
\label{gamma1action1} \\
\frac{\Omega_{12}}{c_{1}\Omega_{11}+d_{1}} &=&\Omega_{12}(\gamma_{1}\tau
_{1},\tau_{2},\frac{{\epsilon }}{c_{1}\tau_{1}+d_{1}}),
\label{gamma1action2} \\
\Omega_{22}-\frac{c_{1}\Omega_{11}^{2}}{c_{1}\Omega_{11}+d_{1}} &=&\Omega
_{22}(\gamma_{1}\tau_{1},\tau_{2},\frac{{\epsilon }}{c_{1}\tau_{1}+d_{1}}).
\label{gamma1action3}
\end{eqnarray}%
$A_{a}(k,l,\tau_{a},\epsilon )$ of (\ref{Aki1}) is a modular form of weight $%
k+l$ for $k+l>2$, whereas $A_{a}(1,1,\tau_{a},\epsilon )=\epsilon
E_{2}(\tau_{a})$ enjoys an exceptional transformation law thanks to (\ref%
{gammaE2}).

Using Lemma \ref{LemGDeps} we then find that 
\begin{eqnarray}
A_{1}(k,l,\gamma_{1}\tau_{1},\frac{{\epsilon }}{c_{1}\tau_{1}+d_{1}}))
&=&(c_{1}\tau_{1}+d_{1})^{(k+l)/2}(A_{1}(\tau_{1},\epsilon )+\kappa \delta
_{k1}\delta_{l1}),  \label{gamma1A1} \\
A_{2}(k,l,\tau_{2},\frac{{\epsilon }}{c_{1}\tau_{1}+d_{1}})) &=&(c_{1}\tau
_{1}+d_{1})^{-(k+l)/2}A_{2}(\tau_{2},\epsilon ),  \label{gamma1A2}
\end{eqnarray}%
where 
\begin{equation}
\kappa =-\frac{\epsilon }{2\pi i}\frac{c_{1}}{c_{1}\tau_{1}+d_{1}}.
\label{kappadef}
\end{equation}%
It follows from Proposition \ref{Propepsperiodgraph} both that

\begin{equation}
1-\kappa \omega _{11}=\frac{c_{1}\Omega _{11}+d_{1}}{c_{1}\tau _{1}+d_{1}},
\label{omegakappa}
\end{equation}%
and 
\begin{eqnarray}
&&\hspace{1.5cm}\Omega _{11}(\gamma _{1}\tau _{1},\tau _{2},\frac{\epsilon }{%
c_{1}\tau _{1}+d_{1}})  \notag \\
&=&\frac{1}{c_{1}\tau _{1}+d_{1}}(a_{1}\tau _{1}+b_{1}+\frac{\epsilon }{2\pi
i}\omega _{11}(\gamma _{1}\tau _{1},\tau _{2},\frac{\epsilon }{c_{1}\tau
_{1}+d_{1}})).  \label{nextOmega11}
\end{eqnarray}%
Consider a necklace $N\in \mathcal{N}_{11}$ of weight $\omega (N)$ and let $%
\mathcal{S}_{11}(N)$ denote the set of all "broken" graphs formed from $N$
by deleting any $n$ edges of type $\overset{1}{\bullet }\overset{1}{%
\longrightarrow }\overset{1}{\bullet }$ for all $n\geq 0$. Every such graph
consists of $n+1$ connected graphs $N_{1},\ldots N_{n+1}$ of type $11$. From
(\ref{gamma1A1}) and (\ref{gamma1A2}) it is therefore follows that 
\begin{equation*}
\omega (N)(\gamma _{1}\tau _{1},\tau _{2},\frac{\epsilon }{c_{1}\tau
_{1}+d_{1}})=\frac{1}{c_{1}\tau _{1}+d_{1}}\sum_{n\geq 0}\kappa
^{n}\sum_{N_{1},\ldots N_{n+1}}\omega (N_{1})\ldots \omega (N_{n+1}).
\end{equation*}%
Summing over all $N$ we then find 
\begin{eqnarray*}
\omega _{11}(\gamma _{1}\tau _{1},\tau _{2},\frac{\epsilon }{c_{1}\tau
_{1}+d_{1}}) &=&\frac{1}{(c_{1}\tau _{1}+d_{1})}\sum_{n\geq 0}\kappa
^{n}\omega _{11}^{n+1} \\
&=&\frac{1}{(c_{1}\tau _{1}+d_{1})}\frac{\omega _{11}}{1-\kappa \omega _{11}}
\\
&=&\frac{\omega _{11}}{c_{1}\Omega _{11}+d_{1}},
\end{eqnarray*}%
where for the last equality we used (\ref{omegakappa}). Now (\ref%
{nextOmega11}) yields 
\begin{eqnarray*}
\Omega _{11}(\gamma _{1}\tau _{1},\tau _{2},\frac{\epsilon }{c_{1}\tau
_{1}+d_{1}}) &=&\frac{1}{c_{1}\tau +d_{1}}(a_{1}\tau _{1}+b_{1}+\frac{\Omega
_{11}-\tau _{1}}{c_{1}\Omega _{11}+d_{1}}) \\
&=&\frac{a_{1}\Omega _{11}+b_{1}}{c_{1}\Omega _{11}+d_{1}},
\end{eqnarray*}%
which is the desired (\ref{gamma1action1}). Similarly from Proposition \ref%
{Propepsperiodgraph} we have 
\begin{equation*}
\Omega _{12}(\gamma _{1}\tau _{1},\tau _{2},\frac{\epsilon }{c_{1}\tau
_{1}+d_{1}})=-\frac{1}{(c_{1}\tau _{1}+d_{1})}\frac{\epsilon }{2\pi i}\omega
_{12}(\gamma _{1}\tau _{1},\tau _{2},\frac{\epsilon }{c_{1}\tau _{1}+d_{1}}).
\end{equation*}%
Breaking necklaces of type $12$ results in products over necklaces of type $%
11$ together with one necklace of type $12$. Hence by a similar argument to
that above we find 
\begin{eqnarray*}
\omega _{12}(\gamma _{1}\tau _{1},\tau _{2},\frac{\epsilon }{c_{1}\tau
_{1}+d_{1}}) &=&\frac{\omega _{12}}{1-\kappa \omega _{11}} \\
&=&\frac{(c_{1}\tau _{1}+d_{1})\omega _{12}}{c_{1}\Omega _{11}+d_{1}},
\end{eqnarray*}%
so that $\Omega _{12}(\gamma _{1}\tau _{1},\tau _{2},\frac{\epsilon }{c\tau
_{1}+d_{1}})$ is as in (\ref{gamma1action2}). Finally, 
\begin{equation*}
\Omega _{22}(\gamma _{1}\tau _{1},\tau _{2},\frac{\epsilon }{c\tau _{1}+d_{1}%
})=\tau _{2}+\frac{1}{c_{1}\tau _{1}+d_{1}}\frac{\epsilon }{2\pi i}\omega
_{22}(\gamma _{1}\tau _{1},\tau _{2},\frac{\epsilon }{c\tau _{1}+d_{1}}).
\end{equation*}%
Breaking necklaces of type $22$ results in products over necklaces of type $%
11$ together with one necklace of type $12$ and another of type $21$. Hence
by a similar argument to that above we find 
\begin{eqnarray*}
&&\frac{1}{(c_{1}\tau _{1}+d_{1})}\frac{\epsilon }{2\pi i}\omega
_{22}(\gamma _{1}\tau _{1},\tau _{2},\frac{\epsilon }{c\tau _{1}+d_{1}}) \\
&=&\frac{\epsilon }{2\pi i}(\omega _{22}+\frac{\kappa \omega _{12}^{2}}{%
1-\kappa \omega _{11}})=\Omega _{22}-\tau _{2}-\frac{c_{1}\Omega _{12}^{2}}{%
c_{1}\Omega _{11}+d_{1}},
\end{eqnarray*}%
leading to (\ref{gamma1action3}). This completes the proof of the Theorem. $%
\qed$

\subsection{Local Invertibility of $F^{\protect\epsilon}$ about the Two Tori
Degeneration Point $\protect\epsilon =0$}

\label{Subsec_Feps_degen}Let $\mathcal{D}_{0}^{\epsilon }$ be the subset of $%
\mathcal{D}^{\epsilon }$ for which $\epsilon =0$. From Theorem \ref%
{Theorem_epsperiod} it is clear that the restriction of $F^{\epsilon }$ to $%
\mathcal{D}_{0}^{\epsilon }$ induces the natural identification 
\begin{eqnarray}
F^{\epsilon }:\mathcal{D}_{0}^{\epsilon } &\overset{\sim }{\rightarrow }&%
\mathbb{H}_{1}\times \mathbb{H}_{1}\subseteq \mathbb{H}_{2}  \notag \\
(\tau _{1},\tau _{2},0) &\mapsto &\left( 
\begin{array}{cc}
\tau _{1} & 0 \\ 
0 & \tau _{2}%
\end{array}%
\right) .  \label{epsdegen}
\end{eqnarray}%
$\mathcal{D}_{0}^{\epsilon }$ corresponds to the set of points where the
genus $2$ Riemann surface degenerates into a pair of genus $1$ surfaces with 
$\Omega ^{(2)}=\mathrm{diag}(\Omega _{11}^{(2)},\Omega _{22}^{(2)})$. We
will consider the invertibility of the map $F^{\epsilon }$ in a neighborhood
of a point in $\mathcal{D}_{0}^{\epsilon }$. First we prepare a Lemma.

Recall (e.g. \cite{FK2}) that a group $H$ of homeomorphisms of a space $X$
is said to act \emph{discontinuously} on $X$ if each point $x\in X$ has a 
\emph{precisely invariant open neighborhood} under the action of $H$ in the
following sense (loc. cit.): the stabilizer Stab$(x)$ of $x$ in $H$ is
finite, and there is an open neighborhood $N$ of $x$ such that $hN\cap
N=\phi $ if $h\notin $ Stab$(x)$ and $hN=N$ if $h\in $ Stab$(x)$.

\begin{lemma}
\label{discontinuity} Suppose that $H$ acts discontinuously on a pair of
spaces $X,Y$, and that $F:X\rightarrow Y$ is a continuous $H$-equivariant
map. Then the following hold:\newline
(a) If $x\in X$ and $F(x)=y$ then there are precisely invariant open
neighborhoods (under the action of $H$) $U\subseteq X$ and $V\subseteq Y$ of 
$x$ and $y$ respectively with $F(U)\subseteq V$; \newline
(b) Suppose further that Stab$(x)=$ Stab$(y)$ and that the restriction of $F$
to $U$ is $1-1$. Then $F$ is $1-1$ on the $H$-invariant domain $%
\bigcup_{h\in H}hU$.
\end{lemma}

\noindent \textbf{Proof.} For part (a), let $V$ be a precisely invariant
open neighborhood of $y$ in $Y$, $U^{\prime}$ a precisely invariant
neighborhood of $x$ in $X$, and set $U = F^{-1}(V) \cap U^{\prime}$. Because 
$F$ is $H$-equivariant then Stab$(x) \subseteq$ Stab$(y)$, and from this it
follows that $U$ is also precisely invariant under the action of $H$. Now
(a) follows.

As for (b), suppose the contrary so that there exist $u_{1},u_{2}\in U$ and $%
h_{1},h_{2}\in H$ such that $h_{1}u_{1}\neq h_{2}u_{2}$ and $%
F(h_{1}u_{1})=F(h_{2}u_{2})$. Thanks to the equivariance of $F$ it is no
loss to assume that $h_{2}=1$, so that $h_{1}u_{1}\neq u_{2}$ and $%
F(h_{1}u_{1})=F(u_{2})$. From the last equality we see that $h_{1}V\cap
V\neq \phi $, so that $h_{1}V=V$ and $h_{1}\in $ Stab$(y)$.

Therefore, $h_{1}\in $ Stab$(x)$ by hypothesis, and therefore $h_{1}U=U$.
But then $h_{1}u_{1}$ and $u_{2}\in U$ are distinct points of $U$ on which $%
F $ takes the same value. This contradicts the assumption that $F$ is $1-1$
on $U$, and completes the proof of the Lemma. $\qed$

We now have

\begin{proposition}
\label{Prop_Fepsinverse} Let $x\in \mathcal{D}_{0}^{\epsilon }$. Then there
exists a $G-$invariant neighborhood $\mathcal{N}_{x}^{\epsilon }\subseteq 
\mathcal{D}^{\epsilon }$ of $x$ throughout which $F^{\epsilon }$ is
invertible.
\end{proposition}

\noindent \textbf{Proof.} Let $x=(\tau _{1},\tau _{2},0)$. From Theorem \ref%
{Theorem_epsperiod}, the Jacobian of $F^{\epsilon }$ at $x$ satisfies 
\begin{equation*}
\left\vert \frac{\partial (\Omega _{11},\Omega _{22},\Omega _{12})}{\partial
(\tau _{1},\tau _{2},\epsilon )}\right\vert _{x}=\left\vert 
\begin{array}{ccc}
1 & 0 & 0 \\ 
0 & 1 & 0 \\ 
0 & 0 & -1%
\end{array}%
\right\vert =-1.
\end{equation*}%
By the inverse function theorem, there exists an open neighborhood of $x$ in 
$\mathcal{D}^{\epsilon }$ throughout which $F^{\epsilon }$ is invertible.
Set $F^{\epsilon }(x)=y$. It follows immediately from (\ref{epsdegen}) that
the stabilizers (in $G$) of $x$ and $y$ are equal.

Next, it is well-known that the action (\ref{eq: modtrans}) of $Sp(2g,%
\mathbb{Z})$ on $\mathbb{H}_{g}$ is discontinuous. In particular, the action
of $G$ on $\mathbb{H}_{2}$ is discontinuous, furthermore from the case $g=1$
it is easy to see that the action of $G$ on $\mathbb{H}_{1}\times \mathbb{H}%
_{1}\times \mathbb{C}$ (and hence also on $D^{\epsilon }$) is also
discontinuous. Choose precisely invariant neighborhoods (under the action of 
$G$) $U,V$ of $x$, respectively $y$ such that the conditions of part (a) of
Lemma \ref{discontinuity} hold. It is clear from the non-vanishing of the
Jacobian that we may also assume that $F^{\epsilon }$ is $1-1$ on $U$. Thus,
we have achieved the hypotheses of part (b) of Lemma \ref{discontinuity}.
That result tells us that the open neighborhood 
\begin{equation*}
\mathcal{N}_{x}^{\epsilon }=\bigcup_{\gamma \in G}\gamma U
\end{equation*}%
has the desired properties. $\qed$

We conclude this section with the explicit form of $\Omega =\Omega
^{(2)}(\tau_{1},\tau_{2},\epsilon )$ to order $\epsilon ^{3}$. We have from (%
\ref{Om11eps}) to (\ref{Om12eps}) that

\begin{eqnarray}
2\pi i\Omega_{11} &=&2\pi i\tau_{1}+\epsilon ^{2}E_{2}(\tau _{2})+O(\epsilon
^{4}),  \label{Om11eps_lead} \\
2\pi i\Omega_{22} &=&2\pi i\tau_{2}+\epsilon ^{2}E_{2}(\tau _{1})+O(\epsilon
^{4}).  \label{Om12eps_lead} \\
2\pi i\Omega_{12} &=&-\epsilon (1+\epsilon ^{2}E_{2}(\tau_{1})E_{2}(\tau
_{2})+O(\epsilon ^{4})).  \label{Om22eps_lead}
\end{eqnarray}%
It is straightforward to check the equivariance properties described in
Theorem \ref{TheoremGequiv} to the given order. In Appendix A more detailed
expansions are provided. We may invert this relationship using Proposition %
\ref{Prop_Fepsinverse} to find to order $\Omega_{12}^{3}$ that

\begin{eqnarray}
\tau_{1} &=&\Omega_{11}-2\pi i\Omega_{12}^{2}E_{2}(\Omega_{22})+O(\Omega
_{12}^{4}),  \label{tau1lead} \\
\tau_{2} &=&\Omega_{22}-2\pi i\Omega_{12}^{2}E_{2}(\Omega_{11})+O(\Omega
_{12}^{4}),  \label{tau2lead} \\
\epsilon &=&-2\pi i\Omega_{12}(1-(2\pi i)^{2}\Omega_{12}^{2} E_{2}(\Omega
_{11})E_{2}(\Omega_{22})+O(\Omega_{12}^{4})).  \label{epslead}
\end{eqnarray}

\section{The $\protect\rho $ Formalism for Self-Sewing a Riemann Surface}

\label{Sect_Rho_g}

\subsection{The General $\protect\rho $ Formalism}

In this section we review the general Yamada construction \cite{Y} for
sewing a Riemann surface of genus $g$ to itself to form a surface of genus $%
g+1$. We consider examples of sewing a Riemann sphere to itself in some
detail where the Catalan series arise in a surprising way. In the next
section, this general formalism will be applied to the construction of a
genus two surface where the Catalan series again plays an important role.

Consider a Riemann surface $\mathcal{S}$ of genus $g$ and let $z_{1},z_{2}$
be local coordinates on $\mathcal{S}$ in the neighborhood of two separated
points $p_{1}$ and $p_{2}$. Consider two disks $\left\vert z_{a}\right\vert
\leq r_{a}$ $\,$for $r_{a}>0$ sufficiently small and $a=1,2$. Note that $%
r_{1},r_{2}$ must be sufficiently small to also ensure that the disks do not
intersect. Introduce a complex parameter ${\rho }$ where $|{\rho }|\leq
r_{1}r_{2}$ and excise the disks%
\begin{equation*}
\{z_{a},\left\vert z_{a}\right\vert <|\rho |r_{\bar{a}}^{-1}\}\subset 
\mathcal{S}
\end{equation*}%
for $a=1,2$ to form a twice-punctured surface 
\begin{equation*}
\hat{\mathcal{S}}=\mathcal{S}\backslash \bigcup_{a=1,2}\{z_{a},\left\vert
z_{a}\right\vert <|\rho |r_{\bar{a}}^{-1}\}.
\end{equation*}%
Here we again use the convention (\ref{bardef}). We define annular regions $%
\mathcal{A}_{a}\subset \hat{\mathcal{S}}$ with $\mathcal{A}_{a}=\{z_{a},|{%
\rho }|r_{\bar{a}}^{-1}\leq \left\vert z_{a}\right\vert \leq r_{a}\}$ and
identify them as a single region $\mathcal{A}=\mathcal{A}_{1}\simeq \mathcal{%
A}_{2}$ via the sewing relation 
\begin{equation}
z_{1}z_{2}=\rho ,  \label{rhosew}
\end{equation}%
to form a compact Riemann surface $\hat{\mathcal{S}}\backslash \{\mathcal{A}%
_{1}\cup \mathcal{A}_{2}\}\cup \mathcal{A}$ of genus $g+1$. The sewing
relation (\ref{rhosew}) can be considered to be a parameterization of a
cylinder connecting the punctured Riemann surface to itself. Using the
Yamada formalism \cite{Y}, and noting the notational differences, the genus $%
g+1$ normalized differential of the second kind $\omega ^{(g+1)}\,$of (\ref%
{omegag}) obeys

\begin{theorem}[Ref. \protect\cite{Y}, Theorem 1, Theorem 4]
\label{Theoremomgp1holo}

\begin{description}
\item[(a)] $\omega ^{(g+1)}$ is holomorphic in $\rho $ for $|\rho
|<r_{1}r_{2}$.

\item[(b)] $\lim_{\rho \rightarrow 0}\omega ^{(g+1)}(x,y)=\omega ^{(g)}(x,y)$
for $x,y\in \hat{\mathcal{S}}$. $\qed$
\end{description}
\end{theorem}

Regarded as a power series in $\rho $, the coefficients of the analytic
expansion of $\omega ^{(g+1)}$ in $\rho $ can be calculated from $\omega
^{(g)}$. Let $\mathcal{C}_{a}(z_{a})\subset \mathcal{A}_{a}$ denote a closed
anti-clockwise oriented contour parameterized by $z_{a}$ surrounding the
puncture at $z_{a}=0$ on $\hat{\mathcal{S}}$. Note that $\mathcal{C}%
_{1}(z_{1})$ may be deformed to $-\mathcal{C}_{2}(z_{2})$. Then similarly to
Lemma \ref{Lemma_inteqn1} we find \cite{Y}

\begin{lemma}
\label{Lemma_inteqn2} 
\begin{equation}
\omega ^{(g+1)}(x,y)=\omega ^{(g)}(x,y)+\frac{1}{2\pi i}\sum_{a=1,2}\oint_{%
\mathcal{C}_{a}(z)}(\omega ^{(g+1)}(y,z)\text{ }\int^{z}\omega
^{(g)}(x,\cdot )),  \label{w2rhointeqn}
\end{equation}%
for $x,y\in \hat{\mathcal{S}}$. $\qed$
\end{lemma}

For $a,b=1,2\,$\ and $k,l=1,2,\ldots $ we define weighted moments 
\begin{equation}
Y_{\bar{a}b}(k,l)=\frac{\rho ^{(k+l)/2}}{\sqrt{kl}}\frac{1}{(2\pi i)^{2}}%
\oint_{\mathcal{C}_{a}(u)}\oint_{\mathcal{C}_{b}(v)}u^{-k}v^{-l}\omega
^{(g+1)}(u,v).  \label{Yijdef}
\end{equation}%
Note that $Y_{ab}(k,l)=Y_{\bar{b}\bar{a}}(l,k)$. We also define $%
Y=(Y_{ab}(k,l))$ to be the infinite matrix indexed by the pairs $a,k\,$and $%
b,l$. We define a set of holomorphic 1-forms on $\hat{\mathcal{S}}$ 
\begin{equation}
a_{a}(k,x)=\frac{\rho ^{k/2}}{2\pi i\sqrt{k}}\oint_{\mathcal{C}%
_{a}(z_{a})}z_{a}^{-k}\omega ^{(g)}(x,z_{a}),  \label{akdef2}
\end{equation}%
and define $a(x)=(a_{a}(k,x))$ and $\bar{a}(x)=(a_{\bar{a}}(k,x))$ to be the
infinite row vectors indexed by $a,k\,$. In a similar way to Lemma \ref%
{LemmaaiXijaj} we then have

\begin{lemma}
\label{Lemma_omg1g} $\omega ^{(g+1)}(x,y)$ for $x,y\in \hat{\mathcal{S}\text{
\ }}$is given by 
\begin{equation}
\omega ^{(g+1)}(x,y)=\omega ^{(g)}(x,y)-a(x)(I-Y)\bar{a}(y)^{T}.\qed
\label{omgplus1}
\end{equation}
\end{lemma}

We next compute the explicit form of $Y$ in terms of the following weighted
moments of $\omega ^{(g)}$ 
\begin{eqnarray}
R_{\bar{a}b}(k,l) &=&-\frac{\rho ^{(k+l)/2}}{\sqrt{kl}}\frac{1}{(2\pi i)^{2}}%
\oint_{\mathcal{C}_{a}(x)}\oint_{\mathcal{C}_{b}(y)}x^{-k}y^{-l}\omega
^{(g)}(x,y)  \notag \\
&=&-\frac{\rho ^{k/2}}{\sqrt{k}}\frac{1}{2\pi i}\oint_{\mathcal{C}%
_{a}(x)}x^{-k}a_{b}(l,x),  \label{Rdef}
\end{eqnarray}%
where $R_{ab}(k,l)=R_{\bar{b}\bar{a}}(l,k)$ and the extra minus sign is
introduced for later convenience. We may consider $R$ as an infinite block
matrix (similar to $Q$ of (\ref{XAQ_def}))%
\begin{equation}
R=(R_{ab}(k,l))=-\left[ 
\begin{array}{ll}
B & A \\ 
A & B^{T}%
\end{array}%
\right] ,  \label{Rformula}
\end{equation}%
with 
\begin{eqnarray}
A(k,l) &=&A(k,l,\rho )=\frac{\rho ^{(k+l)/2}}{(2\pi i)^{2}\sqrt{kl}}\oint_{%
\mathcal{C}_{1}(x)}\oint_{\mathcal{C}_{1}(y)}x^{-k}y^{-l}\omega ^{(g)}(x,y) 
\notag \\
B(k,l) &=&B(k,l,\rho )=\frac{\rho ^{(k+l)/2}}{(2\pi i)^{2}\sqrt{kl}}\oint_{%
\mathcal{C}_{2}(x)}\oint_{\mathcal{C}_{1}(y)}x^{-k}y^{-l}\omega ^{(g)}(x,y).
\label{ABdef}
\end{eqnarray}

Similarly to Proposition \ref{PropXij} we find:

\begin{proposition}
\label{PropYij} $Y_{ab}(k,l)$ is given in terms of $R$ by 
\begin{equation}
I-Y=(I-R)^{-1}.  \label{oneminusY}
\end{equation}%
Here 
\begin{equation*}
(I-R)^{-1}=\sum_{n\geq 0}R^{n}
\end{equation*}%
and is convergent in $\rho $ for $|\rho |<r_{1}r_{2}$. $\qed$
\end{proposition}

Likewise, similarly to Theorem \ref{Theorem_Det} we may define $\det (I-R)$
and find:

\begin{theorem}
\label{Theorem_Det_rho} $\det (I-R)$ is non-vanishing and holomorphic in $%
\rho $ for $|\rho |<r_{1}r_{2}$. $\qed$
\end{theorem}

We can define a standard basis of cycles $\{a_{1},b_{1},\ldots $ $%
a_{g+1},b_{g+1}\}$ on the sewn genus $g+1$ surface as follows, where the set 
$\{a_{1},b_{1},\ldots $ $a_{g},b_{g}\}$ is the original basis. Then $a_{g+1}$
is defined as the contour $\mathcal{C}_{2}$ on $\hat{\mathcal{S}}$ whereas $%
b_{g+1}$ is defined to be a path chosen in $\hat{\mathcal{S}}$ from $%
z_{1}=z_{0}$ to $z_{2}=\rho /z_{0}$ which points are identified on the sewn
surface. Integrating (\ref{omgplus1}) along a $b_{r}$ cycle on $\mathcal{S}$
for $r=1,\ldots g$ gives $g$ holomorphic 1-forms for $x\in \hat{\mathcal{S}}$%
\begin{equation}
\nu_{r}^{(g+1)}(x)=\nu_{r}^{(g)}(x)-a(x)(I-R)^{-1}\bar{\alpha}_{r}^{T},
\label{nu1rho}
\end{equation}%
where $\bar{\alpha}_{r}=(\alpha_{r,\bar{a}}(k))$ with 
\begin{equation}
\alpha_{r,a}(k)=\oint_{b_{r}}a_{a}(x,k).  \label{alphari}
\end{equation}%
We then find from (\ref{nui}), (\ref{period}), (\ref{nu1rho}) and (\ref%
{alphari}) that for $r,s=1,\ldots g$%
\begin{equation*}
2\pi i\Omega_{rs}^{(g+1)}=2\pi i\Omega_{rs}^{(g)}-\alpha_{r}(I-R)^{-1}\bar{%
\alpha}_{s}^{T}.\quad
\end{equation*}

The remaining normalized holomorphic one form $\nu_{g+1}^{(g+1)}$ can be
expressed in terms of the normalized differential of the third kind $\omega
_{p_{2}-p_{1}}^{(g)}$ of (\ref{omp2p1}) with weighted moments 
\begin{equation}
\beta_{a}(k)=\frac{\rho ^{k/2}}{\sqrt{k}}\frac{1}{2\pi i}\int_{\mathcal{C}%
_{a}(z_{a})}(\omega_{p_{2}-p_{1}}^{(g)}+(-1)^{1+a}\frac{dz_{a}}{z_{a}}%
)z_{a}^{-k}.  \label{betai}
\end{equation}%
Then by Cauchy's theorem we find that \cite{Y}

\begin{lemma}[op. cite, Corollary 5]
\label{Lemma_nugplus1} The normalized holomorphic one form $\nu
_{g+1}^{(g+1)}$ is given by 
\begin{equation}
\nu _{g+1}^{(g+1)}(x)=\omega _{p_{2}-p_{1}}^{(g)}(x)+\frac{1}{2\pi i}%
\sum_{a=1,2}\oint_{\mathcal{C}_{a}(z)}\omega ^{(g+1)}(x,z)\int^{z}(\omega
_{p_{2}-p_{1}}^{(g)}+(-1)^{1+a}\frac{dz_{a}}{z_{a}}).\qed
\label{nugplusoneinteqn}
\end{equation}
\end{lemma}

Hence integrating (\ref{nugplusoneinteqn}) over a $b_{r}$ cycle and using (%
\ref{omgplus1}) we find for $r=1,\ldots ,g$ that 
\begin{equation*}
2\pi i\Omega_{rg+1}^{(g+1)}=\int_{p_{1}}^{p_{2}}\nu_{r}^{(g)}-\alpha
_{r}(I-R)^{-1}\bar{\beta}^{T}.
\end{equation*}

Finally $\Omega_{g+1g+1}^{(g+1)}$ is described in \cite{Y}:

\begin{lemma}[op. cite. Lemma 5]
\label{Lemma_omg1g1} $\Omega_{g+1g+1}^{(g+1)}\,$is given by 
\begin{eqnarray*}
2\pi i\Omega_{g+1g+1}^{(g+1)} &=&\log (\frac{\rho }{z_{0}^{2}}%
)+\int_{z_{1}^{-1}(z_{0})}^{z_{2}^{-1}(z_{0})}\omega_{p_{2}-p_{1}}^{(g)} \\
&&+\sum_{a=1,2}\frac{1}{2\pi i}\oint_{\mathcal{C}_{a}}\nu
_{g+1}^{(g+1)}(z)\int_{z_{a}^{-1}(z_{0})}^{z}(\omega
_{p_{2}-p_{1}}^{(g)}+(-1)^{1+a}\frac{dz_{a}}{z_{a}}),
\end{eqnarray*}%
where the logarithmic branch is determined by the choice of the cycle $%
b_{g+1}$ as a path in $\hat{\mathcal{S}}$ \ from $z_{1}=z_{0}$ to $%
z_{2}=\rho /z_{0}$. $\qed$
\end{lemma}

Substituting $\nu_{g+1}^{(g+1)}\,$from (\ref{nugplusoneinteqn}) one
eventually obtains \cite{Y} 
\begin{equation*}
2\pi i\Omega_{g+1g+1}^{(g+1)}=\log \rho +C_{0}-\beta (I-R)^{-1}\bar{\beta}%
^{T},
\end{equation*}%
where 
\begin{equation*}
C_{0}=\lim_{u\rightarrow 0}[\int_{z_{1}^{-1}(u)}^{z_{2}^{-1}(u)}\omega
_{p_{2}-p_{1}}^{(g)}-2\log u].
\end{equation*}%
However from (\ref{omp2p1prime}) we may express $C_{0}$ in terms of the
prime form%
\begin{eqnarray*}
C_{0} &=&\lim_{u\rightarrow 0}\log \frac{%
K^{(g)}(z_{2}^{-1}(u),p_{2})K^{(g)}(z_{1}^{-1}(u),p_{1})}{%
u^{2}K^{(g)}(z_{2}^{-1}(u),p_{1})K^{(g)}(z_{1}^{-1}(u),p_{2})} \\
&=&-\log (-z_{1}^{\prime }(p_{1})z_{2}^{\prime
}(p_{2})K^{(g)}(p_{2},p_{1})^{2}),
\end{eqnarray*}%
where $\frac{d}{du}z_{a}^{-1}(u)|_{u=0}=1/z_{a}^{\prime }(p_{a})$ and using $%
K^{(g)}(p_{2},p_{1})=-K^{(g)}(p_{1},p_{2})$. We therefore find altogether
that

\begin{theorem}
\label{Theorem_periodg_rho} The genus $g+1$ period matrix for $|\rho
|<r_{1}r_{2}$ is given by 
\begin{eqnarray}
2\pi i\Omega_{rs}^{(g+1)} &=&2\pi i\Omega_{rs}^{(g)}-\alpha_{r}(I-R)^{-1}%
\bar{\alpha}_{s}^{T},\quad r,s=1,\ldots ,g,  \label{Omg1rs} \\
2\pi i\Omega_{rg+1}^{(g+1)} &=&\int_{p_{1}}^{p_{2}}\nu_{r}^{(g)}-\beta
(I-R)^{-1}\bar{\alpha}_{r}^{T},\quad r=1,\ldots ,g,  \label{Omg1rg1} \\
2\pi i\Omega_{g+1g+1}^{(g+1)} &=&\log (\frac{-\rho }{z_{1}^{\prime
}(p_{1})z_{2}^{\prime }(p_{2})K^{(g)}(p_{2},p_{1})^{2}})-\beta (I-R)^{-1}%
\bar{\beta}^{T},  \label{Omg1g1g1}
\end{eqnarray}%
where $\Omega ^{(g+1)}$ is holomorphic in $\rho $ for $0<|\rho |<r_{1}r_{2}$
and the logarithmic branch is determined by the choice of the cycle $%
b_{g+1}. $ $\qed$
\end{theorem}

We finally obtain the following holomorphic properties for $\omega ^{(g+1)}$
, $\Omega ^{(g+1)}$ and $\det (I-R)$. The proof follows a similar argument
to that for Propositions \ref{Prop_hol_omega} and \ref{Prop_Det_holomorphic}.

\begin{proposition}
\label{Prop_hol_omega_Det_rho} Suppose that $\omega ^{(g)}$ is a holomorphic
function of a complex parameter $\mu $ for $\left\vert \mu \right\vert <S$.
Then $\det (I-R)$ is non-vanishing and both $\omega ^{(g+1)}$ and $\det
(I-R) $ are holomorphic in $\mu $ for $\left\vert \mu \right\vert <S$ with $%
|\rho |<r_{1}r_{2}$ whereas $\Omega ^{(g+1)}$ is holomorphic in $\mu $ for $%
\left\vert \mu \right\vert <S$ with $0<|\rho |<r_{1}r_{2}$. $\qed$
\end{proposition}

\subsection{Self-Sewing a Sphere to form a Torus}

It is instructive to consider two separate examples of sewing a Riemann
sphere to itself to form a torus. The first is mainly illustrative whereas
the second is related to some later genus two considerations wherein the
Catalan numbers arise in an interesting and surprising way. In both cases $%
\omega ^{(1)}$ is given by (\ref{omega1}) with an appropriately identified
modular parameter $\tau $.

\subsubsection{Simplest Case}

Let $\mathcal{S}_{0}=\mathbb{C}\cup \{\infty \}$ be the Riemann sphere with
bilinear form (\ref{omega0}). Choose local coordinates $z_{1}=z$ in the
neighborhood of the origin and $z_{2}=1/z^{\prime }$ for $z^{\prime }$ in
the neighborhood of the point at infinity. Identify the annular regions $%
|q|r_{\bar{a}}^{-1}\leq \left\vert z_{a}\right\vert \leq r_{a}$ for a
complex parameter $q$ obeying $|q|\leq r_{1}r_{2}$ via the sewing relation 
\begin{equation}
z=qz^{\prime }.  \label{spheresew1}
\end{equation}%
Note that the annular regions do not intersect on the sphere provided $%
r_{1}r_{2}<1$ so that $|q|<1$. We then find \cite{Y}

\begin{proposition}
\label{Prop_torus1} $q=\exp (2\pi i\tau )$ where $\tau $ is the torus
modular parameter.
\end{proposition}

\noindent \textbf{Proof. }The 1-forms (\ref{akdef2}) are 
\begin{eqnarray*}
a_{1}(k,x) &=&\sqrt{k}q^{k/2}x^{-k-1}dx, \\
a_{2}(k,x) &=&-\sqrt{k}q^{k/2}x^{k-1}dx,
\end{eqnarray*}%
so that $A(k,l)=0$ and $B(k,l)=q^{k}\delta _{k,l}\,$ in (\ref{Rformula})
giving 
\begin{equation*}
I-R=\mathrm{diag}(1-q,1-q,\ldots ,1-q^{k},1-q^{k},\ldots )
\end{equation*}%
Hence Lemma \ref{Lemma_omg1g} gives for $x,y\in \hat{\mathcal{S}_{0}}$%
\begin{equation*}
\omega ^{(1)}(x,y)=\{\frac{xy}{(x-y)^{2}}+\sum_{k\geq 1}\frac{kq^{k}}{1-q^{k}%
}[(\frac{x}{y})^{k}+(\frac{y}{x})^{k}]\}\frac{dxdy}{xy}.
\end{equation*}%
Under the conformal map $z\rightarrow \log z$ we then verify $\omega
^{(1)}(u,v)=P_{2}(\tau ,u-v)dudv$ with $u=\log x$ and $v=\log y$ using (\ref%
{P2exp}) where $q=\exp (2\pi i\tau )$. The sewing relation (\ref{spheresew1}%
) is then just the standard torus periodicity relation $\log z=\log
z^{\prime }+2\pi i\tau $.

Alternatively, we may apply (\ref{Omg1g1g1}) of Theorem \ref%
{Theorem_periodg_rho} using $z_{2}=1/z-1/p_{2}$ and then consider $%
p_{2}\rightarrow \infty $. Then $K^{(0)}(p_{2},0)=p_{2}$ with $\omega
_{p_{2}-0}^{(1)}(x)=(\frac{1}{x-p_{2}}-\frac{1}{x})dx$ so that $\beta
_{a}(k)=0$ and $z_{1}^{\prime }(0)z_{2}^{\prime
}(p_{2})K^{(0)}(p_{2},0)^{2}=-1$ independent of $p_{2}$. This implies that $%
2\pi i\tau =2\pi i\Omega_{11}^{(1)}=\log q$ again. $\qed$

\begin{remark}
\label{Rem_simple_Dehn_twist} \ The modular transformation $\tau \rightarrow
\tau +1$ is generated by a continuous variation in the sewing parameter $%
\exp (i\theta )q$ for $0\leq \theta \leq 2\pi $. This corresponds to a Dehn
twist $b_{1}\rightarrow a_{1}+b_{1}$ in the $b_{1}$ cycle chosen in Lemma %
\ref{Lemma_omg1g1} and Theorem \ref{Theorem_periodg_rho} so that $2\pi
i\Omega _{11}^{(1)}=\log q$ is evaluated on the next logarithmic branch.
\end{remark}

\begin{remark}
\label{Rem_eta_sq} \ $\omega ^{(1)}(x,y)$ and $\det(1-R)=\prod_{k\geq
1}(1-q^{k})^{2}$ are clearly holomorphic for $\left\vert q\right\vert <1$ as
expected from Theorems \ref{Theoremomgp1holo} and \ref{Theorem_Det_rho}.
\end{remark}

\subsubsection{General Self-Sewing of a Sphere and the Catalan Series}

\label{Subsubsect_Catalan}For $z\in \mathcal{S}_{0}$ choose local
coordinates $z_{1}=z$ in the neighborhood of the origin and $z_{2}=z^{\prime
}-w$ for $z^{\prime }$ in the neighborhood of $w\in \mathcal{S}_{0}$.
Identify the annuli $|{\rho }|r_{2}^{-1}\leq \left\vert z\right\vert \leq
r_{1}$ and $|{\rho }|r_{1}^{-1}\leq \left\vert z^{\prime }-w\right\vert \leq
r_{2}$ for $|{\rho }|\leq r_{1}r_{2}$ via the sewing relation 
\begin{equation}
z(z^{\prime }-w)=\rho .  \label{spheresew2}
\end{equation}%
The two annular regions do not intersect provided $|w|>r_{1}+r_{2}\geq
r_{1}+|{\rho }|r_{1}^{-1}\,\geq 2|\rho |^{1/2}$. The lower bound occurs for $%
r_{1}=r_{2}=|\rho |^{1/2}$ and is realized when the two annuli become
degenerate (infinitesimally thin) and touch at the point $z_{1}=-z_{2}=w/2$
with $w^{2}=-4\rho $. Thus defining 
\begin{equation*}
\chi =-\frac{\rho }{w^{2}},
\end{equation*}%
then $\chi =\frac{1}{4}$ is the degenerate point. We define the Catalan
series\footnote{%
The Catalan series is more usually defined to be $1+f(\chi )$.} to be the
series $f(\chi )$ convergent for $|\chi |<\frac{1}{4}$ satisfying 
\begin{equation}
\chi =\frac{f}{(1+f)^{2}}.  \label{chi(f)}
\end{equation}%
Thus 
\begin{eqnarray}
f(\chi ) &=&\frac{1-\sqrt{1-4\chi }}{2\chi }-1=\sum_{n\geq 1}\frac{1}{n}%
\binom{2n}{n+1}\chi ^{n}  \notag \\
&=&\chi +2\chi ^{2}+5\chi ^{3}+14\chi ^{4}+O\left( \chi ^{5}\right) ,
\label{Catalan}
\end{eqnarray}%
The coefficients $\frac{1}{n}\binom{2n}{n+1}$ are the Catalan numbers which
occur in a remarkably wide range of combinatorial settings e.g. \cite{St}.

\begin{proposition}
\label{Prop_torus2} For the sewing described by (\ref{spheresew2}), the
torus modular parameter is $q=f(\chi )$, the Catalan series.
\end{proposition}

\noindent \textbf{Proof.} Define the M\"{o}bius transformation 
\begin{equation}
z\mapsto \gamma .z=\frac{w}{1+f}(\frac{z-f}{z-1}),  \label{mobius}
\end{equation}%
where $f=f(\chi )$. Then with $z=\gamma .Z$ and $z^{\prime }=\gamma
.Z^{\prime }$ the sewing relation (\ref{spheresew2}) becomes on using (\ref%
{chi(f)}) 
\begin{equation*}
Z=fZ^{\prime }.
\end{equation*}%
Thus we recover the earlier sewing relation of (\ref{spheresew1}) with
modular parameter $q=f(\chi )$.

This result can be verified from (\ref{Omg1g1g1}) of Theorem \ref%
{Theorem_periodg_rho} as follows. With $\omega ^{(0)}$ of (\ref{omega0}) the
basis of 1-forms (\ref{akdef2}) is given by 
\begin{eqnarray}
a_{1}^{(0)}(k,x) &=&\sqrt{k}\rho ^{k/2}x^{-k-1}dx,  \notag \\
a_{2}^{(0)}(k,x) &=&\sqrt{k}\rho ^{k/2}(x-w)^{-k-1}dx,  \label{akgenuszero}
\end{eqnarray}%
with 
\begin{eqnarray}
R^{(0)} &=&-\left[ 
\begin{array}{ll}
B^{(0)} & A^{(0)} \\ 
A^{(0)} & B^{(0)T}%
\end{array}%
\right] ,  \notag \\
A^{(0)}(k,l) &=&0,\quad B^{(0)}(k,l)=\frac{(-\chi )^{(k+l)/2}}{\sqrt{kl}}%
\frac{(-1)^{k+1}(k+l-1)!}{(k-1)!(l-1)!},  \notag \\
\omega _{w-0}^{(0)}(x) &=&(\frac{1}{x-w}-\frac{1}{x})dx,  \notag \\
K^{(0)}(w,0) &=&w,  \notag \\
\beta ^{(0)}(k) &=&\frac{(-\chi )^{k/2}}{\sqrt{k}}[-1,(-1)^{k}],
\label{Genuszerodata}
\end{eqnarray}%
where the $0$ superscript indicates the genus of the sphere. After some
calculation, we find that $\tau $ is given by%
\begin{equation*}
2\pi i\tau =2\pi i\Omega _{11}^{(1)}=\log \chi +2\sum_{k\geq 1}\frac{1}{k}%
\chi ^{k}\sum\limits_{n\geq 1}S_{n,k}(\chi ),
\end{equation*}%
where $S_{1,k}(\chi )=1$ and 
\begin{eqnarray}
S_{n,k}(\chi ) &=&\sum_{k_{n-1},\ldots k_{1}\geq 1}\chi ^{k_{n-1}+\ldots
+k_{1}}\binom{k+k_{n-1}-1}{k_{n-1}}\binom{k_{n-1}+k_{n-2}-1}{k_{n-2}}  \notag
\\
&&\ldots \binom{k_{2}+k_{1}-1}{k_{1}},  \label{Snk}
\end{eqnarray}%
for $n>1$. We will show below that 
\begin{equation}
\sum\limits_{n\geq 1}S_{n,k}(\chi )=(1+f(\chi ))^{k},  \label{Catalan^k}
\end{equation}%
which implies $\sum_{k\geq 1}\frac{1}{k}\chi ^{k}\sum\limits_{n\geq
1}S_{n,k}(\chi )=-\log (1-\chi (1+f))=$ $\log (1+f)$ from (\ref{Catalan}).
Therefore $2\pi i\tau =\log \chi +2\log (1+f)=\log f$ so that $q=f$ as
claimed.

It remains to prove (\ref{Catalan^k}). Since $\sum_{k_{1}\geq 1}\chi ^{k_{1}}%
\binom{k_{2}+k_{1}-1}{k_{1}}$ $=(1-\chi )^{-k_{2}}-1$ we find for $n>1$ that 
\begin{eqnarray*}
S_{n,k}(\chi ) &=&\sum_{k_{n-1}\geq 1}\chi ^{k_{n-1}}\binom{k+k_{n-1}-1}{%
k_{n-1}}\ldots \\
&&\sum_{k_{2}\geq 1}(\frac{\chi }{1-\chi })^{k_{2}}\binom{k_{3}+k_{2}-1}{%
k_{2}}-S_{n-1,k}(\chi ).
\end{eqnarray*}%
Repeating this process leads to 
\begin{equation*}
\sum_{n=1}^{N}S_{n,k}(\chi )=\left( \left[ \frac{1}{1-\frac{\chi }{1-\chi
/\cdots }}\right]_{N}\right) ^{k},
\end{equation*}%
where $[\frac{1}{1-\frac{\chi }{1-\chi /\cdots }}]_{N}$ denotes the $N^{th}$
term in the continued fraction expansion of $F=1/(1-\chi F)$ whose solution
from (\ref{Catalan}) is $F=1+f$. $\qed$

\begin{remark}
\label{Rem_general_Dehn_twist} \ The modular transformation $\tau
\rightarrow \tau +1$ is generated by a continuous variation in the sewing
parameter $\exp (i\theta )\rho $ for $0\leq \theta \leq 2\pi $.
\end{remark}

Using Lemma \ref{Lemma_omg1g} and comparing to $\omega ^{(1)}$ of (\ref%
{omega1}) results in novel expressions for Eisenstein series $E_{n}(q)$ for $%
q=f(\chi )$. Thus, for example, one finds

\begin{proposition}
\label{Prop_CatalanEisenstein} 
\begin{equation}
E_{2}(q=f(\chi ))=-\frac{1}{12}+\frac{2\chi }{1-4\chi }(1+B^{(0)})^{-1}(1,1),
\label{E2chi}
\end{equation}%
where $(1,1)$ refers to the $(k,l)=(1,1)$ element of the infinite matrix $%
(1+B^{(0)})^{-1}$.
\end{proposition}

\noindent \textbf{Proof. }From (\ref{omgplus1}) and (\ref{akgenuszero}) we
have \textbf{\ }%
\begin{equation*}
\omega ^{(1)}(x,y)=\frac{dxdy}{(x-y)^{2}}-a^{(0)}(x)(I-R^{(0)})^{-1}(\bar{a}%
^{(0)})^{T}(y).
\end{equation*}%
But $\omega ^{(1)}(x,y)=\omega ^{(1)}(u,v)=P_{2}(u-v,\tau )dudv$ with $\tau =%
\frac{1}{2\pi i}\log f(\chi )$ from Proposition \ref{Prop_torus2} with $%
x=\gamma .e^{u}$ and $y=\gamma .e^{v}$ using (\ref{mobius}). Then, on
substituting for $u,v$ into $\omega ^{(1)}(x,y)$ one eventually finds using (%
\ref{Genuszerodata}) that 
\begin{eqnarray*}
\omega ^{(1)}(x,y) &=&\frac{e^{u-v}dudv}{(e^{u-v}-1)^{2}}-\sum_{k,l\geq
1}(1+B^{(0)})^{-1}(k,l)\sqrt{kl}(-\frac{\chi }{1-4\chi })^{(k+l)/2} \\
&&.(e^{u}-1)^{k-1}(e^{v}-1)^{l-1}[(\frac{1-f}{e^{u}-f})^{k+1}(\frac{1-f}{%
1-fe^{v}})^{l+1} \\
&&+(-1)^{k+l}(\frac{1-f}{1-fe^{u}})^{k+1}(\frac{1-f}{e^{v}-f}%
)^{l+1}]e^{u+v}dudv.
\end{eqnarray*}%
using $1-f=(1+f)\sqrt{1-4\chi }$. Expanding in $u,v$ we then find that 
\begin{equation}
\omega ^{(1)}(x,y)=[\frac{1}{(u-v)^{2}}-\frac{1}{12}+\frac{2\chi }{1-4\chi }%
(1+B^{(0)})^{-1}(1,1)+O(u,v)]dudv,  \label{omuv}
\end{equation}%
from which the result follows on comparison with (\ref{P2}). $\qed$

\section{Self-Sewing a Torus to form a Genus Two Riemann Surface}

\label{Sect_Rho_Torus}

\subsection{The Genus Two Period Matrix in the $\protect\rho $ Formalism}

We now apply the $\rho $-formalism to sew a twice punctured torus with
modulus $\tau $ and punctures separated by $w$ to form a genus two Riemann
surface with period matrix $\Omega ^{(2)}(\tau ,w,\rho )$. We will see that $%
\Omega ^{(2)}$ is holomorphic for $(\tau ,w,\rho )$ in an appropriate domain 
$\mathcal{D}^{\rho }$. We again provide a description of $\Omega ^{(2)}$ in
terms of a sum of weights of \emph{necklaces}. There is a holomorphic
mapping $F^{\rho }:\mathcal{D}^{\rho }\rightarrow \mathbb{H}_{2}$, and we
describe its equivariance properties with respect to a certain group. The
logarithmic contribution $\log (-\rho /K^{2})$ to $\Omega _{22}^{(2)}$ in (%
\ref{Omg1g1g1}) gives rise to a subtle analytic structure which we discuss
in some detail. Finally, we prove that $F^{\rho }$ is invertible in a
certain domain.

Consider a framed torus (cf. Subsection \ref{Subsect_closeddisk}) $\mathcal{S%
}=\mathbb{C}/\Lambda $ where $\Lambda \subseteq \mathbb{C}$ is a lattice
with positively oriented basis $(\sigma ,\varsigma )$ and modulus $\tau
=\sigma /\varsigma \in \mathbb{H}_{1}$. Define annuli $\mathcal{A}%
_{a},a=1,2, $ centered at $z=0$ and $z=w$ of \ $\mathcal{S}$ with local
coordinates $z_{1}=z$ and $z_{2}=z-w$ respectively. Take the outer radius of 
$\mathcal{A}_{a}$ to be $r_{a}<\frac{1}{2}D(\Lambda )$ and the inner radius
to be $|{\rho }|/r_{\bar{a}}$, with $|\rho |\leq r_{1}r_{2}<\frac{1}{4}%
D(\Lambda )^{2} $ (cf. Lemma \ref{Lemmadisk}). Identifying the annuli
according to the sewing relation (\ref{rhosew}) $z_{1}z_{2}=\rho $ gives
rise to a compact Riemann surface of genus $2$.

As in the remarks following Lemma \ref{Lemmadisk}, we now take $\Lambda
=\Lambda _{\tau }$ with basis $(2\pi i\tau ,2\pi i)$ and with $w$ in the
fundamental parallelogram for $\Lambda _{\tau }$ with sides $(2\pi i\tau
,2\pi i)$. As with the sphere example above, the two annuli must not
intersect. This requires the inequalities $|w-\lambda |>r_{1}+r_{2}\geq
2|\rho |^{1/2}$ to hold for $\lambda \in \Lambda _{\tau }$. Thus we find 
\begin{equation*}
2|\rho |^{1/2}<|w|<D(\Lambda _{\tau })-2|\rho |^{1/2}.
\end{equation*}%
Notice that this implies $|\rho |<\frac{1}{16}D(\Lambda _{\tau })^{2}$,
which refines the inequality satisfied by $\rho $ discussed above. As a
result of this discussion, we see that the relevant domain in the $\rho $%
-formalism is the following\footnote{%
The footnote relating to (\ref{Deps}) concerning notation applies here too}: 
\begin{equation}
\mathcal{D}^{\rho }=\{(\tau ,w,\rho )\in \mathbb{H}_{1}\times \mathbb{C}%
\times \mathbb{C}\ |\ |w-\lambda |>2|\rho |^{1/2}>0,\ \lambda \in \Lambda
_{\tau }\}.  \label{Drho}
\end{equation}%
We may apply Theorem \ref{Theorem_periodg_rho} to determine $\Omega
^{(2)}(\tau ,w,\rho )$. We find:

\begin{theorem}
\label{Theorem_rhoperiod}Sewing determines a holomorphic map%
\begin{eqnarray}
F^{\rho }:\mathcal{D}^{\rho } &\rightarrow &\mathbb{H}_{2},  \notag \\
(\tau ,w,\rho ) &\mapsto &\Omega ^{(2)}(\tau ,w,\rho ).  \label{Frhomap}
\end{eqnarray}
\end{theorem}

\begin{proposition}
$\Omega ^{(2)}=\Omega ^{(2)}(\tau ,w,\rho )$ is given by%
\begin{eqnarray}
2\pi i\Omega _{11}^{(2)} &=&2\pi i\tau -\rho \sigma ((I-R)^{-1}(1,1)),
\label{Om11rho} \\
2\pi i\Omega _{12}^{(2)} &=&w-\rho ^{1/2}\sigma ((\beta (I-R)^{-1}(1)),
\label{Om12rho} \\
2\pi i\Omega _{22}^{(2)} &=&\log (-\frac{\rho }{K(\tau ,w)^{2}})-\beta
(I-R)^{-1}\bar{\beta}^{T},  \label{Om22rho}
\end{eqnarray}%
where the branch of the $\log $ function in (\ref{Om22rho}) is determined by
the choice of the cycle $b_{2}$. Here, $R=R(\tau ,w,{\rho })=(R_{ab}(k,l))$
is an infinite matrix with indices $k,l=1,2,3,\ldots $ and $a,b=1,2$; $\beta
=\beta (\tau ,w,{\rho })=(\beta _{a}(k))$ is an infinite row vector; $(1,1)$
and $(1)$ are the $(1,1)$- and $(1)$- (block) entries of a matrix; $\sigma
(M)$ denotes sum over the entries of a finite matrix; and 
\begin{eqnarray}
R(k,l) &=&-\frac{\rho ^{(k+l)/2}}{\sqrt{kl}}\left[ 
\begin{array}{cc}
D(k,l,\tau ,w) & C(k,l,\tau ) \\ 
C(k,l,\tau ) & D(l,k,\tau ,w)%
\end{array}%
\right] ,  \label{Rg2} \\
\beta (k) &=&\frac{\rho ^{k/2}}{\sqrt{k}}(P_{k}(\tau ,w)-E_{k}(\tau
))[-1,(-1)^{k}],  \label{betag2}
\end{eqnarray}%
with notation as in Section \ref{Sect_ellipticfunctions}.
\end{proposition}

\noindent \textbf{Proof.} \ Since $\omega ^{(1)}(x,y)=P_{2}(x-y)dxdy$ from (%
\ref{P1Pnexpansion}) we find that the set of 1-forms (\ref{akdef2}) with
periods $(2\pi i\tau ,2\pi i)$ is given by

\begin{eqnarray*}
a_{1}(k,x) &=&a_{1}(k,x,\tau ,\rho )=\sqrt{k}\rho ^{k/2}P_{k+1}(\tau ,x)dx,
\\
a_{2}(k,x) &=&a_{2}(k,x,\tau ,\rho )=a_{1}(k,x-w).
\end{eqnarray*}%
The matrices $A(k,l)$, $B(k,l)\,$ in (\ref{ABdef}) are given directly from
the expansions (\ref{Pkexpansion}) and (\ref{Pkzwexpansion}) which are
convergent on $\mathcal{D}^{\rho }$ resulting in (\ref{Rg2}). $\alpha
_{1,a}(k)$ of (\ref{alphari}) is independent of $a=1,2$ with 
\begin{equation*}
\alpha_{1,a}(k)=\oint_{b_{1}}a_{a}(k,\cdot )=\rho ^{1/2}\delta_{k,1}.
\end{equation*}%
Hence $2\pi i\Omega_{11}^{(2)}$ is as stated from (\ref{Omg1rs}) of Theorem %
\ref{Theorem_periodg_rho} for $(\tau ,w,\rho )\in \mathcal{D}^{\rho }$. From
Example \ref{Exampleom1xy} we know that $\omega_{w-0}^{(1)}(x)=(P_{1}(\tau
,x-w)-P_{1}(\tau ,x))dx$ and the prime form is $K^{(1)}(x,y)=K(\tau ,x-y)$.
We obtain the given moments (\ref{betag2}) of $\omega_{w-0}^{(1)}(x)$ from (%
\ref{P1Pnexpansion}). Hence since $\nu ^{(1)}(x)=dx,$ we find $2\pi i\Omega
_{12}^{(2)}$ is as given from (\ref{Omg1rg1}) of Theorem \ref%
{Theorem_periodg_rho} for $(\tau ,w,\rho )\in \mathcal{D}^{\rho }$. Finally
applying (\ref{Omg1g1g1}) with $K^{(1)}(w,0)=K(\tau ,w)$ we obtain (\ref%
{Om22rho}) for $(\tau ,w,\rho )\in \mathcal{D}^{\rho }$.

$\Omega _{ij}^{(2)}(\tau ,w,\rho )$ is holomorphic in $\rho $ for $0<|{\rho }%
|<r_{1}r_{2}$ from Theorem \ref{Theorem_periodg_rho}. Proposition \ref%
{Prop_hol_omega_Det_rho} then states that $\Omega _{ij}^{(2)}(\tau ,w,\rho )$
is also holomorphic in $\tau \in \mathbb{H}_{1}$. We also need to show that $%
\Omega _{ij}^{(2)}(\tau ,w,\rho )$ is holomorphic in $w$. Since $%
Y=(I-R)^{-1} $ converges for $|{\rho }|<r_{1}r_{2}$ then, following an
argument similar to that in Proposition \ref{Prop_hol_omega}, we find that $%
\Omega _{ij}^{(2)}(\tau ,w,\rho )$ is continuous in $w$ for $(\tau ,w,\rho
)\in \mathcal{D}^{\rho }$. The Weierstrass functions $P_{k}(\tau ,w)$ for $%
k\geq 1 $ are holomorphic in $w$. Hence the $\rho $ expansion coefficients $%
\Omega _{ij}^{(2)}(\tau ,w,\rho )$ are holomorphic functions in $w$ since
they consist of finite sums and products of these Weierstrass functions.
Hence, by Lemma \ref{Lem_fxyhol}, $\Omega _{ij}^{(2)}(\tau ,w,\rho )$ is
holomorphic in $w$. Finally, by Hartog's Theorem, $\Omega _{ij}^{(2)}$ is
holomorphic on $\mathcal{D}^{\rho }$. $\qed$

\subsection{Necklace Expansion for $\Omega ^{(2)}$}

We introduce a graphical interpretation for the $\rho$ period matrix
formulas analogous to that described earlier for the $\epsilon$-expansion.
Consider the set of \emph{necklaces} $\mathcal{N}=\{N\}$: they are connected
graphs with $n\geq 2$ nodes, $n-2$ of which have valency $2$ and two of
which have valency $1$, together with an orientation, say from left to
right. Furthermore, each vertex carries two labels $k,a$ with $k$ a positive
integer and $a=1$ or $2$. A typical necklace in the $\rho$-formalism looks
as follows: 
\begin{equation*}
\ \overset{k_{1},a_{1}}{\bullet }\longrightarrow\overset{k_{2},a_{2}}{%
\bullet }\longrightarrow\overset{k_{3},a_{3}}{\bullet }\longrightarrow 
\overset{k_{4},a_{4}}{\bullet }
\end{equation*}%
We define the degenerate necklace $N_{0}$ to be a single node with no edges.
Next we define a weight function 
\begin{equation*}
\omega :\mathcal{N}\longrightarrow \mathbb{C}[P_{2}(\tau ,w),P_{3}(\tau
,w),E_{2}(\tau ),E_{4}(\tau ),E_{6}(\tau ),\rho ^{1/2}].
\end{equation*}%
If $N\in \mathcal{N}$ has edges $E$ labelled as \ $\overset{k,a}{\bullet }%
\overset{}{\longrightarrow }\overset{l,b}{\bullet }$ \ then we define 
\begin{eqnarray*}
\omega (E) &=&R_{ab}(k,l,\tau ,w,{\rho }), \\
\omega (N) &=&\prod \omega (E),
\end{eqnarray*}%
with $R_{ab}(k,l)$ as in (\ref{Rg2}) and where the product is taken over all
edges of $N$. We further define $\omega (N_{0})=1$.

The necklaces with prescribed end nodes labelled $(k,a;l,b)$ look as
follows: 
\begin{equation*}
\overset{k,a}{\bullet }\longrightarrow\overset{k_{1},a_{1}}{\bullet }\ldots%
\overset{k_{2},a_{2}}{\bullet }\longrightarrow\overset{l,b}{\bullet }\hspace{%
10mm}\mbox{(type $(k,a;l,b)$)}
\end{equation*}%
We set 
\begin{equation*}
\mathcal{N}_{k,a;l,b}=\{\mbox{isomorphism classes of
necklaces of type}\ (k,a;l,b) \}.
\end{equation*}%
As in Lemma \ref{Lemma_cheq} we obtain

\begin{lemma}
\label{Lemma_necklace_rho}\noindent\ We have for $k,l\geq 1$ 
\begin{equation*}
(I-R)_{ab}^{-1}(k,l)=\sum_{N\in \mathcal{N}_{k,a;l,b}}\omega (N).\ \ \ \ \ \
\ \qed
\end{equation*}
\end{lemma}

Finally it is convenient to define 
\begin{eqnarray}
\omega_{11} &=&\sum_{a,b=1,2}\sum_{N\in \mathcal{N}_{1,a;1,b}}\omega (N), 
\notag \\
\omega_{\beta 1} &=&\sum_{a,b=1,2}\sum_{k\geq 1}\beta_{a}(k)\sum_{N\in 
\mathcal{N}_{k,a;1,b}}\omega (N),  \notag \\
\omega_{1\bar{\beta}} &=&\sum_{a,b=1,2}\sum_{k\geq 1}\bar{\beta}%
_{b}(k)\sum_{N\in \mathcal{N}_{1,a;k,b}}\omega (N),  \notag \\
\omega_{\beta \bar{\beta}} &=&\sum_{a,b=1,2}\sum_{k,l\geq 1}\beta_{a}(k)\bar{%
\beta}_{b}(l)\sum_{N\in \mathcal{N}_{k,a;l,b}}\omega (N).  \label{om_weights}
\end{eqnarray}%
Note that $R_{ab}(k,l)=R_{\bar{b}\bar{a}}(l,k)$, so that $\omega_{\beta
1}=\omega_{1\bar{\beta}}$. Then from Theorem \ref{Theorem_rhoperiod} we have

\begin{proposition}
\label{Prop_rhoperiod_graph} 
\begin{eqnarray*}
2\pi i\Omega _{11}^{(2)} &=&2\pi i\tau -\rho \omega _{11}, \\
2\pi i\Omega _{12}^{(2)} &=&w-\rho ^{1/2}\omega _{\beta 1}, \\
2\pi i\Omega _{22}^{(2)} &=&\log (-\frac{\rho }{K(\tau ,w)^{2}})-\omega
_{\beta \bar{\beta}}.\ \ \ \ \ \ \ \qed
\end{eqnarray*}
\end{proposition}

\subsection{Equivariance of $F^{\protect\rho }$}

\label{Subsec_Frho}In Subsection \ref{Subsec_Feps} we defined a subgroup $%
G\subset Sp(4,\mathbb{Z})$ which preserves the domain $\mathcal{D}^{\epsilon
}$, and proved the equivariance of $F^{\epsilon }$ under the action of $G$.
In this section we wish to establish analogous equivariance properties in
the $\rho $-formalism. With this in mind, one might expect that the map $%
F^{\rho }$ occurring in Theorem \ref{Theorem_rhoperiod} is the correct
analog of $F^{\epsilon }$. However, because of the logarithmic branch
structure of $\Omega _{22}^{(2)},$ it is necessary to lift $F^{\rho }$ to a
single-valued function $\hat{F}^{\rho }$ on a certain covering space $%
\mathcal{\hat{D}}^{\rho }$ for $\mathcal{D}^{\rho }$ before the correct
analogs can be established.

\subsubsection{Some Heisenberg and Jacobi-type groups}

In this subsection we consider some groups relevant to our enterprise, and
start with certain subgroups of $Sp(4,\mathbb{Z})$. For $(a,b,c)\in \mathbb{Z%
}^{3}$ set 
\begin{equation}
\mu (a,b,c)=\left( 
\begin{array}{cccc}
1 & 0 & 0 & b \\ 
a & 1 & b & c \\ 
0 & 0 & 1 & -a \\ 
0 & 0 & 0 & 1%
\end{array}%
\right) ,  \label{mudef}
\end{equation}%
with 
\begin{equation*}
A=\mu (1,0,0),\ B=\mu (0,1,0),\ C=\mu (0,0,1).
\end{equation*}%
The matrices (\ref{mudef}) form a subgroup $\hat{H}\subseteq Sp(4,\mathbb{Z})
$ which is a $2$-step nilpotent group with center isomorphic to $\mathbb{Z}$
and generated by $C$, and central quotient isomorphic to $\mathbb{Z}^{2}$.
Note that we have the presentation 
\begin{equation}
\hat{H}=\langle A,B,C\ |\ [A,B]C^{-2}=[A,B,C]=1\rangle .
\label{Hhatpresentation}
\end{equation}%
The `left' modular group $\Gamma _{1}$ (\ref{Gamma1Gamma2}) is also a
subgroup of $Sp(4,\mathbb{Z})$, and indeed it normalizes $\hat{H}$ according
to the conjugation formula 
\begin{equation}
\gamma ^{-1}\mu (u,v,w)\gamma =\mu ((u,v)\gamma ,w),\ \gamma \in \Gamma _{1}.
\label{gammaconjform}
\end{equation}%
Here it was convenient to abuse notation, taking 
\begin{equation}
\gamma =\left( 
\begin{array}{cccc}
a & 0 & b & 0 \\ 
0 & 1 & 0 & 0 \\ 
c & 0 & d & 0 \\ 
0 & 0 & 0 & 1%
\end{array}%
\right) \in \Gamma _{1}\ \mbox{and}\ (u,v)\gamma =(u,v)\left( 
\begin{array}{cc}
a & b \\ 
c & d%
\end{array}%
\right) .  \label{gamma1formrho}
\end{equation}%
In this way we get a subgroup $L=\hat{H}\Gamma _{1}\subseteq Sp(4,\mathbb{Z})
$ which is a split extension of $SL(2,\mathbb{Z})$ by $\hat{H}$. Note that $%
Z(L)=\langle C\rangle $, and that the central quotient $J=L/Z(L)\cong 
\mathbb{Z}^{2}\rtimes SL(2,\mathbb{Z})$ is the Jacobi group which figures in
the transformation laws of Jacobi forms (\cite{EZ}).

Let $H = \langle A, B \rangle$ be the subgroup of $\hat{H}$ generated by $A$
and $B$. It follows from (\ref{Hhatpresentation}) that $H$ has the
presentation 
\begin{eqnarray*}
H = \langle A, B \ | \ [A, B] = C^{\prime}, [A, C^{\prime}] = [B,
C^{\prime}] = 1 \rangle,
\end{eqnarray*}
where we have set $C^{\prime}= C^2$. We call $H$ the \emph{Heisenberg}
group, though $H$ and $\hat{H}$ (which are \emph{not} isomorphic) are often
confused in this regard. We see from (\ref{gammaconjform}) that $|\hat{H}:H|
= 2$ and that $\Gamma_1$ normalizes $H$. Thus $L_0 = H \Gamma_1$ is a
subgroup of $L$ of index $2$.

\begin{lemma}
\label{lemmaLaction} $L$ acts on $\mathcal{D}^{\rho}$ as follows: 
\begin{eqnarray}
\mu(a, b, c).(\tau, w, \rho) &=& (\tau, w+2 \pi i a \tau + 2 \pi i b, \rho),
\label{muaction} \\
\gamma. (\tau, w, \rho) &=& \left( \frac{a \tau+b}{c\tau + d}, \frac{w}{c
\tau + d}, \frac{\rho}{(c \tau + d)^2} \right).  \label{gammaonDrhoaction}
\end{eqnarray}
The kernel of the action is $Z(L)$, so that the effective action is that of $%
J = L/Z(L)$.
\end{lemma}

\noindent \textbf{Proof} Let us first work with the larger domain whereby we
allow the triple $(\tau, w, \rho)$ to lie in $\mathbb{H}_1 \times \mathbb{C}
\times \mathbb{C}$. Then it is easy to see that the first equality defines
an action of $\hat{H}$ with kernel $\langle C \rangle$, and that the second
equality defines a faithful action of $SL(2, \mathbb{Z})$.

Next we show that these two actions jointly define an action of the group $L$%
. To this end it is useful to rewrite (\ref{gammaonDrhoaction}) more
functorially in terms of the cocycle $j(\gamma ,\tau )=c\tau +d$, which
satisfies 
\begin{equation}
j(\gamma _{1}\gamma _{2},\tau )=j(\gamma _{1},\gamma _{2}\tau )j(\gamma
_{2},\tau ),\ \ \gamma _{1},\gamma _{2}\in \Gamma _{1}.
\label{cocycleidentity}
\end{equation}%
Thus 
\begin{equation*}
\gamma .(\tau ,w,\rho )=\left( \gamma \tau ,\frac{w}{j(\gamma ,\tau )},\frac{%
\rho }{j(\gamma ,\tau )^{2}}\right) ,
\end{equation*}%
and we have to show that 
\begin{equation}
\gamma ^{-1}\mu (x,y,z)\gamma .(\tau ,w,\rho )=\mu ((x,y)\gamma ,z).(\tau
,w,\rho ).  \label{relationpreserved}
\end{equation}%
The right-hand-side of (\ref{relationpreserved}) is equal to 
\begin{equation*}
(\tau ,w+2\pi i((ax+cy)\tau +bx+dy),\rho ).
\end{equation*}%
The left-hand-side is equal to 
\begin{eqnarray*}
&&\gamma ^{-1}\mu (x,y,z).\left( \gamma \tau ,\frac{w}{j(\gamma ,\tau )},%
\frac{\rho }{j(\gamma ,\tau )^{2}}\right)  \\
&=&\gamma ^{-1}.\left( \gamma \tau ,\frac{w}{j(\gamma ,\tau )}+2\pi
i(x\gamma \tau +y),\frac{\rho }{j(\gamma ,\tau )^{2}}\right)  \\
&=&\gamma ^{-1}.\left( \gamma \tau ,\frac{w+2\pi i(x(a\tau +b)+y(c\tau +d))}{%
j(\gamma ,\tau )},\frac{\rho }{j(\gamma ,\tau )^{2}}\right)  \\
&=&\left( \tau ,\frac{w+2\pi i(x(a\tau +b)+y(c\tau +d))}{j(\gamma ,\tau
)j(\gamma ^{-1},\gamma \tau )},\frac{\rho }{(j(\gamma ,\tau )j(\gamma
^{-1},\gamma \tau ))^{2}}\right)  \\
&=&\left( \tau ,w+2\pi i(x(a\tau +b)+y(c\tau +d)),\rho \right) ,
\end{eqnarray*}%
where we used (\ref{cocycleidentity}) to get the last equality. This
confirms (\ref{relationpreserved}).

It remains to show that the action of $L$ preserves $\mathcal{D}^{\rho }$,
and for this it is enough to prove it for a set of generators. Bearing in
mind the definition of $\mathcal{D}^{\rho }$ (\ref{Drho}), the result is
clear for $\mu (x,y,z)$. To prove it for $\gamma \in \Gamma _{1}$, we must
show that if $(\tau ,w,\rho )\in \mathcal{D}^{\rho }$ then 
\begin{equation}
\left\vert \frac{w}{j(\gamma ,\tau )}-\lambda \right\vert >2\left\vert \frac{%
\rho }{(j(\gamma ,\tau ))^{2}}\right\vert ^{1/2}>0  \label{gammaactionineq}
\end{equation}%
for all $\lambda \in \Lambda _{\gamma \tau }$. But $\lambda =\frac{1}{%
j(\gamma ,\tau )}\lambda ^{\prime }$ for some $\lambda ^{\prime }\in \Lambda
_{\tau }$, whence (\ref{gammaactionineq}) reduces to $|j(\gamma ,\tau
)|^{-1}|w-\lambda ^{\prime }|>2|j(\gamma ,\tau )|^{-1}|\rho |^{1/2}>0$. This
follows from the fact that $(\tau ,w,\rho )\in \mathcal{D}^{\rho }$, and the
proof of the Lemma is complete. \ \ \ \ \ \ $\qed$

\subsubsection{Some covering spaces}

One sees that projection onto the first coordinate 
\begin{eqnarray*}
pr_{1} &:&\mathcal{D}^{\rho }\rightarrow \mathbb{H}_{1}, \\
&&(\tau ,w,\rho )\mapsto \tau ,
\end{eqnarray*}%
is \emph{locally trivial} with contractible base $\mathbb{H}_{1}$. From the
long exact sequence associated to a fibration we obtain an exact sequence $%
0=\pi _{2}(\mathbb{H}_{1})\rightarrow \pi _{1}(F)\rightarrow \pi _{1}(%
\mathcal{D}^{\rho })\rightarrow \pi _{1}(\mathbb{H}_{1})=0$, where $F$ is
the fiber. Thus, we have $\pi _{1}(\mathcal{D}^{\rho })\cong \pi _{1}(F).$
From Lemma 6.5, there is a free action of $\mathbb{Z}^{2}=\hat{H}/Z(L)$ on
each fiber $pr_{1}^{-1}(\tau ).$ Furthermore, from the definition of $%
\mathcal{D}^{\rho }$ we see that 
\begin{eqnarray*}
\pi _{1}(\mathcal{D}^{\rho }/\mathbb{Z}^{2}) &\cong &\pi _{1}(\mathbb{C}%
/\Lambda _{\tau }\setminus \{0\})\times \pi _{1}(\mathbb{C}\setminus \{0\})
\\
&\cong &H\times \mathbb{Z}.
\end{eqnarray*}%
Here, $H$ is the Heisenberg group of the previous subsection.

We need to describe this identification carefully. Consider the usual
realization of $\mathbb{C}/\Lambda _{\tau }$ as the fundamental
parallelogram for $\Lambda _{\tau }$ with identification of sides, and let $%
\alpha ,\beta $ be the cycles along the sides with periods $2\pi i\tau ,2\pi
i$ respectively. Define $\delta $ to be a closed \emph{clockwise} contour
about an interior point of the parallelogram with local coordinate $w=0$.
Then there is an isomorphism of groups 
\begin{eqnarray*}
&&\pi _{1}(\mathbb{C}/\Lambda _{\tau }\setminus \{0\})\overset{\cong }{%
\rightarrow }H, \\
&&\alpha \mapsto A,\beta \mapsto B,\delta \mapsto C^{\prime }.
\end{eqnarray*}%
Similarly, let $\eta $ denote a closed \emph{anti-clockwise} contour about $%
\rho =0$ in the complex plane. Then $\pi _{1}(\mathbb{C}\setminus
\{0\})=\langle \eta \rangle .$

Let $\tilde{\mathcal{D}}^{\rho }$ be a universal covering space of $\mathcal{%
D}^{\rho }$ with covering projection 
\begin{equation*}
p_{1}:\tilde{\mathcal{D}}^{\rho }\rightarrow \mathcal{D}^{\rho }.
\end{equation*}%
There is a free action of the fundamental group $H\times \mathbb{Z}$ on $%
\tilde{\mathcal{D}}^{\rho }$, and we define 
\begin{equation*}
\hat{\mathcal{D}}^{\rho }=\tilde{\mathcal{D}}^{\rho }/\langle \eta
^{-2}\delta \rangle .
\end{equation*}%
Thus we have a sequence of covering projections 
\begin{equation}
\tilde{\mathcal{D}}^{\rho }\overset{p_{3}}{\longrightarrow }\hat{\mathcal{D}}%
^{\rho }\overset{p_{4}}{\longrightarrow }\mathcal{D}^{\rho }\overset{p_{2}}{%
\longrightarrow }\mathcal{D}^{\rho }/\mathbb{Z}^{2},  \label{seqcovproj}
\end{equation}%
where $p_{1}=p_{4}\circ p_{3}$. The action of $\Gamma _{1}$ on $\mathcal{D}%
^{\rho }$ lifts (modulo the fundamental group) to an action on the universal
cover. That is, there is a group $G$ acting on $\tilde{\mathcal{D}}^{\rho }$
where $G$ fits into a short exact sequence 
\begin{equation*}
1\rightarrow H\times \mathbb{Z}\rightarrow G\rightarrow \Gamma
_{1}\rightarrow 1.
\end{equation*}%
We have $Z(G)=Z(H)\times \mathbb{Z}$, in particular $\eta ^{-2}\delta \in
Z(G)$. It follows that $G$ acts on $\hat{\mathcal{D}}^{\rho }$, and there is
a sequence of surjective group maps 
\begin{equation}
G\rightarrow G/\langle \eta ^{-2}\delta \rangle \rightarrow L\rightarrow
\Gamma _{1}  \label{sesG_def}
\end{equation}%
in which the four groups act on the corresponding spaces in (\ref{seqcovproj}%
).

\subsubsection{ Lifting the logarithm $l(x)$}

From (\ref{Om22rho}), the logarithmic contribution to $\Omega _{22}^{(2)}$
is 
\begin{equation}  \label{ldef}
l(x)=\log (-\frac{\rho }{K(\tau ,w)^{2}}), \ x = (\tau, w, \rho) \in 
\mathcal{D}^{\rho}.
\end{equation}%
The remaining parts $\omega_{11}$, $\omega_{\beta 1}$ and $\omega_{\beta 
\bar{\beta}}$ of $\Omega ^{(2)}$ are single-valued on $\mathcal{D}^{\rho }$
since they are expressible in terms of the Weierstrass functions and
Eisenstein series. Now $K(\tau ,w)^{2}=-\theta _{1}(\tau ,w)^{2}/\eta (\tau
)^{6}$ is a Jacobi form of weight $-2$ and index $1$ \cite{EZ}, so that $%
\exp l(x)=\frac{- \rho }{K(\tau ,w)^{2}}$ is single-valued. The way it
transforms under the Jacobi group $J$ can be read-off of Lemma \ref%
{lemmaLaction}. We find that 
\begin{eqnarray}
\exp l((a,b).x) &=&\exp (2\pi a^{2}i\tau +2aw)\exp l(x),\quad (a,b)\in 
\mathbb{Z}^{2}  \label{ab_rho_expl} \\
\exp l(\gamma_{1}.x) &=&\exp (-\frac{1}{2\pi i}\frac{c_{1}w^{2}}{c_{1}\tau
+d_{1}})\exp l(x),\quad \gamma_{1}\in \Gamma_{1},  \label{g1_rho_expl}
\end{eqnarray}%
where $(a, b)$ is the image of $\mu(a, b, c)$ in $J$. For a given choice of
the branch $l(x)$, we therefore find that 
\begin{equation}
l((a,b).x)=l(x)+2\pi a^{2}i\tau +2aw+2\pi iN(a,b),\quad (a,b)\in \mathbb{Z}%
^{2},  \notag
\end{equation}%
for some $N(a,b) \in \mathbb{Z}.$

Let $\tilde{l}(\tilde{x})$ be a lifting of $l(x)$ to a single-valued
function on $\mathcal{\tilde{D}}^{\rho }$. $K(\tau ,w)^{2}$ is holomorphic
for $(\tau ,w)\in \mathbb{H}_{1}\mathbb{\times C}$ with a zero of order two
for each $w\in \Lambda _{\tau }$ (see (\ref{Kthetaeta})). Let $\tilde{x}\in 
\tilde{\mathcal{D}}^{\rho }$ and $p_{1}(\tilde{x})=x=(\tau ,w,\rho ).$ Using
(\ref{ab_rho_expl}) we find: 
\begin{eqnarray*}
\tilde{l}(\alpha .\tilde{x}) &=&\tilde{l}(\tilde{x})+2\pi i\tau +2w+2\pi
iN_{\alpha }, \\
\tilde{l}(\beta .\tilde{x}) &=&\tilde{l}(\tilde{x})+2\pi iN_{\beta }, \\
\tilde{l}(\eta .\tilde{x}) &=&\tilde{l}(\tilde{x})+2\pi i, \\
\tilde{l}(\delta .\tilde{x}) &=&\tilde{l}(\tilde{x})+4\pi i,
\end{eqnarray*}%
for some $N_{\alpha },N_{\beta }\in \mathbb{Z}$. In particular, note that by
composing these transformations we confirm the relation $[\alpha ,\beta
]=\delta $. We may define new generators $\alpha ^{\prime }=\alpha \eta
^{-N_{\alpha }}$, $\beta ^{\prime }=\beta \eta ^{-N_{\beta }}$ which satisfy
the same relations and for which $N_{\alpha ^{\prime }}=N_{\beta ^{\prime
}}=0$. Relabelling, we then obtain

\begin{lemma}
\label{Lemma_Pi_action} With previous notation, we have 
\begin{equation}
\tilde{l}(\alpha ^{a}\beta ^{b}\gamma ^{c}\delta ^{d}.\tilde{x})=\tilde{l}(%
\tilde{x})+2\pi ia^{2}\tau +2aw+2\pi i(c+2(ab+d)),  \label{ltilde_alpha}
\end{equation}%
for $a,b,c,d\in \mathbb{Z}$. In particular, 
\begin{eqnarray*}
\tilde{l}(\eta^{-2}\delta.\tilde{x}) = \tilde{l}(\tilde{x}),
\end{eqnarray*}
so that $\tilde{l}$ pushes down to a single-valued function $\hat{l}$ on $%
\hat{\mathcal{D}}^{\rho}.$ \ \ \ \ $\qed$
\end{lemma}

From (\ref{g1_rho_expl}) we find for $\gamma \in \Gamma _{1}$ that 
\begin{equation}
l(\gamma _{1}.x)=l(x)-\frac{1}{2\pi i}.\frac{c_{1}w^{2}}{c_{1}\tau +d_{1}}%
+2\pi iN(\gamma _{1}),\quad  \label{gam1_cocycle}
\end{equation}%
for some $N(\gamma _{1})\in \mathbb{Z}$. It is easy to see that (\ref%
{gam1_cocycle}) is consistent with respect to the composition of $\gamma
_{1},\gamma _{2}\in \Gamma _{1}$ with a trivial cocycle condition $N(\gamma
_{1}\gamma _{2})=N(\gamma _{1})+N(\gamma _{2})$. Thus the extension (\ref%
{sesG_def}) splits, in particular $G$ contains the subgroup $L=\hat{H}\Gamma
_{1}$ and $G=L\times \mathbb{Z}.$ Since $L\cap \langle \eta ^{-2}\delta
\rangle =1$ there is an injection 
\begin{equation*}
L\longrightarrow G/\langle \eta ^{-2}\delta \rangle ,
\end{equation*}%
and through this map $L$ acts on $\mathcal{\hat{D}}^{\rho }.$ We can now
read-off from (\ref{ltilde_alpha}) and (\ref{gam1_cocycle}) that $\hat{l}$
transforms as followings:

\begin{theorem}
\label{Theorem_l_L} The action of $L$ on $\mathcal{\hat{D}}^{\rho }$
satisfies 
\begin{eqnarray}
\hat{l}(\mu (a,b,c).\hat{x}) &=&\hat{l}(\hat{x})+2\pi ia^{2}\tau +2aw+2\pi
i(ab+c),\quad \mu (a,b,c)\in \hat{H},  \notag \\
&&  \label{l_mu} \\
\hat{l}(\gamma _{1}.\hat{x}) &=&\hat{l}(\hat{x})-\frac{1}{2\pi i}.\frac{%
c_{1}w^{2}}{c_{1}\tau +d_{1}},\quad \gamma _{1}\in \Gamma _{1},
\label{l_gamma1}
\end{eqnarray}%
where $\hat{x}\in \hat{\mathcal{D}}^{\rho },p_{4}(\hat{x})=(\tau ,w,\rho )$.
\ \ \ \ \ $\qed$
\end{theorem}

\subsubsection{Equivariance of $\hat{F}^{\protect\rho }$ and $F^{\protect%
\rho }$}

After the results of the previous subsection we know that $F^{\rho }$ lifts
to a single-valued holomorphic function $\hat{F}^{\rho }$ on $\mathcal{\hat{D%
}}^{\rho }$: 
\begin{eqnarray}
\hat{F}^{\rho }:\mathcal{\hat{D}}^{\rho } &\rightarrow &\mathbb{H}_{2}, 
\notag \\
\hat{x} &\mapsto &\hat{\Omega}^{(2)}(\hat{x}).  \label{Frhohat}
\end{eqnarray}%
By Proposition \ref{Prop_rhoperiod_graph} we have 
\begin{eqnarray*}
2\pi i\hat{\Omega}_{11}^{(2)}(\hat{x}) &=&2\pi i\tau -\rho \omega _{11}(x),
\\
2\pi i\hat{\Omega}_{12}^{(2)}(\hat{x}) &=&w-\rho ^{1/2}\omega _{\beta 1}(x),
\\
2\pi i\hat{\Omega}_{22}^{(2)}(\hat{x}) &=&\hat{l}(\hat{x})-\omega _{\beta 
\bar{\beta}}(x),
\end{eqnarray*}%
where $(\tau ,w,\rho )=x=p_{4}(\hat{x})$.

\begin{theorem}
\label{Theorem_rho_equiv} $\hat{F}^{\rho }$ is equivariant with respect to
the action of $L$. Thus, there is a commutative diagram for $\gamma \in L$ 
\begin{equation*}
\begin{array}{ccc}
\mathcal{\hat{D}}^{\rho } & \overset{\hat{F}^{\rho }}{\rightarrow } & 
\mathbb{H}_{2} \\ 
\gamma \downarrow &  & \downarrow \gamma \\ 
\mathcal{\hat{D}}^{\rho } & \overset{\hat{F}^{\rho }}{\rightarrow } & 
\mathbb{H}_{2}%
\end{array}%
\end{equation*}
\end{theorem}

\noindent \textbf{Proof.} \ It suffices to consider the separate actions of $%
\hat{H}$ and $\Gamma _{1}$, and we first consider that of $\hat{H}$. From (%
\ref{eq: modtrans}), elements of $\hat{H}$ act as follows: 
\begin{equation}
\mu (a,b,c):\hat{\Omega}\mapsto \left( 
\begin{array}{cc}
\hat{\Omega}_{11}, & \hat{\Omega}_{12}+a\hat{\Omega}_{11}+b \\ 
\hat{\Omega}_{12}+a\hat{\Omega}_{11}+b, & \hat{\Omega}_{22}+a^{2}\hat{\Omega}%
_{11}+2a\hat{\Omega}_{12}+ab+c%
\end{array}%
\right) .  \label{Om_mu}
\end{equation}%
We must show that the matrix in (\ref{Om_mu}) coincides with $\hat{F}^{\rho
}(\mu (a,b,c).\hat{x})$. Set $x=(\tau ,w,\rho )$. Using (\ref{muaction}),
the periodicity of $P_{k}(\tau ,w)$ in $w$ for $k>1$, and the
quasi-periodicity of $P_{1}(\tau ,w)$ (\ref{P1quasiperiod}), we find that $%
R(k,l)$ and $\beta (k)$ satisfy 
\begin{eqnarray}
R(k,l)((a,b).x) &=&R(k,l)(x),  \label{R_mu} \\
\beta (k)((a,b).x) &=&\beta (k)(x)+a\rho ^{1/2}\delta _{k,1}.
\label{beta_mu}
\end{eqnarray}%
Thus $\omega _{11}$, $\omega _{\beta 1}$ and $\omega _{\beta \bar{\beta}}$
satisfy%
\begin{eqnarray*}
\omega _{11}((a,b).x) &=&\omega _{11}(x), \\
\omega _{\beta 1}((a,b).x) &=&\omega _{\beta 1}(x)+a\rho ^{1/2}\omega
_{11}(x), \\
\omega _{\beta \bar{\beta}}((a,b).x) &=&\omega _{\beta \bar{\beta}%
}(x)+a^{2}\rho \omega _{11}(x)+2a\rho ^{1/2}\omega _{\beta 1}(x).
\end{eqnarray*}%
We therefore find 
\begin{equation*}
\hat{\Omega}_{11}(\mu (a,b,c).\hat{x})=\tau -\frac{\rho }{2\pi i}\omega
_{11}((a,b).x)=\hat{\Omega}_{11}(\hat{x}).
\end{equation*}%
Similarly, we have 
\begin{eqnarray*}
\hat{\Omega}_{12}(\mu (a,b,c).\hat{x}) &=&\frac{1}{2\pi i}(w+2\pi ia\tau
+2\pi ib-\rho ^{1/2}\omega _{\beta 1}((a,b).x)) \\
&=&\hat{\Omega}_{12}(\hat{x})+a\hat{\Omega}_{11}(\hat{x})+b.
\end{eqnarray*}%
Now application of (\ref{l_mu}) yields 
\begin{eqnarray*}
\hat{\Omega}_{22}(\mu (a,b,c).\hat{x}) &=&\frac{1}{2\pi i}(\hat{l}(\mu
(a,b,c).\hat{x})-\omega _{\beta \bar{\beta}}((a,b).x)) \\
&=&\frac{1}{2\pi i}(\hat{l}(\hat{x})+2\pi ia^{2}\tau +2aw+2\pi i(ab+c) \\
&&-\omega _{\beta \bar{\beta}}(x)-a^{2}\rho \omega _{11}(x)-2a\rho
^{1/2}\omega _{\beta 1}(x)) \\
&=&\hat{\Omega}_{22}+a^{2}\hat{\Omega}_{11}+2a\hat{\Omega}_{12}+ab+c.
\end{eqnarray*}%
This establishes equivariance of $\hat{F}^{\rho }$ with respect to $\hat{H}$.

As in the $\epsilon $-formalism, the exceptional transformation law (\ref%
{gammaE2}) for $E_{2}$ plays a critical r\^{o}le in establishing $\Gamma
_{1} $-equivariance of $\hat{F}^{\rho }$. Consider the action (\ref%
{gamma1Omega}) of $\Gamma _{1}$ on $\mathbb{H}_{2}$. Since $E_{2}$ appears
only in $R(1,1)$ and $\beta (1)$, (\ref{gammaonDrhoaction}) implies that 
\begin{eqnarray}
R(k,l)(\gamma _{1}.x) &=&R(k,l)(x)+\kappa \delta _{k,1}\delta _{l,1},
\label{Rgam1} \\
\beta (k)(\gamma _{1}.x) &=&\beta (k)(x)-\kappa \frac{w}{\rho ^{1/2}}\delta
_{k,1},  \label{betagam1} \\
\kappa &=&\frac{c_{1}}{c_{1}\tau +d_{1}}\frac{\rho }{2\pi i}.
\label{anotherkappadef}
\end{eqnarray}%
We then have 
\begin{equation*}
\hat{\Omega}_{11}(\gamma _{1}.\hat{x})=\frac{1}{c_{1}\tau +d_{1}}(a_{1}\tau
+b_{1}-\frac{1}{2\pi i}\frac{\rho }{c_{1}\tau +d_{1}}\omega _{11}(\gamma
_{1}.x)).
\end{equation*}%
Similarly to Theorem \ref{TheoremGequiv} in the $\epsilon $ formalism, (\ref%
{Rgam1}) implies that the transformation under $\gamma _{1}$ of the weight $%
\omega (N)$ for $N\in $ $\mathcal{N}_{1,a;1,b}\ $is the sum of the weights
of the product over all possible necklaces in $\mathcal{N}_{1,a;1,b}$ formed
from $N$ by deleting the edges of type $\overset{1,a_{1}}{\bullet }\overset{}%
{\longrightarrow }\overset{1,a_{2}}{\bullet }$ and multiplying by a $\kappa $
factor for each such deletion. From Proposition \ref{Prop_rhoperiod_graph}
and (\ref{anotherkappadef}) we obtain 
\begin{equation*}
1-\kappa \omega _{11}=(c_{1}\hat{\Omega}_{11}+d_{1})/(c_{1}\tau +d_{1}).
\end{equation*}%
Then we find, much as before, that 
\begin{eqnarray}
\omega _{11}(\gamma _{1}.x) &=&\sum_{n\geq 0}\kappa ^{n}\omega _{11}^{n+1}(x)
\notag \\
&=&\frac{(c_{1}\tau +d_{1})\omega _{11}(x)}{c_{1}\hat{\Omega}_{11}+d_{1}}.
\label{omega11gamma1}
\end{eqnarray}%
Then $\hat{\Omega}_{11}(\gamma _{1}.\hat{x})$ is as given in (\ref%
{gamma1Omega}).

We next have 
\begin{equation*}
\hat{\Omega}_{12}(\gamma_{1}.\hat{x})=\frac{1}{c_{1}\tau +d_{1}}\frac{1}{%
2\pi i}(w-\rho ^{1/2}\omega_{\beta 1}(\gamma_{1}.x)).
\end{equation*}%
(\ref{Rgam1}) and (\ref{betagam1}) imply that $\sum_{k}\beta_{a}(k)\omega
(N) $ for $N\in \mathcal{N}_{k,a;1,b}$ transforms under $\gamma_{1}$ as a
sum over the product with $\kappa $ factors of weights of necklaces in $%
\mathcal{N}_{1,a;1,b}$ and at most one necklace in $\mathcal{N}_{k,a;1,b}$
with a $\beta_{a}(k)$ factor. Then one finds 
\begin{eqnarray*}
\rho ^{1/2}\omega_{\beta 1}(\gamma_{1}.x) &=&\frac{\rho ^{1/2}\omega _{\beta
1}(x)-\kappa w\omega_{11}(x)}{1-\kappa \omega_{11}(x)} \\
&=&w-2\pi i(c_{1}\tau +d_{1}).\frac{\hat{\Omega}_{12}}{c_{1}\hat{\Omega}%
_{11}+d_{1}}.
\end{eqnarray*}%
This implies $\hat{\Omega}_{12}(\gamma_{1}.\hat{x})$ is as given in (\ref%
{gamma1Omega}).

Finally, using (\ref{l_gamma1}) we have 
\begin{equation*}
\hat{\Omega}_{22}(\gamma_{1}.\hat{x})=\frac{1}{2\pi i}(\hat{l}(\hat{x}%
)+\kappa \frac{w^{2}}{\rho }-\omega_{\beta \bar{\beta}}(\gamma_{1}.x)).
\end{equation*}%
Using a similar argument, (\ref{Rgam1}) and (\ref{betagam1}) imply that 
\begin{eqnarray*}
\omega_{\beta \bar{\beta}}(\gamma_{1}.x) &=&\omega_{\beta \bar{\beta}}(x)+%
\frac{\kappa \omega_{\beta 1}^{2}(x)-2\kappa \frac{w}{\rho ^{1/2}}\omega
_{\beta 1}(x)+\kappa ^{2}\frac{w^{2}}{\rho }\omega_{11}(x)}{1-\kappa \omega
_{11}} \\
&=&\omega_{\beta \bar{\beta}}(x)-\kappa \frac{w^{2}}{\rho }+2\pi i\frac{c_{1}%
\hat{\Omega}_{12}^{2}}{c_{1}\hat{\Omega}_{11}+d_{1}}.
\end{eqnarray*}%
Hence $\hat{\Omega}_{22}(\gamma_{1}.\hat{x})$ is as given in (\ref%
{gamma1Omega}) and hence $\Gamma_{1}$ acts equivariantly. This completes the
proof of the Theorem. \ \ \ \ $\qed$

\begin{remark}
\label{Rem_simple_Dehn_twist_rho} \ Much the same as in Remark \ref%
{Rem_simple_Dehn_twist}, $\Omega_{22}\rightarrow \Omega_{22}+1$ is generated
by $C$ corresponding to a Dehn twist in the connecting cylinder.
\end{remark}

We may also choose a branch of $l(x)$ and consider the equivariance of $%
F^{\rho }$ on $\mathcal{D}^{\rho }$ under the action of the subgroup $\Gamma
_{1}.$

\begin{corollary}
\label{Cor_rho_G1_equiv} For any choice of branch for $l(x),F^{\rho }$ is
equivariant with respect to the action of $\Gamma _{1}$. Thus, there is a
commutative diagram for $\gamma \in \Gamma _{1},$ 
\begin{equation*}
\begin{array}{ccc}
\mathcal{D}^{\rho } & \overset{F^{\rho }}{\rightarrow } & \mathbb{H}_{2} \\ 
\gamma \downarrow &  & \downarrow \gamma \\ 
\mathcal{D}^{\rho } & \overset{F^{\rho }}{\rightarrow } & \mathbb{H}_{2}%
\end{array}%
\ \ \ \ \ \ \qed
\end{equation*}
\end{corollary}

\subsection{Local Invertibility About the Two Tori Degeneration Limiting
Point}

\label{Subsec_Frhoinvert}We now consider the invertibility of the $F^{\rho }$
map in the neighborhood of a degeneration point. In the $\rho $-formalism
this degeneration\ is more subtle than that of the $\epsilon $-formalism
discussed in Subsection \ref{Subsec_Feps_degen}. Define the $\Gamma _{1}$%
-invariant parameter 
\begin{equation}
\chi =-\frac{\rho }{w^{2}}.  \label{chidef}
\end{equation}%
We will show that $\rho ,w\rightarrow 0$ $\,$for fixed $\chi $ is the $2$%
-torus degeneration limit. From (\ref{Drho}) we have $|w|>2|\rho |^{1/2}$
(for $\lambda =0$) so that $|\chi |<\frac{1}{4}$. (Similarly to Subsection %
\ref{Subsubsect_Catalan}, $\chi =\frac{1}{4}$ is a singular point where two
degenerate annuli touch at $z_{1}=-z_{2}=w/2.$) To describe this limit more
precisely, we introduce the domain 
\begin{equation}
\mathcal{D}^{\chi }=\{(\tau ,w,\chi )\in \mathbb{H}_{1}\times \mathbb{%
C\times C}\ |\ (\tau ,w,-w^{2}\chi )\in \mathcal{D}^{\rho },0<|\chi |<\frac{1%
}{4}\}.  \label{Dchidomain}
\end{equation}%
Thanks to Theorem \ref{Theorem_rhoperiod} and Corollary \ref%
{Cor_rho_G1_equiv}, there is a $\Gamma _{1}$-equivariant holomorphic map 
\begin{eqnarray}
F^{\chi }:\mathcal{D}^{\chi } &\rightarrow &\mathbb{H}_{2},  \notag \\
(\tau ,w,\chi ) &\mapsto &\Omega ^{(2)}(\tau ,w,-w^{2}\chi ).  \label{Fchi}
\end{eqnarray}%
Let 
\begin{equation*}
\mathcal{D}_{0}^{\chi }=\{(\tau ,0,\chi )\in \mathbb{H}_{1}\times \mathbb{%
C\times C}|0<|\chi |<\frac{1}{4}\}
\end{equation*}%
denote the space\ of limit points where $\rho ,w\rightarrow 0$ $\,$for fixed 
$\chi \neq 0$. Then

\begin{proposition}
\label{PropOmrhodegen}For $(\tau ,w,\chi )\in \mathcal{D}^{\chi }\cup 
\mathcal{D}_{0}^{\chi }$ we have 
\begin{eqnarray}
2\pi i\Omega_{11}^{(2)} &=&2\pi i\tau +w^{2}(1-4\chi )G(\chi )+O(w^{4}), 
\notag \\
2\pi i\Omega_{12}^{(2)} &=&w\sqrt{1-4\chi }(1+w^{2}(1-4\chi )E_{2}(\tau
)G(\chi )+O(w^{4})),  \notag \\
2\pi i\Omega_{22}^{(2)} &=&\log f(\chi )+w^{2}(1-4\chi )E_{2}(\tau
)+O(w^{4}),  \label{Omrhodegen}
\end{eqnarray}%
where 
\begin{equation*}
G(\chi )=\frac{1}{12}+E_{2}(q=f(\chi )),
\end{equation*}%
and $f(\chi )$ is the Catalan series (\ref{Catalan}).
\end{proposition}

\noindent \textbf{Proof.} Note that $P_{n}(\tau ,w)=\frac{1}{w^{n}}%
(1+w^{2}E_{2}(\tau )(\delta _{n,2}-\delta _{n,1})+O(w^{4}))$ from (\ref{P2})
and (\ref{P1}). Then (\ref{Rg2}) and (\ref{betag2}) imply 
\begin{eqnarray*}
R(k,l) &=&R^{(0)}(k,l)+w^{2}\chi E_{2}(\tau )\delta _{k,1}\delta
_{l,1}+O(w^{4}), \\
\beta (k) &=&\beta ^{(0)}(k)(1-w^{2}E_{2}(\tau )\delta _{k,1})+O(w^{4}), \\
\log (-\frac{\rho }{K(\tau ,w)^{2}})) &=&\log \chi +w^{2}E_{2}(\tau
)+O(w^{4})\text{.}
\end{eqnarray*}%
using (\ref{Genuszerodata}). For $w=0$ these expressions are exactly those
found in Proposition \ref{Prop_torus2} for a torus formed from a sphere by
sewing an annulus centered at $z=0$ to another centered at $z=w$. Expanding (%
\ref{Om11rho}) with $\rho =-\chi w^{2}$ to order $w^{2}$ implies 
\begin{equation*}
2\pi i\Omega _{11}^{(2)}=2\pi i\tau +w^{2}\chi \sigma
((I-R^{(0)})^{-1}(1,1))+O(w^{4}).
\end{equation*}%
But (\ref{E2chi}) implies 
\begin{equation}
\sigma ((I-R^{(0)})^{-1}(1,1))=2(I+B^{(0)})^{-1}(1,1)=\frac{(1-4\chi )}{\chi 
}G(\chi ),  \label{Gchi}
\end{equation}%
leading to the stated result for $\Omega _{11}^{(2)}$. Similarly, for $%
\Omega _{12}^{(2)}$ we find 
\begin{eqnarray*}
2\pi i\Omega _{12}^{(2)} &=&w[1-(-\chi )^{1/2}\sigma (\beta
^{(0)}(1-R^{(0)})^{-1}(1))]. \\
&&[1+w^{2}E_{2}(\tau )\chi \sigma ((1-R^{(0)})^{-1}(1,1))+O(w^{4})].
\end{eqnarray*}%
After some algebra and using (\ref{chi(f)}), (\ref{Snk}) and (\ref{Catalan^k}%
) one finds 
\begin{eqnarray*}
1-(-\chi )^{1/2}\sigma (\beta ^{(0)}(1-R^{(0)})^{-1}(1)) &=&1-2\chi
\sum\limits_{n\geq 1}S_{n,1}(\chi ) \\
&=&1-2\chi (1+f(\chi )) \\
&=&\sqrt{1-4\chi }.
\end{eqnarray*}%
The stated result for $\Omega _{12}^{(2)}$ then follows on using (\ref{Gchi}%
) again. Finally, for $\Omega _{22}^{(2)}$ we find as above, using
Proposition \ref{Prop_torus2}, that 
\begin{eqnarray*}
2\pi i\Omega _{22}^{(2)} &=&\log \chi -\beta ^{(0)}(1-R^{(0)})^{-1}\bar{\beta%
}^{(0)T} \\
&&+E_{2}(\tau )w^{2}[1-(-\chi )^{1/2}\sigma (\beta
^{(0)}(1-R^{(0)})^{-1}(1))]^{2}+O(w^{4}) \\
&=&\log f(\chi )+w^{2}(1-4\chi )E_{2}(\tau )+O(w^{4}).\ \ \ \qed
\end{eqnarray*}

The restriction of $F^{\chi }$ to $\mathcal{D}_{0}^{\chi }$ induces the
natural identification 
\begin{eqnarray}
F^{\chi }:\mathcal{D}_{0}^{\chi } &\overset{\sim }{\rightarrow }&\mathbb{H}%
_{1}\times \mathbb{H}_{1}\subseteq \mathbb{H}_{2}  \notag \\
(\tau ,0,\chi ) &\mapsto &\left( 
\begin{array}{cc}
\tau & 0 \\ 
0 & \frac{1}{2\pi i}\log f(\chi )%
\end{array}%
\right) ,  \label{rhodegen}
\end{eqnarray}%
i.e. $\Omega ^{(2)}=\mathrm{diag}(\Omega_{11}^{(2)},\Omega_{22}^{(2)})$. The
invertibility of the $\Gamma_{1}$-equivariant map $F^{\chi }$ in a
neighborhood of a point in $\mathcal{D}_{0}^{\chi }$ then follows:

\begin{proposition}
\label{Prop_Frhoinverse} Let $x\in \mathcal{D}_{0}^{\chi }$. Then there
exists a $\Gamma_{1}-$invariant neighborhood $\mathcal{N}_{x}^{\chi
}\subseteq \mathcal{D}^{\chi }$ of $x$ throughout which $F^{\chi }$ is
invertible.
\end{proposition}

\noindent \textbf{Proof.} The proof is very similar to Proposition \ref%
{Prop_Fepsinverse}. Let $x=(\tau ,0,\chi )$. The Jacobian of the $F^{\chi }$
map at $x$\ is from (\ref{Omrhodegen}) 
\begin{eqnarray*}
\left\vert \frac{\partial (\Omega _{11},\Omega _{12},\Omega _{22})}{\partial
(\tau ,w,\chi )}\right\vert _{x} &=&\left\vert 
\begin{array}{ccc}
1 & 0 & 0 \\ 
0 & \frac{1}{2\pi i}\sqrt{1-4\chi } & 0 \\ 
0 & 0 & \frac{1}{2\pi i}\frac{f^{\prime }(\chi )}{f(\chi )}%
\end{array}%
\right\vert \\
&=&\frac{1}{4\pi ^{2}\chi }\neq 0,
\end{eqnarray*}%
using $f^{\prime }(\chi )=f(\chi )/(\chi \sqrt{1-4\chi })$. By the inverse
function theorem, there exists an open neighborhood of $x\in \mathcal{D}%
_{0}^{\chi }$ throughout which $F^{\chi }$ is invertible. The result then
follows by choosing precisely invariant neighborhoods (under the action of $%
\Gamma _{1}$) $U,V$ of $x$, respectively $y=F^{\chi }(x)$ such that the
conditions of part (a) of Lemma \ref{discontinuity} hold. The open
neighborhood 
\begin{equation*}
\mathcal{N}_{x}^{\chi }=\bigcup_{\gamma \in \Gamma _{1}}\gamma U
\end{equation*}%
has the required properties. \ \ \ \ \ $\qed$

\section{Mapping Between the $\protect\epsilon $ and $\protect\rho $
Parameterizations}

\label{Sect_Epsilon_Rho}We have described in the previous sections two
separate parameterizations for the genus two period matrix $\Omega ^{(2)}$
based on either sewing two punctured tori in the $\epsilon $-formalism or
sewing a twice-punctured torus to itself in the $\rho $-formalism. In this
final section we show that there is a 1-1 mapping between suitable $\Gamma
_{1}-$invariant domains in both parameterizations.

\begin{theorem}
\label{Theorem_epsrho_11map} There exists a 1-1 holomorphic mapping between $%
\Gamma_{1}-$invariant open domains $\mathcal{I}^{\chi }\subset \mathcal{D}%
^{\chi }$ and $\mathcal{I}^{\epsilon }\subset \mathcal{D}^{\epsilon }$ where 
$\mathcal{I}^{\chi }$ and $\mathcal{I}^{\epsilon }$ are open neighborhoods
of a $2$-torus degeneration point.
\end{theorem}

\noindent \textbf{Proof.} From Proposition \ref{Prop_Fepsinverse} there
exists a $G$-invariant (and thus $\Gamma _{1}-$invariant) open domain $%
\mathcal{N}^{\epsilon }\subset \mathcal{D}^{\epsilon }$ such that the
holomorphic map $F^{\epsilon }:\mathcal{N}^{\epsilon }\rightarrow
F^{\epsilon }(\mathcal{N}^{\epsilon })$ is invertible with $F^{\epsilon }(%
\mathcal{N}^{\epsilon })$ an open neighborhood of a given $2$-torus
degeneration point $\Omega ^{(2)}=\mathrm{diag}(\Omega _{11}^{(2)},\Omega
_{22}^{(2)})$. Similarly, from Proposition \ref{Prop_Frhoinverse} there
exists a $\Gamma _{1}$-invariant open domain $\mathcal{N}^{\chi }\subset 
\mathcal{D}^{\chi }$ such that the holomorphic map $F^{\chi }:\mathcal{N}%
^{\chi }\rightarrow F^{\chi }(\mathcal{N}^{\chi })$ is invertible with $%
F^{\chi }(\mathcal{N}^{\chi })$ an open neighborhood of $\mathrm{diag}%
(\Omega _{11}^{(2)},\Omega _{22}^{(2)})$. Define a $\Gamma _{1}$-invariant
open neighborhood of $\mathrm{diag}(\Omega _{11}^{(2)},\Omega _{22}^{(2)})$
by $\mathcal{I}^{\Omega }=F^{\epsilon }(\mathcal{N}^{\epsilon })\cap $ $%
F^{\chi }(\mathcal{N}^{\chi }).$ Hence, defining $\Gamma _{1}$-invariant
open domains $\mathcal{I}^{\chi }=(F^{\chi })^{-1}(\mathcal{I}^{\Omega })$
and $\mathcal{I}^{\epsilon }=(F^{\epsilon })^{-1}(\mathcal{I}^{\Omega })$,
we find $(F^{\epsilon })^{-1}\circ F^{\chi }:\mathcal{I}^{\chi }\rightarrow 
\mathcal{I}^{\epsilon }$ is holomorphic and 1-1. \ \ \ \ $\qed$

We conclude by displaying the explicit form of the 1-1 mapping to order $%
w^{3}$, obtained by combining the expansions of (\ref{tau1lead})-(\ref%
{epslead}) and Proposition \ref{PropOmrhodegen}: 
\begin{eqnarray}
\tau _{1}(\tau ,w,\chi ) &=&\tau +\frac{1}{2\pi i}w^{2}(1-4\chi )\frac{1}{12}%
+O(w^{4}),  \notag \\
\tau _{2}(\tau ,w,\chi ) &=&\frac{1}{2\pi i}\log (f(\chi ))+O(w^{4}),  \notag
\\
\epsilon (\tau ,w,\chi ) &=&-w\sqrt{1-4\chi }(1+w^{2}E_{2}(\tau )(1-4\chi
)+O(w^{4})).  \label{rhoepsilon}
\end{eqnarray}%
It is then straightforward to check that these relations are invariant under
the action of $\Gamma _{1}$ to the given order using (\ref{gammaE2}), (\ref%
{gam1eps}) and (\ref{gammaonDrhoaction}).

\textbf{Acknowledgement.} 
The authors wish to thank Harold Widom and
Alexander Zuevsky for useful discussions.

\section{Appendix}

In this appendix we supply the explicit form of the genus two period matrix $%
\Omega ^{(2)}$ of Theorem \ref{Theorem_epsperiod} in the $\epsilon$%
-formalism to $O(\epsilon ^{9})$ and of Theorem \ref{Theorem_rhoperiod} in
the $\rho$-formalism to $O(\rho ^{5})$.

\begin{align*}
2\pi i\Omega_{11}^{(2)}(\tau_{1},\tau_{2},\epsilon )& =2\pi i\Omega
_{22}^{(2)}(\tau_{2},\tau_{1},\epsilon ) \\
& =2\pi i\tau_{1}+F_{{2}}{\epsilon }^{2}+E_{{2}}{F_{{2}}}^{2}{\epsilon }%
^{4}+({E_{{2}}}^{2}{F_{{2}}}^{3}+6\,E_{{4}}F_{{2}}F_{{4}}){\epsilon }^{6} \\
& +({E_{{2}}}^{3}{F_{{2}}}^{4}+12\,E_{{2}}E_{{4}}{F_{{2}}}^{2}F_{{4}}+10\,E_{%
{6}}F_{{2}}F_{{6}}+30\,E_{{6}}{F_{{4}}}^{2}){\epsilon }^{8} \\
& +O(\epsilon ^{10}), \\
& \\
2\pi i\Omega_{12}^{(2)}(\tau_{1},\tau_{2},\epsilon )& =-\epsilon \lbrack
1+E_{{2}}F_{{2}}{\epsilon }^{2}+({E_{{2}}}^{2}{F_{{2}}}^{2}+3\,E_{{4}}F_{{4}%
}){\epsilon }^{4} \\
& +({E_{{2}}}^{3}{F_{{2}}}^{3}+9\,E_{{2}}E_{{4}}F_{{2}}F_{{4}}+5\,E_{{6}}F_{{%
6}}){\epsilon }^{6} \\
& +({E_{{2}}}^{4}{F_{{2}}}^{4}+15\,{E_{{2}}}^{2}E_{{4}}{F_{{2}}}^{2}F_{{4}%
}+5\,E_{{2}}E_{{6}}F_{{2}}F_{{6}}+30\,E_{{2}}E_{{6}}{F_{{4}}}^{2} \\
& +30\,{E_{{4}}}^{2}F_{{2}}F_{{6}}+9\,{E_{{4}}}^{2}{F_{{4}}}^{2}+7\,E_{{8}%
}F_{{8}}){\epsilon }^{8}]+O(\epsilon ^{11}),
\end{align*}%
where for brevity's sake of we have defined $E_{k}=E_{k}(\tau_{1})$ and $%
F_{k}=E_{k}(\tau_{2})$.

Similarly, in the $\rho$-formalism we find that $\Omega ^{(2)}(\tau ,w,\rho)$
to $O({\rho }^{4})$ is as follows: 
\begin{eqnarray*}
2\pi i\Omega_{11}^{(2)} &=&2\pi i\tau -2\,{\rho }+2\left( \,P_{{2}}+E_{{2}%
}\right) {\rho }^{2}-2\,\left( P_{{2}}+E_{{2}}\right) ^{2}{\rho }^{3} \\
&&+2(\,\left( P_{{2}}+E_{{2}}\right) ^{3}+2\,{P_{{3}}}^{2}){\rho }%
^{4}+O(\rho ^{5}), \\
&& \\
2\pi i\Omega_{12}^{(2)} &=&w+2\,P_{{1}}\rho -2\,P_{{1}}\left( P_{{2}}+E_{{2}%
}\right) {\rho }^{2} \\
&&+2[\,P_{{1}}\left( P_{{2}}+E_{{2}}\right) ^{2}+\,P_{{3}}\left( P_{{2}}-E_{{%
2}}\right) ]{\rho }^{3}-2[\,P_{{3}}\left( P_{{4}}+E_{{4}}\right) \\
&&+\,P_{{1}}\left( P_{{2}}+E_{{2}}\right) ^{3}+2\,P_{{1}}{P_{{3}}}^{2}+\,P_{{%
3}}({P_{{2}}}^{2}-{E_{{2}}}^{2})]{\rho }^{4}+O(\rho ^{5}), \\
&& \\
2\pi i\Omega_{22}^{(2)} &=&\log (-\frac{\rho }{K^{2}(w,\tau )})-2\,{P_{{1}}}%
^{2}\rho +[2\,{P_{{1}}}^{2}\left( P_{{2}}+E_{{2}}\right) +\left( P_{{2}}-E_{{%
2}}\right) ^{2}]{\rho }^{2} \\
&&-[2\,{P_{{1}}}^{2}\left( P_{{2}}+E_{{2}}\right) ^{2}+2/3\,{P_{{3}}}%
^{2}+4\,P_{{1}}P_{{3}}\left( P_{{2}}-E_{{2}}\right) ]{\rho }^{3} \\
&&+[1/2\,{P_{{4}}}^{2}+1/2\,{E_{{4}}}^{2}+3\,\left( P_{{4}}-E_{{4}}\right)
\left( P_{{2}}-E_{{2}}\right) ^{2}+2\,{P_{{1}}}^{2}\left( P_{{2}}+E_{{2}%
}\right) ^{3} \\
&&-E_{{4}}P_{{4}}+4\,P_{{3}}P_{{1}}(P_{{1}}P_{{3}}+E_{{4}}+P_{{4}}+{P_{{2}}}%
^{2}-{E_{{2}}}^{2})]{\rho }^{4}+O(\rho ^{5}),
\end{eqnarray*}%
where $E_{k}=E_{k}(\tau )$ and $P_{k}=P_{k}(w,\tau )$.

\end{document}